\documentclass[11pt]{article}
\usepackage[utf8]{inputenc}
\usepackage{amsmath,amssymb,epsfig,bbm}
\usepackage{stmaryrd,mathabx}
\usepackage{comment}
\usepackage{color}
\usepackage[T1]{fontenc}

\usepackage{interval}

\usepackage[textsize=small]{todonotes}
\usepackage{enumitem}
\usepackage{varwidth}
\setlist{nolistsep}
\usepackage{hyperref}

\newtheorem{defi}{Definition}[section]
\newtheorem{prop}[defi]{Proposition}
\newtheorem{theo}[defi]{Theorem}
\newtheorem{theofr}[defi]{Théorème}
\newtheorem{conj}[defi]{Conjecture}
\newtheorem{lemm}[defi]{Lemma}
\newtheorem{lemmfr}[defi]{Lemme}
\newtheorem{coro}[defi]{Corollary}
\newtheorem{rema}[defi]{Remark}
\newtheorem{exem}[defi]{Example}
\newtheorem{exems}[defi]{Examples}

\newcommand{\bdefi}{\begin{defi}}
\newcommand{\edefi}{\end{defi}}
\newcommand{\bprop}{\begin{prop}}
\newcommand{\eprop}{\end{prop}}
\newcommand{\btheo}{\begin{theo}}
\newcommand{\etheo}{\end{theo}}
\newcommand{\btheofr}{\begin{theofr}}
\newcommand{\etheofr}{\end{theofr}}
\newcommand{\blemm}{\begin{lemm}}
\newcommand{\elemm}{\end{lemm}}
\newcommand{\blemmfr}{\begin{lemmfr}}
\newcommand{\elemmfr}{\end{lemmfr}}
\newcommand{\brema}{\begin{rema}}
\newcommand{\erema}{\end{rema}}
\newcommand{\bexer}{\begin{exem}}
\newcommand{\eexer}{\end{exem}}
\newcommand{\bexems}{\begin{exems}}
\newcommand{\eexems}{\end{exems}}
\newcommand{\bconj}{\begin{conj}}
\newcommand{\econj}{\end{conj}}
\newcommand{\bcoro}{\begin{coro}}
\newcommand{\ecoro}{\end{coro}}
\newcommand{\dem}{\noindent{\bf Proof. }}

\usepackage{mathrsfs}
\renewcommand\mathcal{\mathscr}

\newcommand{\F}{{\cal F}}

\renewcommand{\H}{{\cal H}}
\newcommand{\I}{{\cal I}}

\renewcommand{\L}{{\cal L}}

\newcommand{\OOO}{{\cal O}}

\newcommand{\R}{{\cal R}}

\newcommand{\maths}[1]{{\mathbb #1}}  

\newcommand{\CC}{\maths{C}}
\newcommand{\DD}{\maths{D}}

\newcommand{\HH}{\maths{H}}

\newcommand{\NN}{\maths{N}}

\newcommand{\PP}{\maths{P}}
\newcommand{\QQ}{\maths{Q}}
\newcommand{\RR}{\maths{R}}

\newcommand{\ZZ}{\maths{Z}}

\newcommand{\aaa}{{\mathfrak a}}
\newcommand{\bbb}{{\mathfrak b}}
\newcommand{\ccc}{{\mathfrak c}}

\newcommand{\mmm}{{\mathfrak m}}

\newcommand{\ppp}{{\mathfrak p}}

\newcommand{\weakstar}{\overset{*}\rightharpoonup}
\newcommand{\ra}{\rightarrow}
\newcommand{\bs}{\backslash}

\newcommand{\wt}[1]{{\widetilde{#1}}}

\newcommand{\ga}{\gamma}
\newcommand{\Ga}{\Gamma}

\newcommand{\cqfd}{\hfill$\Box$}

\newcommand{\bigO}{\operatorname{O}}

\newcommand{\card}{{\operatorname{Card}}}
\newcommand{\CAT}{\operatorname{CAT}}

\newcommand{\covol}{\operatorname{covol}}

\newcommand{\diam}{{\operatorname{diam}}}

\newcommand{\Haar}{\operatorname{Haar}}

\newcommand{\id}{\operatorname{id}}

\newcommand{\Leb}{\operatorname{Leb}}

\newcommand{\mult}{\operatorname{mult}}
\newcommand{\Nr}{\operatorname{{\tt N}}}
\newcommand{\ols}{\operatorname{OL}}

\newcommand\Perp{\operatorname{Perp}}

\renewcommand{\Re}{{\operatorname{Re}}}

\newcommand{\sg}{\operatorname{sg}}

\newcommand{\Vol}{\operatorname{Vol}}

\newcommand{\hdr}{{\HH}^2_\RR}
\newcommand{\htr}{{\HH}^3_\RR}

\newcommand{\PSL}{\operatorname{PSL}}

\newcommand{\GL}{\operatorname{GL}}

\newcounter{fig}

\usepackage{tcolorbox}

\title{Pair correlations of logarithms of complex lattice points}
\author{Jouni Parkkonen \and Fr\'ed\'eric Paulin} 
\date{\today}

\begin{document}
\bibliographystyle{../alphanum}
\maketitle
\begin{abstract} 
We study the correlations of pairs of complex logarithms of
$\ZZ$-lattice points in the complex line at various scalings, proving
the existence of pair correlation functions. We prove that at the
linear scaling, the pair correlations exhibit level repulsion, as it
sometimes occurs in statistical physics. We prove total loss of mass
phenomena at superlinear scalings, and Poissonian behaviour at
sublinear scalings. The case of Euler weights has applications to the
pair correlation of the lengths of common perpendicular geodesic arcs
from the maximal Margulis cusp neighborhood to itself in the Bianchi
orbifold $\PSL_2(\ZZ[i]) \bs\htr$.
\end{abstract}

\section{Introduction}
\label{sec:intro}

When studying the asymptotic distribution of a sequence of finite
subsets of $\RR$, finer information is sometimes given by the
statistics of the spacing (or gaps) between pairs or $k$-tuples of
elements, seen at an appropriate scaling. These problems often arise
in quantum chaos, including energy level spacings or clusterings, and
in statistical physics, including molecular repulsion or interstitial
distribution.  See for instance \cite{Montgomery73, Berry88, RudSar98,
  BocZah05, MarStr13, LarSto20, HofKal21, ParPau22b}.  This paper may
be seen as a complex version of our paper \cite{ParPau22a} where we
study the pair correlation of logarithms of pairs of natural integers,
though new phenomena occur, including the necessity to take limits of
the underlying spaces, as we now explain.

The general setting for our study may be described as follows. Let $E$
be an abelian locally compact group. Let $\F=(F_N, \;\omega_N)
_{N\in\NN}$ be a sequence of finite subsets $F_N$ of $E$, endowed with
a {\em weight function} $\omega_N: F_N \ra\; ]\,0,+\infty\,[\,$ (or
multiplicity function when its values are positive integers). When
studying the asymptotic distribution of differences of elements of
$F_N$, looking at them at various scalings is often desirable. As
explained by Gromov (see for instance \cite{Gromov93}), scaling some
metric space (unless this space has a nice family of homotheties, as
the Euclidean space $\RR^n$ does) sometimes requires to change the
space, especially at the limit. We thus introduce a sequence
$(E_N)_{N\in\NN}$ of abelian locally compact groups converging for the
pointed Hausdorff-Gromov convergence to an abelian locally compact
group $E_\infty$ (see for instance \cite{Gromov99a}). Let
$\Haar_{E_\infty}$ be a Haar measure on $E_\infty$. Let $\psi:N\mapsto
\psi(N)$ be a {\it scaling function}, that is, for every
$N\in\NN-\{0\}$, let $\psi(N):E\ra E_N$ be any map, typically a
dilating homeomorphism for appropriate distances, that we think of as
``scaling'' the space $E$. Let $\psi': \NN-\{0\}\ra[1,+\infty[$ be an
appropriately chosen function, called a {\em renormalising
  function}. The {\it pair correlation measure of $\F$ at time $N$
with scaling $\psi(N)$} is the measure on $E_N$ with finite support
\begin{equation}\label{eq:defiPCmeasure}
\R^{\F,\psi}_N=\sum_{x,y\in F_N, x\neq y}\;\omega_N(x)\,\omega_N(y)\,
\Delta_{\psi(N)(y-x)}\;,
\end{equation}
where $\Delta_z$ denotes the unit Dirac mass at $z$ in any measurable
space.  When the sequence of measures $(\R^{\F,\psi}_N)_{N\in\NN}$,
renormalized by $\psi'(N)$, converges (see Section
\ref{sect:latlogpaircor2} for background definitions) for the pointed
Hausdorff-Gromov weak-star convergence to a measure
$g_{\F,\psi}\,\Haar_{E_\infty}$ absolutely continuous with respect to the Haar
measure $\Haar_{E_\infty}$ of $E_\infty$, the Radon-Nikodym derivative
$g_{\F,\psi}$ is called the asymptotic {\it pair correlation
  function} of $\F$ {\em for the scaling} $\psi$ {\em and
  renormalisation} $\psi'$. When $g_{\F,\psi}$ is a positive constant,
we say that $\F$ has a {\em Poissonian behaviour for the scaling}
$\psi$ {\em and renormalisation} $\psi'$. When $g_{\F,\psi}$ vanishes
on a neighbourhood of $0$ in
$E_\infty$, we say that the pair $(\F,\psi)$ {\it exhibits a strong
  level repulsion}. The standard level repulsion only requires
$g_{\F,\psi}$ to vanish at $0$.

\medskip
Recall that a {\it $\ZZ$-lattice} in $\CC$ is a discrete (free
abelian) subgroup of $(\CC,+)$ generating $\CC$ as an $\RR$-vector
space. Let $\Lambda$ be a {\it $\ZZ$-grid} in $\CC$ (or {\it affine
  (Euclidean) lattice} in the terminology of
\cite{MarStr13,ElBMarVin15}), that is, a translate
$\Lambda=a+\vec\Lambda$ of a $\ZZ$-lattice $\vec\Lambda$ in the
Euclidean space $\CC$ for some $a\in\CC$ (modulo $\vec\Lambda$), see
for instance \cite{AkaEinSha16a}.  We denote by
$\covol_{\vec\Lambda}=\Vol(\CC/ \vec\Lambda)$ the area of a
fundamental parallelogram for $\vec\Lambda$. We denote by
\[
\operatorname{Sys}_{\vec\Lambda}=
\min\big\{|z|: z\in \vec\Lambda-\{0\}\big\} >0
\]
the {\it systole} of the $\ZZ$-lattice $\vec\Lambda$.  Recall that the
complex logarithm is an isomorphism of abelian topological groups
$\log:\CC^\times\ra E=\CC/(2\pi i\ZZ)$. Given $N\in\NN-\{0\}$ and a
function $\psi:\NN-\{0\} \ra\,]0,+\infty[$, we again denote by $\psi(N)$
    the scaling map from $E$ to $E_N=\CC/(2\pi i\psi(N)\ZZ)$
    defined by $z\!\!\mod 2\pi i\ZZ\mapsto \psi(N)z\!\!\mod 2\pi
    i\psi(N)\ZZ$. In Sections \ref{sect:latlogpaircor1} and
    \ref{sect:latlogpaircor2}, we study the pair correlations of the
    family of the complex logarithms of grid points
\[
\L_\Lambda=\big(L_N=\{\log z: z\in\Lambda,\;\;0<|z|\leq N\},\;
\omega_N=1\big)_{N\in\NN}
\]
without multiplicities. In order to simplify the statements in this
introduction, we only consider power scalings $\psi:N\mapsto N^\alpha$
for $\alpha\ge 0$, and we denote them by $\id^\alpha$.

\btheo\label{theo:intro1} Let $\alpha\ge 0$ and let $\Lambda$ be a
$\ZZ$-grid.  As $N\ra+\infty$, the normalized pair correlation
measures $\frac{1}{N^{4-2\alpha}}\; \R_N^{\,\L_\Lambda,\,\id^\alpha}$
on the cylinder $E_N=\CC/(2\pi i N^\alpha \ZZ)$ converge for the
pointed Hausdorff-Gromov weak-star convergence to the measure
$g_{\L_\Lambda,\,\id^\alpha}\;\Leb_{E_\infty}$ on $E_\infty=\CC/(2\pi
i \ZZ)$ if $\alpha=0$ and $E_\infty=\CC$ otherwise, with pair
correlation function given by
$$
g_{\L_\Lambda,\,\id^\alpha}:z\mapsto\begin{cases}
\frac{\covol_{\vec\Lambda}^2}{2\pi^3}\;e^{-2\,|\Re\,z|} &
\textrm{if } \alpha=0,\\
\frac{\pi}{2\covol_{\vec\Lambda}^2}&\textrm{if } 0<\alpha<1,\\
\frac{1}{\covol_{\vec\Lambda}\,|z|^4}
\sum\limits_{p\in\vec\Lambda\,:\,|p|\leq |z|}|p|^2
     &\textrm{if }\alpha=1,\\
     0&\textrm{if } \alpha>1\;.
\end{cases}
$$
The convergence is uniform on every compact subset of
$\ZZ$-grids $\Lambda$ for the Chabauty topology.
\etheo

The renormalisation by $\frac{1}{N^{4-2\alpha}}$ in Theorem
\ref{theo:intro1} is naturally chosen in order for the pair
correlation function to be finite. We refer to Theorems
\ref{theo:logpaircorrelnoscal} and \ref{theo:logpaircorrelscal} for
more complete versions of Theorem \ref{theo:intro1}, with more general
scaling functions, as well as for error terms. These error terms, as
well as the ones in Theorems \ref{theo:latlogpaircorrelphi} and
\ref{theo:logpaircorrelphipsi}, constitute the main technical parts of
this paper.

A standard scaling function in dimension $n$ is by the inverse of
$n$-root of the {\it average volume gap}, which is the quotient of the
volume of the ball of smallest radius containing $F_N$ by the number
of elements in $F_N$. See for instance
\cite{Montgomery73,RudSar98,BocZah05, LarSto20,HofKal21}, though these
references are in dimension $n=1$.  For the family $\L_\Lambda$, this
average volume gap is equivalent to $\frac{(\ln N)^2}{N^2}$, up to a
positive multiplicative constant. As we shall see in Theorem
\ref{theo:logpaircorrelscal}, the corresponding scaling function
$\psi:N\mapsto \frac{N}{\ln N}$ gives, as for $\psi:N\mapsto N^\alpha$
for $0<\alpha<1$ in the above theorem, a Poissonian behavior (see also
\cite{Vanderkam99a,ElBMarVin15} for a similar behaviour).

There is a phase transition from a Poissonian behaviour when
$0<\alpha<1$ to a total loss of mass at $\alpha>1$.  In fact, the
support of the measure itself converges to infinity for
$\alpha>1$. The transition occurs at the linear scaling, where an
exotic pair correlation function $g_{\L_\Lambda,\,\id}$ appears, which
has a discontinuity along every circle (centered at $0$) through a
grid point.  Since $g_{\L_\Lambda,\,\id}(z)$ vanishes when $z\in
\;\stackrel{\circ}{B} \!\!(0, \operatorname{Sys}_{\vec\Lambda})$, the
pair $(\L_\Lambda,\,\id)$ exhibits a strong level repulsion. Hence
$g_{\L_\Lambda,\,\id}$ has near $z=0$ a behaviour similar to the case
$\alpha>1$. Note that (see Lemma \ref{lem:liminftytheta})
$g_{\L_\Lambda,\,\id}(z)$ converges to
$\frac{\pi}{2\covol_{\vec\Lambda}^2}$ when $z$ goes to $\infty$,
corresponding to the Poissonian behaviour of $0<\alpha<1$.

The figure below gives the graph of the pair correlation function
$g_{\L_\Lambda,\,\psi}$ of $\L_\Lambda$ for the $\ZZ$-grid (which is a
$\ZZ$-lattice) $\Lambda=\vec\Lambda= \ZZ[i]$ of the Gaussian integers
at the linear scaling $\psi: N\mapsto N$ in the ball of center $0$ and
radius $5$. The blue lines on the bounding box represent the limit
$\frac{\pi}{2\covol_{\vec\Lambda}^2}=\frac{\pi}{2}$ at $+\infty$ of
$g_{\L_\Lambda,\,\psi}$.  We refer to the end of Section
\ref{sect:latlogpaircor2} for further illustrations, also in the case
of the Eisenstein integers.

\begin{center}
\includegraphics[width=13cm]{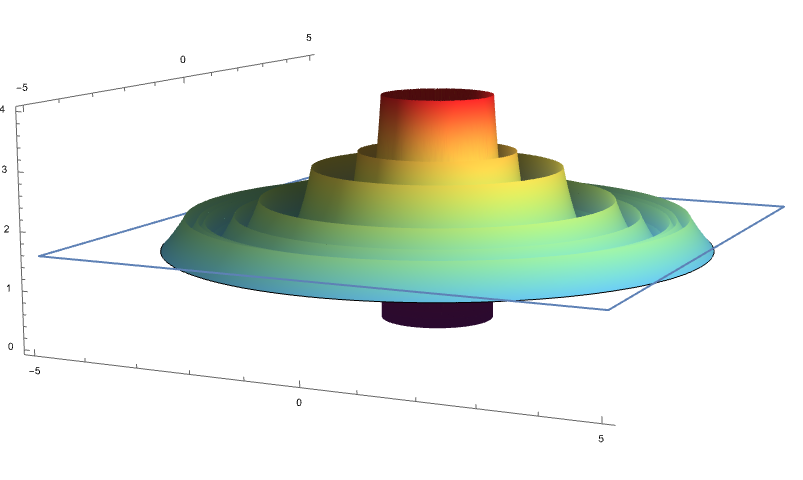}
\end{center}

We now give some existence results of pair correlation functions of
logarithms of lattice points with weights, restricting to integral
lattices with an arithmetic weight motivated by geometric
applications. Let $K$ be an imaginary quadratic number field $K$, with
discriminant $D_K$, whose ring of integers $\OOO_K$ is principal. We
fix a nonzero ideal $\Lambda$ in $\OOO_K$, and we denote by
$\varphi_K: \OOO_K-\{0\} \ra\NN$ the {\it Euler function} $a\mapsto
\card\big( (\OOO_K/a\OOO_K)^\times \big)$ of $K$. In the products
below, $\ppp$ runs over the prime ideals of $\OOO_K$.  The following
result describes the asymptotic behaviour of the pair correlation
measures associated with the family
\begin{equation}\label{eq:defiLLambdavarphi}
\L_\Lambda^{\varphi_K}=\big(L_N=\{\log z: z\in\Lambda,\;\;0<|z|\leq N\},\;
\omega_N=\varphi_K\circ\exp\big)_{N\in\NN}\;.
\end{equation}

\btheo\label{theo:intro2} (1) As $N\ra+\infty$, the pair correlation
measures $\R^{\L_\Lambda^{\varphi_K},1}_N$ on the constant cylinder
$E=\CC/(2\pi i \ZZ)$, renormalized to be probability measures,
weak-star converge to the probability measure
$g_{\L_\Lambda^{\varphi_K},1}\; \Leb_E$, with pair correlation
function independent of $\Lambda$ given by
$g_{\L_\Lambda^{\varphi_K},1}:z'\mapsto
\frac{1}{\pi}\,e^{-\,4\,|\Re\;z'|}$.

\medskip
(2) As $N\ra+\infty$, the normalized pair correlation measures
$\frac{1}{N^6}\,\R^{\L_{\OOO_K}^{\varphi_K},\,\id^1}_N$ on the varying
cylinder $E_N=\CC/(2\pi i\,N\, \ZZ)$ converge for the pointed
Hausdorff-Gromov weak-star convergence to the measure
$g_{\L_{\OOO_K}^{\varphi_K},\,\id^1}\; \Leb_\CC$, with pair
correlation function
\begin{equation}\label{eq:gLvarphiid}
g_{\L_{\OOO_K}^{\varphi_K},\,\id^1}: z\mapsto  \frac{2}{|z|^8\,\sqrt{|D_K|}\,}
\prod_{\ppp}\big(1-\frac{2}{\Nr(\ppp)^2}\big)\;
\sum_{\substack{k\in\OOO_K\\|k|\leq |z|}} |k|^6
\prod_{\ppp\,\mid\, k\OOO_K}\big(1+\frac{1}{\Nr(\ppp)(\Nr(\ppp)^2-2)}\big)\;.
\end{equation}
\etheo

We refer to Theorems \ref{theo:latlogpaircorrelphi} and
\ref{theo:logpaircorrelphipsi} for more complete versions of Theorem
\ref{theo:intro2}, including possible congruence restrictions, and for
error terms.  The proof of Theorem \ref{theo:intro2} (2) uses Theorem
1.3 of \cite{ParPau22c} that describes the asymptotic behaviour in
angular sectors in $\CC$ for the Euler function of $K$.  For the
readers convenience, we briefly review these results in Section
\ref{sect:MerMir}.  In order to simplify the treatment, we only
consider the constant and linear scaling in Theorem \ref{theo:intro2}.

\begin{center}
\includegraphics[width=13cm]{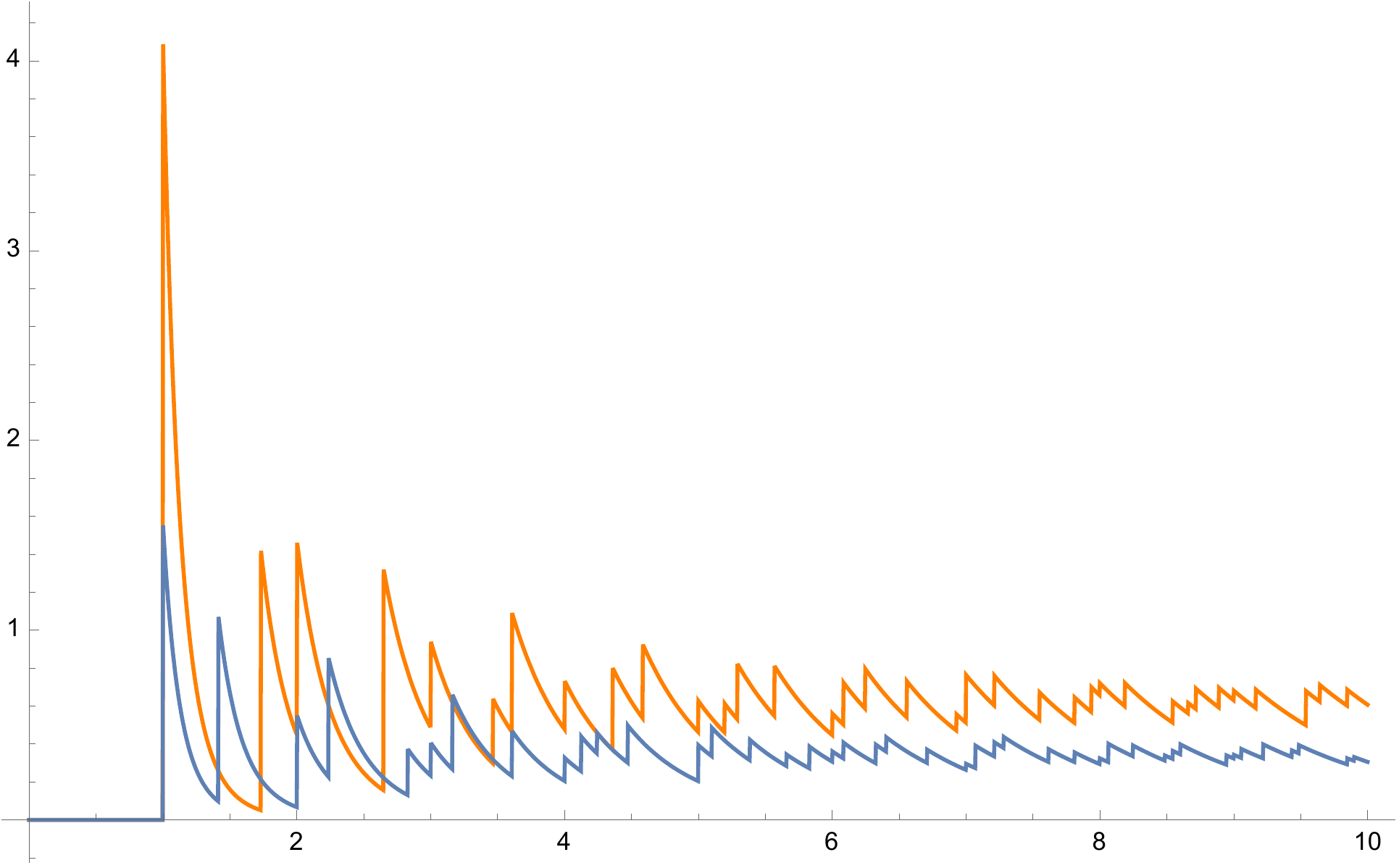}
\end{center}

The pair correlation functions at the linear scaling are radially
symmetric by Theorem \ref{theo:intro2} (2).  The figure above compares
the radial profiles of the pair correlation functions
$g_{\L_\Lambda^{\varphi_K} ,\,\id^1}$ for $K=\QQ(i)$ and $\Lambda=
\OOO_K =\ZZ[i]$ in blue and $K=\QQ(i\sqrt 3)$ and $\Lambda =\OOO_K
=\ZZ[\frac{1+i\sqrt 3}2]$ in orange.  The radial profiles of the pair
correlation functions converge to a limit
\[
\frac{\pi}{|D_K|} \;\prod_{\ppp}
\big(1-\frac{2}{\Nr(\ppp)^2} \big) \big(1+\frac{1}{\Nr(\ppp)^2
  (\Nr(\ppp)^2-2)}\big)
\]
at infinity, where $\ppp$ ranges over the prime ideals of $\OOO_K$ by
Proposition \ref{prop:horizasymptot}. This limit is approximately
$0.346$ for the blue curve and $0.634$ for the orange one.

The radial profiles of the pair correlation functions in the weighted
and unweighted cases are similar to certain radial distribution
functions in statistical physics, see for example
\cite[Sect.~II]{ZerPri27}, \cite[Fig.~7]{SanLop16}, \cite[page
  199]{Chandler87} or \cite[page 18]{Bohigas91}. See also
\cite{MatEtAl}. The unfolding technique (see for instance
\cite[p.~14]{Bohigas91} and \cite[\S 3, \S 5]{MarStr13}), though
guiding the very first step of the proofs of Theorem \ref{theo:intro1} 
and \ref{theo:intro2}, falls short of giving a complete answer, in
particular when varying the scalings and weights and for the error
term analysis.

\medskip
As explained in Section \ref{sect:geometricmotivation}, our motivation
for introducing the weights by the Euler function comes from
hyperbolic geometry. We prove in Proposition
\ref{prop:corrhoromodular} that the pair correlation measures of the
lengths (counted with multiplicity) of the common perpendiculars
between the maximal Margulis cusp neighbourhood and itself in the
(one-cusped) Bianchi orbifold $\PSL_2(\OOO_K)\bs\htr$ are closely
related to the pair correlation measures of the weighted family
$\L_{\OOO_K}^{\varphi_K}$. Theorem \ref{theo:intro2} implies a
pair correlation result for the lengths of common perpendiculars of
cusps neighborhoods in the Bianchi orbifold $\PSL_2(\OOO_K)/\htr$, see
Corollary \ref{coro:comperphoromodular} for a  precise statement and a 
version with congruences.

\medskip
\noindent{\small {\it Acknowledgements: } This research was supported
  by the French-Finnish CNRS IEA BARP and PaCap. }

\medskip\noindent {\bf Notation.} We introduce here some of the
notation used throughout the paper.

All our measures are Borel, positive, regular measures on locally
compact spaces. The pushforward of a measure $\mu$ by a mapping $f$ is
denoted by $f_*\mu$, and its total mass by $\|\mu\|$. We denote by
$\operatorname{Leb}_K$ the restriction of Lebesgue's measure of $\CC$
to any Borel subset $K$ of $\CC$. For every smooth manifold with
boundary $Y$ and every $k\in\NN$, we denote by $C^k_{\rm c}(Y)$ the
set of complex-valued $C^k$ functions with compact support on $Y$.

We equivariantly identify the space $\operatorname{Grid}_2$ of
$\ZZ$-grids in the real Euclidean plane $\CC$, endowed with the
Chabauty topology and the affine action of $\GL_2(\RR)\ltimes \RR^2$
with the homogeneous space $(\GL_2(\RR)\ltimes \RR^2)/ (\GL_2(\ZZ)
\ltimes \ZZ^2)$, which smoothly fibers by the map $a+\vec\Lambda
\mapsto \vec\Lambda$ over the space of $\ZZ$-lattices $\GL_2(\RR)/
\GL_2(\ZZ)$, with fibers the elliptic curves $\CC/\vec\Lambda$.

We will use the following indexing sets in Sections
\ref{sect:latlogpaircor1}, \ref{sect:latlogpaircor2} and
\ref{sect:latlogpaircor3}. Given a $\ZZ$-grid $\Lambda$, for every
$N\in\NN-\{0\}$, let
$$
I_N=I_{N,\Lambda}=\{(m,n)\in\Lambda^2\;:\;
0<|m|,|n|\leq N,\;m\neq n\}\;,
$$
$$
I^-_N=\{(m,n)\in I_N\;:\;|m|\leq |n|\}\;\;\;{\rm and}\;\;\;
I^+_N=\{(m,n)\in I_N\;:\; |n|\leq |m|\}\;.
$$

Given a subset $b$ of the set of ambient parameters, for every
positive function $g$ of a variable in $\NN-\{0\}$, we will denote by
$\bigO_b(g)$ (and $\bigO(g)$ when $b$ is empty) any function $f$ on
$\NN-\{0\}$ such that there exists a constant $C'$ depending only on
the parameters in $b$ and a constant $N_0$ possibly depending on the
all the parameters (including the ones in $b$) such that for every
$N\geq N_0$, we have $|f(N)|\leq C\;|g(N)|$.

\section{Pair correlation of grid points without weight or scaling}
\label{sect:latlogpaircor1}

In this section, we work on the constant cylinder $E=\CC/(2\pi i\ZZ)$,
endowed with its quotient Riemann surface structure, with its quotient
additive abelian locally compact group structure, and with its Haar
measure $d\Leb_E(x'+iy')=dx'dy'$ where $x'\in\RR$ and $y'\in\RR/(2\pi
\ZZ)$. We endow the multiplicative group $\CC^\times$ with its Riemann
surface structure as an open subset of $\CC$ and with the restriction
of the Lebesgue measure $\Leb_\CC$ of $\CC$. The logarithm map $\log :
\CC^\times\ra E$ defined by $\rho \,e^{i\theta}\mapsto \ln \rho +
i\theta$ is a biholomorphic group isomorphism, whose inverse is the
exponential map $z'=x'+iy'\mapsto \exp (z')=e^{x'}e^{iy'}$. The real
part map $\Re:E\ra\RR$ defined by $x'+iy'\mapsto x'$ is a smooth
(trivial) fibration, and
\begin{equation}\label{eq:pushLebEparRe}
  \Re_*\Leb_E=2\pi\Leb_\RR\;.
\end{equation}
Note that for every $z\in\CC-\{0\}$, we have
\begin{equation} \label{eq:relatlognormRenlog}
  \ln(|z|^2)=2\;\Re(\log z)\;.
\end{equation}
Since $d\Leb_\CC(\rho e^{i\theta})=\rho \,d\rho \,d\theta$, we have
\begin{equation}\label{eq:poussmeslog}
d(\log_*\Leb_\CC)(z')=e^{2\,\Re (z')}\;d\Leb_E(z')\;.
\end{equation}

Let $\Lambda=a+\vec\Lambda$ be a $\ZZ$-grid.  We choose a $\ZZ$-basis
$(v_1,v_2)$ of $\vec\Lambda$ such that the (weak) fundamental
parallelogram
$$
\F_{\vec\Lambda}=\big\{s\,v_1+t\,v_2\;:\;
s,t\in[-\frac{1}{2},\frac{1}{2}]\,\big\}
$$
for the action of $\vec\Lambda$ on $\CC$ has smallest diameter. We then
denote by
$$
\diam_{\vec\Lambda}=\diam(\F_{\vec\Lambda})=\max\{|v_1+v_2|,\;|v_1-v_2|\}
$$
the diameter of $\F_{\vec\Lambda}$, which is the length of a longest
diagonal of the parallelogram $\F_{\vec\Lambda}$. We denote by
\[
\covol_{\vec\Lambda}=\Vol(\CC/ \vec\Lambda)=
\operatorname{Area}(\F_{\vec\Lambda})= |\det(v_1,v_2)\,|
\]
the area of the elliptic curve $\CC/ \vec\Lambda$ for the measure
induced by the Lebesgue measure on $\CC$, or the area of the
parallelogram $\F_{\vec\Lambda}$ (which does not depend on the choice
of the $\ZZ$-basis $(v_1,v_2)$ of $\vec\Lambda$).  We will use several
times the following well known result, having a more precise error
term that we won't need, and we only give a proof in order to
explicit the dependence on the parameter $\Lambda$.

\blemm\label{lem:liminftytheta} For every $k\in\NN$, as $x\ra+\infty$,
we have
\begin{equation}\label{eq:gaussRquat}
  \sum_{p\in\Lambda\,:\,|p|\leq x} |p|^k =
  \frac{2\pi}{(k+2)\covol_{\vec\Lambda}} \;x^{k+2} +\bigO_k\big(\,
  \frac{k+\diam_{\vec\Lambda}}{\covol_{\vec\Lambda}} \,x^{k+1}\big)\;.
\end{equation}
\elemm

\noindent \begin{minipage}{8.5cm}
\dem The case $k=0$ is the standard Gauss counting result of lattice
points in discs. With $A_x=\{p\in\Lambda\,:\,|p|\leq x\}$ and
$B_x=\bigcup_{p\in A_x}(p+\F_{\vec\Lambda})$, so that $\operatorname{Area}
(B_x)= \card(A_x)\;\operatorname{Area} (\F_{\vec\Lambda})$, we have
$$
B\big(0,x-\diam_{\vec\Lambda}\big) \subset \;B_x\;
\subset B\big(0,x+\diam_{\vec\Lambda}\big)\;,
$$
so that the result for $k=0$ follows by computing the area of the two above
discs.
\end{minipage}
\begin{minipage}{6.3cm}
\begin{center}
\begin{picture}(0,0)%
\includegraphics{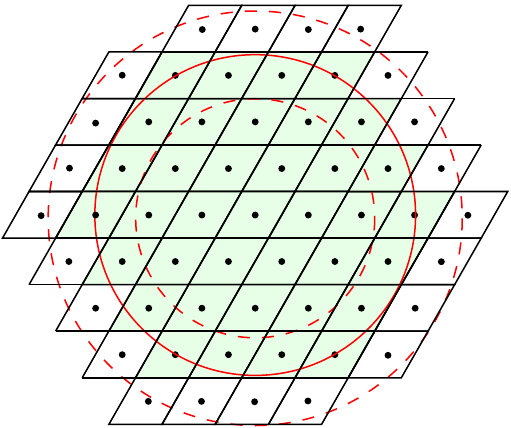}%
\end{picture}%
\setlength{\unitlength}{3812sp}%
\begingroup\makeatletter\ifx\SetFigFont\undefined%
\gdef\SetFigFont#1#2#3#4#5{%
  \reset@font\fontsize{#1}{#2pt}%
  \fontfamily{#3}\fontseries{#4}\fontshape{#5}%
  \selectfont}%
\fi\endgroup%
\begin{picture}(2535,2107)(538,-2111)
\put(1715,-1140){\makebox(0,0)[lb]{\smash{{\SetFigFont{8}{9.6}{\rmdefault}{\mddefault}{\updefault}{\color[rgb]{0,0,0}$0$}%
}}}}
\end{picture}%

\end{center}
\end{minipage}

\medskip
Assume now that $k\geq 1$. We use Abel's summation formula
\[
\sum_{1\leq n\leq x}a_nf(n)=\big(\sum_{1\leq n\leq x} a_n\big)\;f(x) -
\int_{1}^x\big(\sum_{1\leq n\leq t} a_n\big)\;f'(t)\;dt
\]
applied to the numerical sequence $\big(a_n =\card\{p\in\Lambda:
n-1<|p|\leq n\}\big)_{n\in\NN}$ and to the smooth functions
$f:[1,+\infty[\;\ra\RR$ defined by $t\mapsto t^k$ or $t\mapsto
(t-1)^k$.  Using the case $k=0$, the result when $k>0$ then follows
from the estimates
\[
\sum_{1\leq n\leq x} a_n\,(n-1)^k\leq \sum_{p\in\Lambda\,:\,|p|\leq x} |p|^k
\leq \sum_{1\leq n\leq x} a_n\,n^k\;. \;\;\;\Box
\]

\medskip
For every $N\in\NN-\{0\}$, the (not normalised) {\em pair correlation
  measure} of the logarithms of nonzero grid  points in $\Lambda$,
with trivial multiplicities and with trivial scaling function, is the
finite measure on the cylinder $E$ defined by
\[
\nu_N =\nu_{N,\Lambda}=\R_N^{\L_\Lambda,1}=
\sum_{(m,\,n)\in I_N}\;\Delta_{\log m-\log n}\;.
\]
Note that for every $k\in\NN-\{0\}$, we have $I_{kN,k\Lambda}
=I_{N,\Lambda}$ and $\nu_{kN,k\Lambda} =\nu_{N,\Lambda}$.

\btheo\label{theo:logpaircorrelnoscal} As $N\ra+\infty$, the measures
$\nu_N$ on $E$, renormalized to be probability measures, weak-star
converge to the measure absolutely continuous with respect to the Haar
measure $\Leb_E$, with Radon-Nikodym derivative the function
$g_{\L_\Lambda,1}:z'\mapsto\frac{1}{2\pi}\,e^{-2\,|\Re(z')|}$, which is
independent of $\Lambda$. Besides, the convergence
\begin{equation}\label{eq:logpaircorrelnoscal}
\frac{\nu_{N,\Lambda}}{\|\nu_{N,\Lambda}\|}\;\;\;\weakstar\;\;\;
g_{\L_\Lambda,1}\;\Leb_E
\end{equation}
is uniform on every compact subset of $\Lambda$ in the space of
$\ZZ$-grids $\operatorname{Grid}_2$. Furthermore, for every $f\in
C^1_{\rm c}(E)$, we have
$$
\frac{\nu_N}{\|\nu_N\|}(f)=
\frac{1}{2\pi}\,\int_{E} f(z')\;e^{-2\,|\Re(z')|}\;d\Leb_E(z')+ \bigO\Big(\,
\frac{\diam_\Lambda}{N}\, (\|f\|_\infty+\|e^zdf(z)\|_\infty)\Big)\;.
$$
\etheo

This result implies the case $\alpha=0$ of Theorem \ref{theo:intro1}
in the introduction, since we will prove in Formula
\eqref{ed:valasympnuN} that $\lim_{N\ra+\infty}\frac{\|\nu_N\|}{N^4}=
\frac{\pi^2}{\covol_{\vec\Lambda}^2}$.

\brema\label{rem:mneqn}{\rm Theorem \ref{theo:logpaircorrelnoscal} is
  still valid if we allow $n=m$ in the definition of the index set
  $I_{N}$ (this correspond to removing the condition $p\neq q$ in the
  definition below of $J_q$), see also Remark (2) in \cite[\S
    3]{ParPau22b} for a general argument. We will use this comment in
  Corollary \ref{coro:mmm} and \ref{coro:applir2d2}, as well as in
  Corollary \ref{coro:comperphoromodular}.}  \erema

\dem
For all $N\in\NN$ and $q\in\Lambda$ with $0<|q|\leq N$, let
\begin{equation}\label{eq:defJsubq}
  J_q=\{p\in\Lambda :  0<|p|\leq |q|,\; p\neq q\}
  \;\;\;\;{\rm and}\;\;\;\;\omega_q=\sum_{p\in J_q}\Delta_{\frac{p}{q}}\;,
\end{equation}
which is a finitely supported measure on the closed unit disc $\DD$ of
$\CC$. Note that the assumptions $0<|p|$ and $0<|q|$ are automatic
when $0\notin \Lambda$, that is, when $\Lambda$ is not a
$\ZZ$-lattice. As $q\ra+\infty$, by Equation \eqref{eq:gaussRquat}
with $k=0$, its total mass, which is nonzero since $-q\in J_q$,
satisfies
\begin{equation}\label{eq:normomegaq}
\|\omega_q\|=\frac{\pi}{\covol_{\vec\Lambda}}\,|q|^2+
\bigO\big(\,\frac{\diam_{\vec\Lambda}}{\covol_{\vec\Lambda}}\,|q|\big)\;,
\end{equation}
for some $\bigO(\cdot)$ uniform in $\Lambda$. Note that we need to
remove $0$ if $0\in\Lambda$ and $q$ from the counting of Equation
\eqref{eq:gaussRquat}, but this is taken care of by the above
$\bigO(\cdot)$. We hence have
$$
\frac{1}{\|\omega_q\|} =\frac{\covol_{\vec\Lambda}}{\pi\,|q|^2}+
\bigO\big(\,\frac{\diam_{\vec\Lambda}\covol_{\vec\Lambda}}{|q|^3}\,\big)\;,
$$
for some $\bigO(\cdot)$ uniform in $\Lambda$. We denote by
$\overline{\omega_q} =\frac{\omega_q}{\|\omega_q\|}$ the
renormalisation of $\omega_q$ to a probability measure on $\DD$.

Let $f\in C^1(\DD)$. Let
\[C_q=\bigcup_{p\in J_q}(p+\F_{\vec\Lambda})\;.
\]
Note that the symmetric difference $(\DD- \frac{C_q}{q})
\cup(\frac{C_q}{q}-\DD)$ is contained in the union of the annulus
$B(0,1+ \frac{\diam_{\vec\Lambda}}{|q|}\big) - B\big(0,
1-\frac{\diam_{\vec\Lambda}} {|q|}\big)$ and (when $0\in\Lambda$) the
parallelogram $\frac{\F_{\vec\Lambda}}{q}$, hence has area at most
\[
\frac{\covol_{\vec\Lambda}}{|q|^2}+
\pi\big((1+\frac{\diam_{\vec\Lambda}}{|q|})^2 -(1-\frac{\diam_{\vec\Lambda}}
{|q|})^2\big)=\bigO\big(\frac{\diam_{\vec\Lambda}}{|q|}\big)\;.
\]
Also note that $\frac{\diam_{\vec\Lambda}\covol_{\vec\Lambda}}
{|q|^3}|\omega_q(f)|= \bigO\big(
\frac{\diam_{\vec\Lambda}\|f\|_\infty}{|q|}\big)$. Therefore
\begin{align*}
&\Big|\frac{1}{\pi}\int_\DD f(z)\;d\Leb_\CC(z)\;-\;
\overline{\omega_q} (f)\,\Big|\\=\;&
\Big|\frac{1}{\pi}\int_{\frac{C_q}{q}} f(z)\;d\Leb_\CC(z)\;-\;
\frac{\covol_{\vec\Lambda}}{\pi\,|q|^2}\omega_q(f)
\Big|+\bigO\big(\frac{\diam_{\vec\Lambda}\|f\|_\infty}{|q|}\big)\;.
\end{align*}
By the mean value inequality, for all $p\in J_q$ and
$z\in\frac{p+\F_{\vec\Lambda}}{q}$, we have
$$
\big|\,f(z)-f(\frac{p}{q})\,\big|\leq
\|df\|_\infty\;\big|\,z-\frac{p}{q}\,\big|\;.
$$
Hence
\begin{align*}
 \Big|\frac{1}{\pi}\int_{\frac{C_q}{q}} f(z)\;d\Leb_\CC(z)\;-\;
\frac{\covol_{\vec\Lambda}}{\pi\,|q|^2}\,\omega_q(f)
\Big| & = \frac{1}{\pi}\Big|\sum_{p\in J_q}\int_{\frac{p+\F_{\vec\Lambda}}{q}}
\big(f(z)-f(\frac{p}{q})\big)\;d\Leb_\CC(z)\Big| \\ &\leq
\frac{1}{\pi}\sum_{p\in J_q}\|df\|_\infty\;\sup_{z\in\frac{p+\F_{\vec\Lambda}}{q}}
\big|\,z-\frac{p}{q}\,\big|\;\operatorname{Area}(\frac{p+\F_{\vec\Lambda}}{q})
 \\ &\leq \frac{1}{\pi}\;\card(J_q)\;\frac{\diam_{\vec\Lambda}}{|q|}\;
 \frac{\covol_{\vec\Lambda}}{|q|^2}\;\|df\|_\infty \\ &
 =\bigO\Big(\frac{\diam_{\vec\Lambda}\;\|df\|_\infty}{|q|}\Big)\;.
\end{align*}
Therefore
\begin{equation}\label{eq:limomegaq}
\overline{\omega_q} (f)=\frac{1}{\pi}\int_\DD f(z)\;d\Leb_\CC(z)+
\bigO\Big(\frac{\diam_{\vec\Lambda}\,(\|f\|_\infty+\|df\|_\infty)}{|q|}\Big)\;.
\end{equation}
In particular, as $q\ra+\infty$, we have $\overline{\omega_q}
\;\;\weakstar\;\;\frac{1}{\pi}\Leb_{\DD}$.

Let us now define
\[
\mu^-_N =\sum_{(m,\,n)\in I^-_N}\;\Delta_{\frac{m}{n}}\;\;
=\sum_{q\in\Lambda,\; 0<|q|\leq N} \;\;\omega_q\;,
\]
which is a finitely supported measure on $\DD$.  By Equations
\eqref{eq:normomegaq} and \eqref{eq:gaussRquat} with $k=2$ and $k=1$,
its total mass is equal to
\begin{align*}
\|\mu^-_N\|=\sum_{q\in\Lambda,\; 0<|q|\leq N} \|\omega_q\|&=
\sum_{q\in\Lambda,\; 0<|q|\leq N} \big(\frac{\pi}{\covol_{\vec\Lambda}}\,|q|^2+
\bigO\big(\,\frac{\diam_{\vec\Lambda}}{\covol_{\vec\Lambda}}\,|q|\big)\;\big)
\\ &=\frac{\pi^2}{2\covol_{\vec\Lambda}^2}\,N^4+ \bigO\Big(
\frac{1+\diam_{\vec\Lambda}}{\covol_{\vec\Lambda}^2}N^3\Big)\;.
\end{align*}
It follows that
\begin{equation}\label{eq:asympinvmuN}
  \frac{1}{\|\mu^-_N\|} =\frac{2\covol_{\vec\Lambda}^2}
{\pi^2\,N^4} + \bigO\big(\frac{(1+\diam_{\vec\Lambda})\covol_{\vec\Lambda}^2}{N^5}
\big)\;.
\end{equation}

Let $f\in C^1(\DD)$. By Equations \eqref{eq:limomegaq},
\eqref{eq:asympinvmuN}, and \eqref{eq:gaussRquat} with $k=1$, we have,
as $N\ra+\infty$,
\begin{align}
  \frac{\mu^-_N(f)}{\|\mu^-_N\|}& =
  \frac{1}{\|\mu^-_N\|}\sum_{q\in\Lambda,\; 0<|q|\leq N}
\;\;\|\omega_q\|\;\overline{\omega_q}(f)\nonumber\\ &=
\frac{1}{\pi}\int_\DD f(z)\;d\Leb_\CC(z)+\sum_{q\in\Lambda,\; 0<|q|\leq N}
\frac{\|\omega_q\|}{\|\mu^-_N\|}
\bigO\Big(\frac{\diam_{\vec\Lambda}\,(\|f\|_\infty+\|df\|_\infty)}{|q|}\Big)
\nonumber\\ & =
\frac{1}{\pi}\int_\DD f(z)\;d\Leb_\CC(z)+
\bigO\Big(\,\frac{\diam_{\vec\Lambda}}{N}\,(\|f\|_\infty+\|df\|_\infty)\Big)\;.
\label{eq:sansfacteur2}
\end{align}

Let $E^\pm=(\pm[0,\infty[\,+i\RR)/(2\pi i \ZZ)$ so that $E=E^-\cup
    E^+$.  Note that $\log:\DD-\{0\}\ra E^-$ and $\log:
    \overline{\CC-\DD}\ra E^+$ are homeomorphisms. Let us define a
    measure with finite support on $E^\pm$ by
$$
\nu^\pm_N =\sum_{(m,\,n)\in I^\pm_N}\;\Delta_{\log\frac{m}{n}}\;,
$$ so that $\nu^-_N=\log_*\mu_N^-=\nu_N\!\mid_{E^-}$, and
$\|\nu_N^-\|=\|\mu_N^-\|$. For every $f\in C^1_{\rm c}(E^-)$, we have
$f\circ\log\in C^1_{\rm c}(\DD-\{0\})$ (hence $f\circ\log$ may be
extended to a $C^1$ function on $\DD$ which vanishes on a neighborhood
of $0$) and, by Equations \eqref{eq:sansfacteur2} and
\eqref{eq:poussmeslog},
\begin{align*}
  \frac{\nu^-_N(f)}{\|\nu^-_N\|}& =
  \frac{\mu^-_N(f\circ\log)}{\|\mu^-_N\|} \\ & =\frac{1}{\pi}
  \int_\DD f\circ\log(z)\;d\Leb_\CC(z)+ \bigO\Big(\,
  \frac{\diam_{\vec\Lambda}}{N}\, (\|f\circ\log\|_\infty+
  \|d(f\circ\log)\|_\infty)\Big) \\ & =\frac{1}{\pi}
  \int_{E^-} f\;e^{2\,\Re (z')}\;d\Leb_E(z')+ \bigO\Big(\,
  \frac{\diam_{\vec\Lambda}}{N}\, (\|f\|_\infty+
    \|e^zdf(z)\|_\infty)\Big)\;.
\end{align*}
Let $\sg:E\ra E$ be the horizontal change of sign map $x'+iy'\mapsto
-x'+iy'$, which maps $E^-$ to $E^+$. Then $\nu_N^+=\sg_*\nu^-_N$ and
$\nu_N=\nu_N^-+\nu_N^+$.  Since $E^-\cap E^+$ has zero measure for the
Haar measure $\Leb_E$ and since $\|\nu_N^\pm\|=
\frac{1}{2}\,\|\nu_N\|+ \bigO(\diam_{\vec\Lambda} N^3)$, the last claim of
Theorem \ref{theo:logpaircorrelnoscal} follows. Note that, as
needed just after the statement of Theorem
\ref{theo:logpaircorrelnoscal}, as $N\ra+\infty$, we have
\begin{equation}\label{ed:valasympnuN}
  \|\nu_N\|\sim2\|\mu_N^-\|\sim\frac{\pi^2}{\covol_{\vec\Lambda}^2}\,N^4\;.
\end{equation}

The first claim follows by approximating continuous functions with
compact support by $C^1$ ones. The uniformity of the convergence on
compact subsets of lattices follows from the uniformity of the
functions $\bigO(\cdot)$ and the fact that the constants
$\covol_{\vec\Lambda}$ and $\diam_{\vec\Lambda}$ vary in a compact subset of
$]0,+\infty[$ when $\Lambda$ varies in a compact subset of
$\operatorname{Grid}_2$.  \cqfd

\medskip
The following picture illustrates the convergence in Formula
\eqref{eq:logpaircorrelnoscal} in Theorem
\ref{theo:logpaircorrelnoscal} when $\Lambda=\vec\Lambda=\ZZ[i]$ is
the ring of Gaussian integers and $N=20$, using as horizontal
coordinates $(x',y')\in E$ with $x'\in\RR$ and $y'\in [-\pi,\pi[$. A
    smooth histogram scaled to a probability density is displayed in
    gold, and the limiting distribution in grey.

\bigskip\bigskip
\begin{center}
\includegraphics[width=14cm]{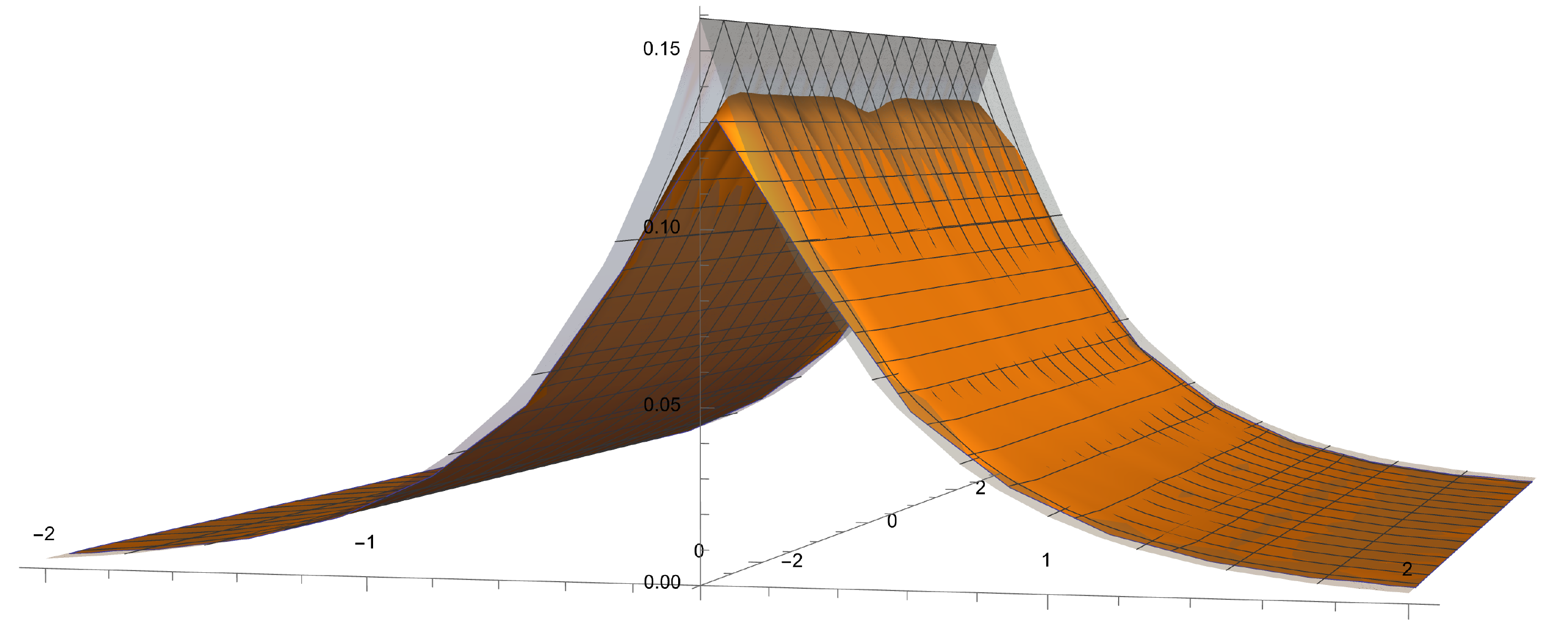}
\end{center}
\bigskip\bigskip

\medskip
\noindent{\bf Arithmetic applications. } (1) Let $K$ be an imaginary
quadratic number field, with discriminant $D_K$, ring of integers
$\OOO_K$ and Dedekind zeta function $\zeta_K$.  We denote by $\I^+_K$
the semigroup of nonzero (integral) ideals of the Dedekind ring
$\OOO_K$ (with unit $\OOO_K$).  We denote by $\Nr (I)=\card
(\OOO_K/I)$ the {\it norm} of an ideal $I\in\I^+_K$, which is
completely multiplicative. The {\it norm} of $a\in \OOO_K-\{0\}$ is
\[
\Nr (a)=\Nr (a\OOO_K)\;.
\]
It coincides with the (relative) norm $N_{K/\QQ}(a)$ of $a$ (see for
instance \cite{Narkiewicz04}), and in particular is equal to $|a|^2$
since $K$ is imaginary quadratic.  The {\em norm} of a fractional
ideal $\mmm$ of $\OOO_K$ is $\frac 1{|c|^2}\Nr(c\mmm)$ for any $c\in
\OOO_K-\{0\}$ such that $c\mmm \subset \OOO_K$.

Let $\mmm$ be a nonzero fractional ideal of $\OOO_K$.  Note that
$\mmm$ is a $\ZZ$-lattice in $\CC$ with
\begin{equation}\label{eq:calccocoldiamideal}
\covol_\mmm=\frac{\sqrt{|D_K|}\,\Nr (\mmm)}{2}
\;\;\;{\rm and}\;\;\;
\diam_\mmm=\bigO(\sqrt{|D_K|\,\Nr (\mmm)}\;)\;,
\end{equation}
since $\OOO_K=\ZZ+\frac{\sqrt{D_K}}{2}\ZZ$ and $\diam_{\OOO_K}=
|1+\frac{\sqrt{D_K}}{2}|$ if $D_K\equiv 0\!\!\mod 4$, and
$\OOO_K=\ZZ+\frac{1+\sqrt{D_K}}{2}\ZZ$ and $\diam_{\OOO_K}=
|\frac{3+\sqrt{D_K}}{2}|$ if $D_K\equiv 1\!\!\mod 4$. In particular,
the Gauss ball counting argument of Equation \eqref{eq:gaussRquat}
with $k=0$ and $x=\sqrt{N}$ gives, as $N'\ra+\infty$,
\begin{align*}
(\card\{m\in\mmm :0<\Nr (m)\leq N'\})^2&=
\Big(\frac{\pi}{\covol_\mmm}\,N'+
\bigO\big(\,\frac{\diam_\mmm}{\covol_\mmm}\,\sqrt{N'}\,\big)\Big)^2\\ &=
\frac{4\pi^2{N'}^2}{|D_K|\Nr (\mmm)^2}\Big(1+\bigO\big(
\frac{\sqrt{|D_K|\,\Nr (\mmm)}\;}{\sqrt{N'}}\big)\Big)\;.
\end{align*}
Hence Theorem \ref{theo:logpaircorrelnoscal} implies the existence of
a pair correlation function (independent of $\mmm$) for the family of
the complex logarithms of nonzero elements of $\mmm$
\[
\L_\mmm=\big(L_{N'}=\{\log n\;:\;n\in\mmm,\;0<\Nr (n)\leq N'\},
\omega_{N'}= 1\big)_{N'\in\NN}
\]
without weights or scaling, as stated in the following result, using
Remark \ref{rem:mneqn}.

\bcoro\label{coro:mmm}
For every $f\in C^1_{\rm c}(E)$, as $N'\ra+\infty$, we have
\begin{align*}
  &\frac{|D_K|\Nr (\mmm)^2}{4\pi^2{N'}^2}
  \sum_{m,n\in\mmm\;:\;0<\Nr (m),\Nr (n)\leq N'}
f(\log m-\log n)\\=\;&
\frac{1}{2\pi}\,\int_{E} f(z')\;e^{-2\,|\Re(z')|}\;d\Leb_E(z')+
\bigO\Big(\,\frac{\sqrt{|D_K|\,\Nr (\mmm)}}{\sqrt{N'}}\,
(\|f\|_\infty+\|e^zdf(z)\|_\infty)\Big)\;.\;\;\; \Box
\end{align*}
\ecoro

\medskip
(2) For every positive integer $d$, let $r_{2,d}:\NN-\{0\}\ra \NN$ be
the integral function where
\[
r_{2,d}(n)= \card\{(x,y)\in\ZZ^2:x^2+d\,y^2 =n\}
\]
is the number of integral solutions of the Diophantine equation
$x^2+d\,y^2 =n$, for every $n\in \NN$. In particular, if $d=1$, then
$r_{2,d}=r_2$ is the well known function counting the sum of two
squares representatives of a given positive integer (see for instance
\cite{Cox13} or \cite[Sect.~16.9]{HarWri08}).  The following result
proves that the map
\[
g_\RR:t\mapsto\frac{1}{2}\,e^{-|t|}
\]
on $\RR$ is the pair correlation function for the family
\[
\L_\NN^{r_{2,d}}=\big(F_N=\{\ln n\;:\; 0<n\leq N, \;r_{2,d}(n)\neq
0\},\; \omega_N=r_{2,d}\circ \exp\big)_{N\in\NN}
\]
of the logarithms of the nonzero natural integers, without scaling but
with weights given by $r_{2,d}$ (removing the zero weights). Other
weights have been considered in \cite{ParPau22a} (including the one
given by the Euler function $\varphi$). Note that the following
corollary holds also when $r_{2,d}(n)$ is replaced by the number of
representations of $n$ by the norm form of any imaginary quadratic
number field, evaluated on any order of their ring of integers (as for
instance the norm form $(x,y)\mapsto x^2-xy+y^2$ of the Eisenstein
integers).

\bcoro \label{coro:applir2d2} As $N\ra+\infty$, we have
\begin{align*}
  \frac{1}{\big(\sum_{0<m\leq N^2}r_{2,d}(m)\big)^2} &\;\;
  \sum_{m,n\in\NN\;:\;0<m,n\leq N^2} r_{2,d}(m)\;r_{2,d}(n)\;
  \Delta_{\ln m-\ln n}
  \;\;\;\weakstar\;\;\; g_\RR\;\Leb_\RR \;.
\end{align*}
\ecoro

\dem Let us consider the $\ZZ$-lattice $\Lambda=\ZZ+i\sqrt{d}\;\ZZ$ in
$\CC$. Using Remark \ref{rem:mneqn}, we remove the assumptions $m\neq
n$ in the summations defining $\R^{\L_\Lambda,1}_N$ as well as
$\R^{\L_\NN^{r_{2,d}},1}_{N^2}$.

By the linearity of $(2\,\Re)_*$ and $2\,\Re$, and by Equation
\eqref{eq:relatlognormRenlog}, for every $N\in\NN-\{0\}$, we have
\begin{align*}
  (2\,\Re)_*\big(\R^{\L_\Lambda,1}_N\big)&
  =\sum_{p,q\in\Lambda\;:\;0<|p|,|q|\leq N}\Delta_{2\,\Re(\log p)-2\,\Re(\log q)}
  \\ & = \sum_{0<m,n\leq N^2}\;\;\sum_{p,q\in\Lambda\;:\;|p|^2=m,|q|^2=n}
  \Delta_{\ln(|p|^2)-\ln(|q|^2)}\\ & =
  \R^{\L_\NN^{r_{2,d}},1}_{N^2}\;.
\end{align*}
The pushforward map $(2\,\Re)_*$ preserves the total mass and is
continuous for the weak-star topology, since the map $2\,\Re:E\ra \RR$
is proper. Hence by Formulas \eqref{eq:logpaircorrelnoscal} and
\eqref{eq:pushLebEparRe}, we have
\begin{align*}
  \frac{\R^{\L_\NN^{r_{2,d}},1}_{N^2}}{\|\R^{\L_\NN^{r_{2,d}},1}_{N^2}\|}=
  (2\,\Re)_*\Big(\frac{\R^{\L_\Lambda,1}_N}{\|\R^{\L_\Lambda,1}_N\|}\Big)
  \;\;\;\weakstar\;\;\;
  &(2\,\Re)_*\Big(\frac{1}{2\pi}\,e^{-\,|2\,\Re(z')|}\;d\Leb_E(z')\Big)\\ &=
  \frac{1}{2}\,e^{-\,|t|}\;d\Leb_\RR(t)\;.
\end{align*}
Corollary \ref{coro:applir2d2} follows. \cqfd

\section{Pair correlation of grid points with scaling without weight}
\label{sect:latlogpaircor2}

In this section, we study the pair correlations of complex logarithms
of grid points at various scaling. We fix a positive scaling function
$\psi: \NN-\{0\}\ra \;]0,+\infty[$ such that ${\displaystyle
      \lim_{+\infty}\;\psi=+\infty}$.  We consider a normalisation
function $\psi': \NN-\{0\}\ra \;]0,+\infty[$ depending on $\psi$,
which will be made precise later on, but which in most cases will
not yield the renormalisation to a probability measure.

We will work on the following family $(E_N)_{N\in \NN-\{0\}}$ of
varying cylinders. For every $N\in \NN-\{0\}$, we consider $E_N=\CC/
(2\pi i\,\psi(N)\,\ZZ)$, endowed with its quotient Riemann surface
structure and its quotient additive abelian locally compact group
structure. Since a real number $\theta$ is well defined modulo
$2\pi\ZZ$ if and only if $\psi(N)\theta$ is well defined modulo
$2\pi\psi(N)\ZZ$, the scaled logarithm map $\psi(N)\log :
\CC^\times\ra E_N$ defined by $\rho \,e^{i\theta}\mapsto \psi(N)\ln
\rho + i\psi(N)\theta$ is a biholomorphic group isomorphism, whose
inverse is the rescaled exponential map $z'=x'+iy'\mapsto \exp
(\frac{z'}{\psi(N)})= e^{\frac{x'}{\psi(N)}} e^{i\frac{y'}{\psi(N)}}$.
The real part map $\Re:\CC\ra\RR$ induces a map again denoted by
$\Re:E_N\ra\RR$, which is a trivial smooth bundle map with fibers
$i\RR/(2\pi i\psi(N)\ZZ)$, such that for every $z\in E$,
\begin{equation}\label{eq:RedilatRe}
\Re(\psi(N)z)=\psi(N)\;\Re(z)\;.
\end{equation}

\smallskip
We consider also $E_N$ as a pointed metric space, with distance the
quotient of the Euclidean distance on $\CC$ and base point its
(additive) identity element $0$. Note that $E_N$ is a proper metric
space.  As ${\displaystyle \lim_{+\infty} \;\psi=+\infty}$, for every
$R>0$, there exists $N_R\in\NN-\{0\}$ such that for every $N\geq N_R$,
the closed ball $B(0,R)$ in $\CC$ injects isometrically by the
canonical projection $p_N:\CC\ra E_N$. Hence the sequence $(E_N)_{N\in
  \NN-\{0\}}$ of proper pointed metric spaces converges to the proper
metric space $\CC$ pointed at $0$ for the pointed Hausdorff-Gromov
convergence (see \cite{Gromov99a} for background).

Any function $f\in C^0_{\rm c}(\CC)$ defines for all $N$ large enough
a function $f_N\in C^0_{\rm c}(E_N)$ as follows. Let $R_f>0$ be such
that the support of $f$ is contained in $B(0,R_f)$. Then for every
$N\geq N_{R_f}$, the function $f_N\in C^0_{\rm c}(E_N)$ is the
function which vanishes outside $p_N(B(0,R_f))$ and coincides with
$f\circ ({p_N}_{\mid B(0,R_f)})^{-1}$ on $p_N(B(0,R_f))$. Note that
$f_N$ is $C^1$ if $f$ is $C^1$.

We say that a sequence $(\mu_N)_{N\in \NN-\{0\}}$ of measures $\mu_N$
on $E_N$ converges to a measure $\mu_\infty$ on $\CC$ for the {\it
  pointed Hausdorff-Gromov weak-star convergence} if for every $f\in
C^0_{\rm c}(\CC)$, the sequence $(\mu_N(f_N))_{N\geq N_{R_f}}$ converges in
$\CC$ to $\mu_\infty(f_\infty)$ (see
\cite[Chap.~3$\frac{1}{2}$]{Gromov99a} for background). We again use
the symbol $\weakstar$ in order to denote this convergence.

\medskip
Let $\Lambda$ be a $\ZZ$-grid in $\CC$. For every $N\in\NN-\{0\}$,
the (not normalised, empirical) {\em pair correlation measure} of the
complex logarithms of points in $\Lambda$ at time $N$ {\em with
  trivial weights and with scaling} $\psi(N)$ is the measure with
finite support in $E_N$ defined by
$$
\R^{\L_\Lambda,\psi}_N=
\sum_{(m,\,n)\in I_N}\;\;\Delta_{\psi(N)\log m-\psi(N)\log n}\;,
$$
and the normalized one is $\frac{1}{\psi'(N)}\;\R^{\L_\Lambda,\psi}_N$.

\btheo\label{theo:logpaircorrelscal} Let $\Lambda=a+\vec\Lambda$ be a
$\ZZ$-grid in $\CC$. Assume that the scaling function $\psi$
satisfies ${\displaystyle \lim_{N\ra+\infty}}\; \frac{\psi(N)}{N}
=\lambda_\psi\in [0,+\infty]$.  As $N\ra+\infty$, the measures
$\R^{\L_\Lambda,\psi}_N$ on $E_N$, normalized by $\psi'(N)$ as given
below, converge for the pointed Hausdorff-Gromov weak-star convergence
to a measure $g_{\L_\Lambda,\psi}\;\Leb_\CC$ on $\CC$, absolutely
continuous with respect to the Lebesgue measure on $\CC$, with
Radon-Nikodym derivative the function
\begin{equation}\label{eq:troiscas}
g_{\L_\Lambda,\psi}:z\mapsto\begin{cases} 0 &{\rm ~~if~~}
\lambda_\psi=+\infty {\rm ~and~} \psi'=\psi\;, \\
\frac{\pi}{2\covol_{\vec\Lambda}^2}
  &{\rm ~~if~~}\lambda_\psi=0{\rm ~and~}
  \psi'(N)=\frac{N^4}{\psi(N)^2}\;,\vspace{.5mm}\\ 
  \frac{1}{\covol_{\vec\Lambda}\,|z|^4}
\sum_{p\in\vec\Lambda\,:\,|p|\leq \frac{|z|}{\lambda_\psi}}\;\;|p|^2 &{\rm
    ~~if~~}\lambda_\psi\neq 0,+\infty{\rm ~and~} \psi'(N)=\psi(N)^2\;.
\end{cases}
\end{equation}
The convergence 
\begin{equation}\label{eq:logpaircorrelscal}
\frac{1}{\psi'(N)}\;\R^{\L_\Lambda,\psi}_N\;\;\;\weakstar\;\;\;
g_{\L_\Lambda,\psi}\;\Leb_\CC\;,
\end{equation}
is uniform on every compact subset of $\ZZ$-grids $\Lambda$ in the
space $\operatorname{Grid}_2$.

Furthermore, if $\lambda_\psi\neq 0,+\infty$, for all $A\geq 1$ and
$f\in C^1_{\rm c} (\CC)$ with support contained in $B(0,A)$, we have
\begin{align*}
&\frac{1}{\psi'(N)}\;\R^{\L_\Lambda,\psi}_N(f_N)=
\int_{z\in\CC}f(z)\;g_{\L_\Lambda,\psi}(z)\;d\Leb_\CC(z)\\ &\;\;\;\;+
\bigO\Big(\;\frac{A^5\;\|f\|_\infty\,\big|\,\lambda_\psi -
  \frac{\psi(N)}{N}\,\big|}{\lambda_\psi^9\,
  \operatorname{Sys}_{\vec\Lambda}^4\,\covol_{\vec\Lambda}^2}+
\frac{A^4\,\diam_{\vec\Lambda}\,\|df\|_\infty}
     {\lambda_\psi^4\covol_{\vec\Lambda}^2\,
       \operatorname{Sys}_{\vec\Lambda}\,\psi(N)}
+\frac{A^2(\diam_{\vec\Lambda} +\frac{A}{\lambda_\psi})\,
\|f\|_\infty}{\lambda_\psi^3\covol_{\vec\Lambda}^2\,\psi(N)}\;\Big)\;.
\end{align*}
\etheo

Note that the pair correlation function $g_{\L_\Lambda,\psi}$ depends
on $\vec\Lambda$ but is independent of $a$.  The above result shows in
particular that renormalizing to probability measures (taking
$\psi'(N)\sim \frac{\pi^2N^4}{\covol_{\vec\Lambda}^2}$) is
inappropriate, as the limiting measure would always be $0$. We will
see during the proof that the above result implies the cases
$\alpha>0$ of Theorem \ref{theo:intro1} in the introduction.

The fact that $g_{\L_\Lambda,\psi}$ vanishes when $\lambda_\psi =
+\infty$ means that the sequence of measures $\big(
\frac{1}{\psi'(N)}\; \R^{\L_\Lambda,\psi}_N \big)_{N\in\NN-\{0\}}$ on
$(E_N)_{N\in\NN-\{0\}}$ has a total loss of mass at infinity. For
error terms when $\lambda_\psi=+\infty$ and $\lambda_\psi= 0$, see
respectively Equation \eqref{eq:errortermlambdaphiinfty} and Equation
\eqref{eq:errortermlambdaphizero}.

\medskip
\dem Let $\Lambda=a+\vec\Lambda$ be a $\ZZ$-grid in $\CC$. We may
assume that $a\in \F_{\vec\Lambda}$. Let $N\in\NN-\{0\}$. Let
\[
E^\pm_N=(\pm[0,\infty[\,+i\RR)/(2\pi i\,\psi(N)\, \ZZ)
\]
(which contains the base point $0$) so that $E_N=E^-_N\cup
E^+_N$. Note that $(E^\pm_N)_{N\in\NN-\{0\}}$ converges for the
pointed Hausdorff-Gromov convergence to the closed halfplane
$\CC^\pm=\pm[0,\infty[\,+i\RR$ and that $\CC^-\cap \CC^+$ has measure
    $0$ for any measure absolutely continuous with respect to the
    Lebesgue measure on $\CC$. Note that if $f\in C^1_{\rm c}
    (\CC^\pm)$, then for $N$ large enough, we have $f_N\in C^1_{\rm c}
    (E^\pm_N)$, with the above notation.

Let $\sg_N:E_N\ra E_N$ be the change of sign map $z'\mapsto -z'$,
which maps $E^-_N$ to $E^+_N$ and converges to the change of sign map
$\sg:z\mapsto -z$ on $\CC$.  The change of variables $(m,n)\mapsto
(n,m)$ in the index set $I_N$ proves that we have
$\R^{\L_\Lambda,\psi}_N \!\mid_{E^-_N}=(\sg_N)_* \big(
\R^{\L_\Lambda,\psi}_N \!  \mid_{E^+_N}\big)$. We will thus only study
the convergence of the measures $\frac{1}{\psi'(N)}
\;\R^{\L_\Lambda,\psi}_N$ on $E^+_N$, and deduce the global result by
the symmetry of $g_{\L_\Lambda,\psi}$ under $\sg$.

For every $p\in\vec\Lambda-\{0\}$, let 
\begin{equation}\label{eq:defiJsubpN}
J_{p,\,N}=\{q\in\Lambda :  0< |q|\leq |p+q|\leq N\}\;,
\end{equation}
and let
\begin{equation}\label{eq:defiomegapN}
\omega_{p,\,N}=\sum_{q\in J_{p,\,N}}\Delta_{\psi(N)\frac{p}{q}}
      {\rm ~~~and~~~}
\mu_N^+=\sum_{p\in\vec\Lambda-\{0\}}\omega_{p,\,N}\;.
\end{equation}
Note that $\omega_{p,\,N}$ is a measure on $\CC$ with finite support,
which vanishes if $|p|>2N$ by the triangle inequality, hence $\mu_N^+$
is also a measure on $\CC$ with finite support.

\blemm\label{lem:massomegapN} As $N\ra+\infty$, we have
$\|\omega_{p,\,N}\|=\card\; J_{p,\,N}=\frac{\pi N^2}{2\covol_{\vec\Lambda}}+
  \bigO\Big(\frac{(|p|+\diam_{\vec\Lambda})\,N}{\covol_{\vec\Lambda}}\Big)$.
\elemm

\dem We may assume that $|p|\leq 2N$. Note that $J_{p,\,N}$ is the
finite set of nonzero grid points in the intersection
\begin{equation}\label{eq:defCpN}
\wt C_{p,\,N}=\{z\in\CC:|z|\leq |p+z|\leq N\}
\end{equation}
of the disc $B(-p,N)$ of radius $N$ centered at $-p$ with the closed
halfplane containing $0$ with boundary the 
perpendicular bisector of $0$ and
$-p$ (see the picture below).
\begin{center}
\begin{picture}(0,0)%
\includegraphics{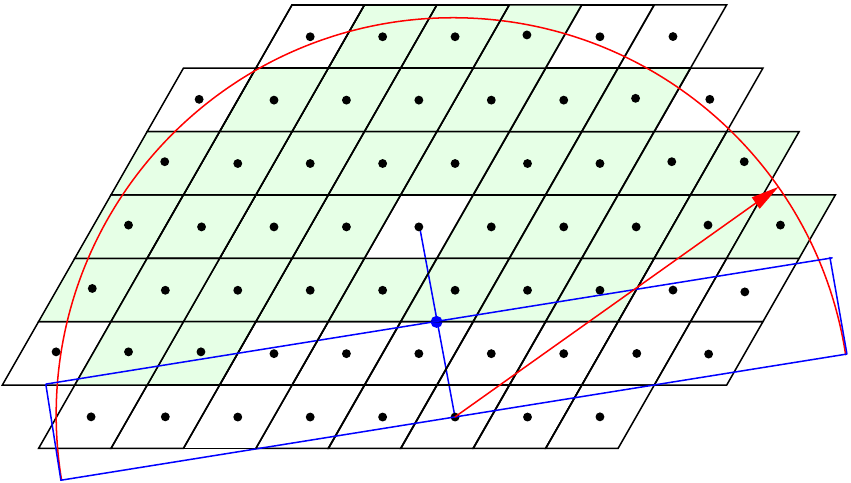}%
\end{picture}%
\setlength{\unitlength}{3812sp}%
\begingroup\makeatletter\ifx\SetFigFont\undefined%
\gdef\SetFigFont#1#2#3#4#5{%
  \reset@font\fontsize{#1}{#2pt}%
  \fontfamily{#3}\fontseries{#4}\fontshape{#5}%
  \selectfont}%
\fi\endgroup%
\begin{picture}(4221,2387)(-731,-1941)
\put(1390,-1755){\makebox(0,0)[lb]{\smash{{\SetFigFont{11}{13.2}{\rmdefault}{\mddefault}{\updefault}{\color[rgb]{0,0,0}$-p$}%
}}}}
\put(1228,-781){\makebox(0,0)[lb]{\smash{{\SetFigFont{11}{13.2}{\rmdefault}{\mddefault}{\updefault}{\color[rgb]{0,0,0}$0$}%
}}}}
\put(3181,-462){\makebox(0,0)[lb]{\smash{{\SetFigFont{11}{13.2}{\rmdefault}{\mddefault}{\updefault}{\color[rgb]{1,0,0}$N$}%
}}}}
\put(1323,-1301){\makebox(0,0)[lb]{\smash{{\SetFigFont{11}{13.2}{\rmdefault}{\mddefault}{\updefault}{\color[rgb]{0,0,1}$-\frac{p}{2}$}%
}}}}
\end{picture}%

\end{center}

Since $\wt C_{p,\,N}$ is contained in a halfdisc of radius $N$ and
contains the complement in this halfdisc of its intersection with a
rectangle of length $2N$ and height $\frac{|p|}{2}$, we have
$\frac{\pi}{2}N^2-|p|\,N\leq \operatorname{Area}(\wt C_{p,\,N})\leq
\frac{\pi}{2}N^2$, so that
\[
\operatorname{Area}(\wt C_{p,\,N})= \frac{\pi}{2}N^2+\bigO(|p|\,N)\;.
\]
Let
\begin{equation}\label{eq:defiCpN}
  C_{p,\,N}=\bigcup_{q\in J_{p,\,N}} (q+\F_{\vec\Lambda})\;.
\end{equation}
By a Gauss counting argument similar to the one in the proof of
Equation \eqref{eq:gaussRquat} with $k=0$, we have
\begin{align*}
  \|\omega_{p,\,N}\|&=\card\; J_{p,\,N}=
  \frac{\operatorname{Area}(C_{p,\,N})}{\covol_{\vec\Lambda}}
  =\frac{\operatorname{Area}(\wt C_{p,\,N})}{\covol_{\vec\Lambda}}+
  \frac{\operatorname{Area}(C_{p,\,N})-
    \operatorname{Area}(\wt C_{p,\,N})}{\covol_{\vec\Lambda}}\\ &=
  \frac{\pi N^2}{2\covol_{\vec\Lambda}}+
  \bigO\Big(\frac{(|p|+\diam_{\vec\Lambda})\,N}{\covol_{\vec\Lambda}}\Big)\;.
\end{align*}
The lemma follows.
\cqfd

\blemm\label{lem:relatRetmu} For every $A>0$ and for every $f\in
C^1_{\rm c}(\CC^+)$ with support contained in $B(0,A)$,  as
$N\ra+\infty$ and uniformly on $\Lambda$ varying in a compact subset
of $\operatorname{Grid}_2$, we have
$$
\big|\;(\R^{\L_\Lambda,\psi}_N)_{\mid E^+_N}(f_N)-\mu_N^+(f)\,\big|=
\bigO\Big(
\frac{A^4\;\|df\|_\infty\;N^4}{\covol_{\vec\Lambda}^{\;2}\,\psi(N)^3}\Big)\;.
$$
\elemm

\dem Let $A$ and $f$ be as in the statement of this lemma. Note that
since $\psi(N)>0$ and by Equation \eqref{eq:relatlognormRenlog}, for
every $(m,n)\in I_N$, we have $(m,n)\in I_N^+$, that is $|n|\leq |m|$,
if and only if $\psi(N)\log m-\psi(N)\log n\in E_N^+$. Hence by the
change of variable $(p,q)\mapsto (m=p+q,n=q)$ (which is a bijection
from $\vec\Lambda\times\Lambda$ to $\Lambda\times\Lambda$), we have
\begin{align*}
(\R^{\L_\Lambda,\psi}_N)_{\mid E^+_N}(f_N)
&=\sum_{(m,\,n)\in I_N^+}\;\;f_N(\psi(N)\log m-\psi(N)\log n)
\\&=\sum_{\substack{p\in\vec\Lambda-\{0\},\,q\in\Lambda\\
    0<|q|\leq |p+q|\leq N}}
\;\;f_N\big(\psi(N)\log (p+q)-\psi(N)\log q\big)\;.
\end{align*}
By the assumption on the support of $f$, if an index $(p,q)$
contributes to the above sum, then $\Re(\psi(N)\log (p+q)-\psi(N)\log
q)\leq A$. Hence by Equations \eqref{eq:RedilatRe} and
\eqref{eq:relatlognormRenlog}, we have $\ln\big|1+\frac{p}{q}\big|\leq
\frac{A}{\psi(N)}$, which tends to $0$ as $N\ra+\infty$, since
${\displaystyle \lim_{+\infty}\;\psi=+\infty}$.  In particular, using
the assumption on $q$, we have
\begin{equation} \label{eq:controlabsp}
  \frac{|p|}{|q|}=\bigO\big(\frac{A}{\psi(N)}\big) \;\;\;{\rm and}\;\;\;
    |p|=\bigO\big(\frac{AN}{\psi(N)}\big)\;,
\end{equation}
so that $\big|\frac{p}{q}\big|<1$ if $N$ is large enough. This allows
to use the principal branch, again denoted by $\log$, of the
complex logarithm in the open ball of center $1$ and radius $1$.  By
the analytic expansion of this branch, we have
\[
\Big|\log(1+\frac{p}{q})-\frac{p}{q}\,\Big|=
\bigO\big(\big|\frac{p}{q}\big|^2\big)=
\bigO\big(\frac{A^2}{\psi(N)^2}\big)\;.
\]
The mean value theorem hence implies that
\begin{align} 
f_N(\psi(N)\log (p+q)-\psi(N)\log q)&=
f\big(\psi(N)\log(1+\frac{p}{q})\big)\nonumber\\
&=f(\psi(N)\frac{p}{q})+\bigO\big(\frac{A^2\|df\|_\infty}{\psi(N)}\big)\;.
\label{eq:meanvalscal}
\end{align}

By Lemma \ref{lem:massomegapN} and Equation \eqref{eq:gaussRquat} with
$k=0$, we have
\begin{align}
  &\card\big\{(p,q)\in\vec\Lambda\times\Lambda\;:\;0<|q|\leq |p+q|\leq N,\;
  |p|=\bigO\big(\frac{AN}{\psi(N)}\big)\big\} \nonumber\\= \;&
  \sum_{p\in\vec\Lambda-\{0\}\;:\;|p|=\bigO\big(\frac{AN}{\psi(N)}\big)}\;
  \card\; J_{p,\,N}=
  \bigO\Big(\card\big\{p\in\vec\Lambda-\{0\}\;:\;|p|=
  \bigO\big(\frac{AN}{\psi(N)}\big)\big\}\,\frac{N^2}{\covol_{\vec\Lambda}}\Big)
  \nonumber\\
  = \; &\bigO\Big(\frac{A^2N^4}{\psi(N)^2\covol_{\vec\Lambda}^2}
  \Big)\;.\label{eq:cardpqbigO}
\end{align}
Similarly, if an index $(p,q)$ contributes to the sum
\[
\mu_N^+(f)=\sum_{\substack{p\in\vec\Lambda-\{0\},\,q\in\Lambda\\
    0<|q|\leq |p+q|\leq N}}
\;\;f\big(\psi(N)\frac{p}{q}\big)\;,
\]
then Equation \eqref{eq:controlabsp} holds. By summing Equation
\eqref{eq:meanvalscal} on the set of elements $(p,q)\in\vec\Lambda
\times\Lambda$ such that $0<|q|\leq |p+q|\leq N$ and $|p|= \bigO
\big(\frac{AN}{\psi(N)}\big)$, and by using Equation \eqref{eq:cardpqbigO},
Lemma \ref{lem:relatRetmu} follows.  \cqfd

\medskip
Let us now study the convergence properties (after renormalization) of
the measures $\omega_{p,\,N}$ and of their sums $\mu^+_N$ as
$N\ra+\infty$. We assume in what follows that $|p|<N$ (which is
possible if $N$ is large enough since we will have
$|p|=\bigO\big(\frac{AN}{\psi(N)}\big)$~). Let $\iota:\CC^\times \ra
\CC^\times$ be the involutive diffeomorphism $z\mapsto \frac{1}{z}$,
which maps $\CC^+-\{0\}$ to $\CC^+-\{0\}$, whose holomorphic
derivative at $z$ is $\frac{1}{z^2}$, hence whose Jacobian at $z$ is
\begin{equation}\label{eq:jaciota}
J\iota(z)=\frac{1}{|z|^4}\;,
\end{equation}
By the equation on the left in Formula \eqref{eq:defiomegapN}, we have
\begin{equation}\label{eq:iotaomegapN}
\iota_*\omega_{p,\,N}=\sum_{q\,\in J_{p,\,N}}\Delta_{\frac{q}{\psi(N)p}}\;.
\end{equation}
When $q$ varies in $J_{p,\,N}$, as seen in the proof of Lemma
\ref{lem:massomegapN}, the above Dirac masses are exactly at the nonzero
points of the $\ZZ$-grid $\Lambda_{p,N}=\frac{1}{\psi(N)p}\,\Lambda$
that belong to the set
\[
\wt Y_{p,N} =\frac{1}{\psi(N)\,p}\,\wt C_{p,N}\;.
\]

\medskip
\noindent \begin{minipage}{8.2cm} ~~~ 
Note that
\begin{equation}\label{eq:covolLambdapN}
  \covol_{\vec\Lambda_{p,N}}=
  \operatorname{Area}\big(\frac{\F_{\vec\Lambda}}{\psi(N)\,p}\big)=
  \frac{\covol_{\vec\Lambda}}{\psi(N)^2\,|p|^2}\;.
\end{equation}
By Equation \eqref{eq:defCpN}, the set $\wt Y_{p,N}$ is the
intersection of the disc $B(-\frac{1}{\psi(N)},\frac{N}{\psi(N)|p|})$
with the closed halfplane containing $0$ with boundary the 
perpendicular bisector 
of $0$ and $-\frac{1}{\psi(N)}$. Let us define
\[
Z_{p,N} =\big\{z\in\CC:\;\Re \;z\geq 0,\;\;
|z|\leq \frac{N}{\psi(N)|p|}\big\}\;.
\]
Note that
\begin{equation}\label{eq:iotaZpN}
\iota(Z_{p,N}) =\big\{z\in\CC:\;\Re \;z\geq 0,\;\;
|z|\geq \frac{\psi(N)|p|}{N}\big\}\;.
\end{equation}
\end{minipage}
\begin{minipage}{6.7cm}
\begin{picture}(0,0)%
\includegraphics{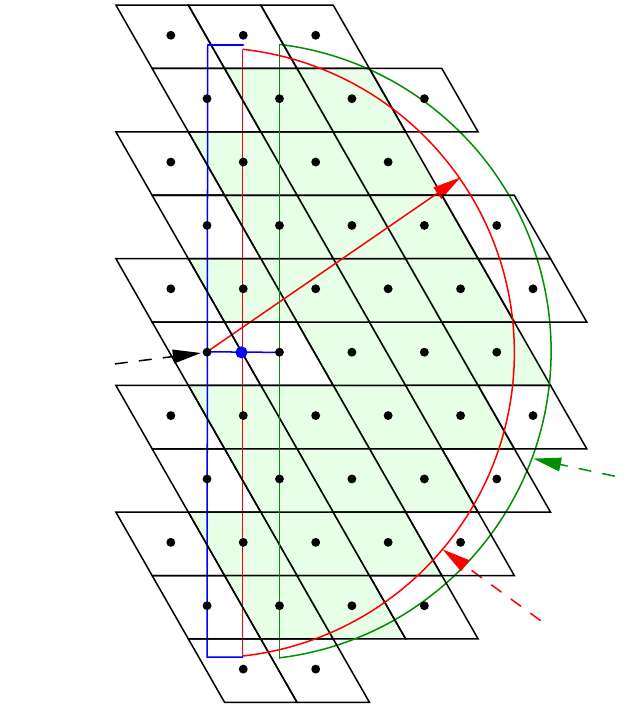}%
\end{picture}%
\setlength{\unitlength}{3812sp}%
\begingroup\makeatletter\ifx\SetFigFont\undefined%
\gdef\SetFigFont#1#2#3#4#5{%
  \reset@font\fontsize{#1}{#2pt}%
  \fontfamily{#3}\fontseries{#4}\fontshape{#5}%
  \selectfont}%
\fi\endgroup%
\begin{picture}(3117,3489)(5680,-2458)
\put(8782,-1391){\makebox(0,0)[lb]{\smash{{\SetFigFont{11}{13.2}{\rmdefault}{\mddefault}{\updefault}{\color[rgb]{0,.56,0}$Z_{p,N}$}%
}}}}
\put(8413,-2102){\makebox(0,0)[lb]{\smash{{\SetFigFont{11}{13.2}{\rmdefault}{\mddefault}{\updefault}{\color[rgb]{1,0,0}$\wt Y_{p,N}$}%
}}}}
\put(8011,164){\makebox(0,0)[lb]{\smash{{\SetFigFont{11}{13.2}{\rmdefault}{\mddefault}{\updefault}{\color[rgb]{1,0,0}$N/(\psi(N)|p|)$}%
}}}}
\put(7111,-826){\makebox(0,0)[lb]{\smash{{\SetFigFont{11}{13.2}{\rmdefault}{\mddefault}{\updefault}{\color[rgb]{0,0,0}$0$}%
}}}}
\put(5695,-793){\makebox(0,0)[lb]{\smash{{\SetFigFont{11}{13.2}{\rmdefault}{\mddefault}{\updefault}{\color[rgb]{0,0,0}$-\frac{1}{\psi(N)}$}%
}}}}
\end{picture}%

\end{minipage}

\medskip
The symmetric difference of $\wt Y_{p,N}$ and $Z_{p,N}$, that we
denote by $\wt Y\!Z_{p,N}$, is contained in the union of the rectangle
$\big[-\frac{1}{2\psi(N)},0\big]\times \big[-\frac{N}{\psi(N)|p|},
  \frac{N}{\psi(N)|p|}\big]$ and the half-annulus
\[
\Big\{z\in\CC:\;\Re \;z\geq 0,\;\;\frac{N}{\psi(N)|p|}
-\frac{1}{\psi(N)}\leq |z|\leq \frac{N}{\psi(N)|p|}\Big\}
\]
(well defined since $|p|<N$). In particular, its area satisfies 
$
\Leb_\CC(\wt Y\!Z_{p,N})=\bigO\big(\frac{N}{\psi(N)^2|p|}\big)$.
Let
\begin{equation}\label{eq:defiYpN}
Y_{p,N} =\frac{1}{\psi(N)\,p}\,C_{p,N}\;,
\end{equation}
so that, as in the proof of Lemma \ref{lem:massomegapN}, the symmetric
difference of $Y_{p,N}$ and $\wt Y_{p,N}$ has area $\bigO\big(
\frac{N\diam_{\vec\Lambda}}{\psi(N)^2|p|^2}\big)$.  The symmetric
difference of $Y_{p,N}$ and $Z_{p,N}$, that we denote by $Y\!Z_{p,N}$,
hence has area $ \Leb_\CC(Y\!Z_{p,N})= \bigO\big(
\frac{N(\diam_{\vec\Lambda} +|p|)}{\psi(N)^2|p|^2}\big)$.  In
particular, for every $\phi\in C^1_{\rm c}(\CC^+-\{0\})$, since
$Z_{p,N}\subset B(0,\frac{N}{\psi(N)|p|})$ and $Y_{p,N}\subset
B(0,\frac{N+\diam_{\vec\Lambda}}{\psi(N)|p|})$, we have
\begin{align}
  & \Big|\;\int_{Z_{p,N}} \phi(z)\;d\Leb_\CC(z)-\int_{Y_{p,N}} \phi(z)
  \;d\Leb_\CC(z)\;\Big|\nonumber\\= \;\;& \bigO\Big(\;
  \frac{N(\diam_{\vec\Lambda} +|p|)\,
    \big\|\phi_{\mid B(0,\frac{N+\diam_{\vec\Lambda}}{\psi(N)|p|})}\big\|_\infty}
       {\psi(N)^2|p|^2}\;\Big)\;.\label{eq:relatYpNZpN}
\end{align}
By Equations \eqref{eq:defiYpN}, \eqref{eq:defiCpN},
\eqref{eq:iotaomegapN} and \eqref{eq:covolLambdapN}, by the mean value
theorem and by Lemma \ref{lem:massomegapN}, we have
\begin{align*}
  &\Big|\;\int_{Y_{p,N}} \phi(z)\;d\Leb_\CC(z)\;-\;
  \frac{\covol_{\vec\Lambda}}{\psi(N)^2\,|p|^2}\;\iota_*\omega_{p,\,N}(\phi)
  \;\Big|\\=\;\; &
  \Big|\;\sum_{q\in J_{p,N}}\int_{\frac{q+ \F_{\vec\Lambda}}{\psi(N)\,p}}
  \Big(\phi(z)-\phi\big(\frac{q}{\psi(N)\,p}\big)\Big)\;d\Leb_\CC(z)\;\Big|
  \\\leq\;\; &(\card\; J_{p,\,N})\;
  \frac{\covol_{\vec\Lambda}}{\psi(N)^2\,|p|^2}\;
  \Big\|d\phi_{\mid B(0,\frac{N+\diam_{\vec\Lambda}}{\psi(N)|p|})}\Big\|_\infty\;
  \frac{\diam_{\vec\Lambda}}{\psi(N)\,|p|}\\=\;\; &
  \bigO\Big(\;\frac{\diam_{\vec\Lambda}\;
    \Big\|d\phi_{\mid B(0,\frac{N+\diam_{\vec\Lambda}}{\psi(N)|p|})}\Big\|_\infty\,N^2}
           {\psi(N)^3\,|p|^3}\;\Big)\;.
\end{align*}
Hence by Equation \eqref{eq:relatYpNZpN}, we have
\begin{align}
&\iota_*\omega_{p,\,N}(\phi) =\frac{\psi(N)^2\,|p|^2}{\covol_{\vec\Lambda}}
\;\int_{Z_{p,N}} \phi(z)\;d\Leb_\CC(z)\nonumber\\&\;\;+
\bigO\Big(\;\frac{\diam_{\vec\Lambda}\,
\Big\|d\phi_{\mid B(0,\frac{N+\diam_{\vec\Lambda}}{\psi(N)|p|})}\Big\|_\infty\,N^2}
{\covol_{\vec\Lambda}\,\psi(N)\,|p|}
+\frac{(\diam_{\vec\Lambda} +|p|)\,\Big\|
\phi_{\mid B(0,\frac{N+\diam_{\vec\Lambda}}{\psi(N)|p|})}\Big\|_\infty\,N}
{\covol_{\vec\Lambda}}\;\Big)\;.\label{eq:iotaomegapNbis}
\end{align}
Let $f\in C^1_{\rm c}(\CC^+-\{0\})$ with support contained in
$B(0,A)$. Note that $f\circ\iota\in C^1_{\rm c}(\CC^+-\{0\})$, that
$\Big\|f\circ\iota_{\mid \{|z|\leq \frac{N+\diam_{\vec\Lambda}}{\psi(N)|p|}\}}
\Big\|_\infty=\Big\|
f_{\mid\{|z|\geq \frac{\psi(N)|p|}{N+\diam_{\vec\Lambda}}\}}\Big\|_\infty$ and that
\[
\Big\|\,d(f\circ\iota)_{\mid \{|z|\leq \frac{N+\diam_{\vec\Lambda}}{\psi(N)|p|}\}}
\Big\|_\infty\leq A^2 \;\big\|\,
df_{\mid\{|z|\geq \frac{\psi(N)|p|}{N+\diam_{\vec\Lambda}}\}}\big\|_\infty
\]
since the support of $f$ is contained in $B(0,A)$. The change of
variable by $\iota$ in the integral of Equation
\eqref{eq:iotaomegapNbis} applied with $\phi=f\circ \iota$, together
with Equations \eqref{eq:iotaZpN} and \eqref{eq:jaciota}, hence give
\begin{align*}
&\omega_{p,\,N}(f) =\frac{\psi(N)^2\,|p|^2}{\covol_{\vec\Lambda}}
  \;\int_{|z|\geq \frac{\psi(N)|p|}{N}} f(z)\;\frac{1}{|z|^4}\;
  d\Leb_\CC(z)\\&\;\;+\bigO\Big(\;\frac{A^2\,\diam_{\vec\Lambda}\,\big\|
  df_{\mid\{|z|\geq \frac{\psi(N)|p|}{N+\diam_{\vec\Lambda}}\}}\big\|_\infty\,N^2}
{\covol_{\vec\Lambda}\,\psi(N)\,|p|}
+\frac{(\diam_{\vec\Lambda} +|p|)\,\big\|
f_{\mid\{|z|\geq \frac{\psi(N)|p|}{N+\diam_{\vec\Lambda}}\}}\big\|_\infty\,N}
{\covol_{\vec\Lambda}}\;\Big)\;.
\end{align*}

For every $z\in \CC^+-\{0\}$, let
\begin{equation}\label{eq:defithetaN}
\theta_N(z)=\frac{1}{|z|^4}\sum_{p\in\vec\Lambda-\{0\}}\;\;|p|^2\;
\mathbbm{1}_{\big\{|z|\geq \frac{\psi(N)|p|}{N}\big\}}(z)
=\frac{1}{|z|^4}\sum_{p\in\vec\Lambda-\{0\}\,:\,|p|\leq
  \frac{N|z|}{\psi(N)}}\;\;|p|^2\,.
\end{equation}
Note that if $z$ and $N$ are fixed, then for $|p|$ large enough, we
have $|z|< \frac{\psi(N)|p|}{N}$, thus the above sum has
only finitely many nonzero terms. Let $\theta_N(0)=0$.

Note that $\theta_N(z)$ vanishes if and only if $|z|< \frac{\psi(N)
  \operatorname{Sys}_{\vec\Lambda}}{N}$, by the definition of the
systole of $\vec\Lambda$.

As seen in the proof of Lemma \ref{lem:relatRetmu}, the only elements
$p\in\vec\Lambda$ that give a nonzero contribution to the sum
$\sum_{p\in\vec\Lambda-\{0\}}\omega_{p,\,N}(f)$ satisfy $p\neq 0$
%
%
and $|p|=\bigO\big(\frac{AN}{\psi(N)}\big)$. By Equation
\eqref{eq:gaussRquat} with $k=0$, we have
\[
\card\big\{\,p\in\vec\Lambda-\{0\}:|p|=
\bigO\big(\frac{AN}{\psi(N)}\big) \big\}
=\bigO\Big(\frac{A^2N^2}{\covol_{\vec\Lambda}\psi(N)^2}\Big)\;
\]
if
%
%
$\lambda_\psi<+\infty$. Otherwise, if $\lambda_\psi=+\infty$, we have
$\bigO\big(\frac{AN}{\psi(N)}\big)\leq
\operatorname{Sys}_{\vec\Lambda}$ if $N$ is large enough, hence if $N$
is large enough, we have
\begin{equation}\label{eq:case1prepa}
  \card\big\{\,p\in\vec\Lambda-\{0\}:|p|=
  \bigO\big(\frac{AN}{\psi(N)}\big) \big\}=0\;.
\end{equation}
Thus, by the right equality in Formula \eqref{eq:defiomegapN}, we have 
{\small \begin{align}
  &\mu_N^+(f)=\sum_{p\in\vec\Lambda-\{0\}}\omega_{p,\,N}(f)=
  \frac{\psi(N)^2}{\covol_{\vec\Lambda}}
  \;\int_{z\in\CC^+} f(z)\;\theta_N(z)\;d\Leb_\CC(z)
  \nonumber\\&\;\;+\bigO\Big(\;\frac{A^4\,\diam_{\vec\Lambda}\,\big\|
df_{\mid\{|z|\geq \frac{\psi(N)\,\operatorname{Sys}_{\vec\Lambda}}{N+\diam_{\vec\Lambda}}\}}
\big\|_\infty\,N^4}
  {\covol_{\vec\Lambda}^2\,\operatorname{Sys}_{\vec\Lambda}\,\psi(N)^3}
+\frac{(\diam_{\vec\Lambda} +\frac{AN}{\psi(N)})\,\big\|
  f_{\mid\{|z|\geq \frac{\psi(N)\,\operatorname{Sys}_{\vec\Lambda}}{N+\diam_{\vec\Lambda}}\}}
  \big\|_\infty\,A^2\,N^3}{\covol_{\vec\Lambda}^2\,\psi(N)^2}\;\Big)\;.
\label{eq:equation5dePCPC}
\end{align}}

\noindent{\bf Case 1. } Let us first assume that $\lambda_\psi
=+\infty$, that is, ${\displaystyle \lim_{N\ra+\infty}}
\frac{N}{\psi(N)}=0$.

For every $A\geq 1$, if $N$ is large enough (uniformally on $\Lambda$
varying in a compact subspace of $\operatorname{Grid}_2$, since then
$\vec\Lambda$ varies in a compact subspace of the space of
$\ZZ$-lattices, on which the systole function $\vec\Lambda\mapsto
\operatorname{Sys}_{\vec\Lambda}$ has a positive lower bound), then
for every $z\in B(0,A)$, we have $\theta_N(z)=0$ by Equation
\eqref{eq:defithetaN}, and $\mu^+_N(f)=0$ by Formulas
\eqref{eq:defiomegapN} and \eqref{eq:case1prepa}, since the sum
defining $\mu^+_N(f)$ is an empty sum.  Thus, whatever the (positive)
normalizing function $\psi'$ is, we have a total loss of mass at
infinity~:
$$
\frac{1}{\psi'(N)}\;\mu^+_N\;\;\weakstar\;\; 0\;.
$$

Assume that the renormalizing function $\psi'$ is such that
$\frac{N^4} {\psi(N)^3\psi'(N)}$ tends to $0$ as $N$ tends to
$\infty$, for instance $\psi'=\psi$, as assumed in the first case of
Equation \eqref{eq:troiscas}. Note that if $\psi(N)=N^\alpha$ with
$\alpha >1$, then we indeed have $\lambda_\psi =+\infty$ and if
$\psi'(N) =N^{4-2\alpha}$ as in the statement of Theorem
\ref{theo:intro1}, we do have $\lim_{N\ra+\infty}
\frac{N^4}{\psi(N)^3\psi'(N)}=0$.

Together with Lemma \ref{lem:relatRetmu}, the above centered formula
proves Formula \eqref{eq:logpaircorrelscal} when $\lambda_\psi=
+\infty$, with a convergence which is uniform on every compact subset
of $\Lambda$ in $\operatorname{Grid}_2$, as well as the case $\alpha
>1$ in Theorem \ref{theo:intro1}. Furthermore, if follows from the
error term in Lemma \ref{lem:relatRetmu} that for every $f\in C^1_{\rm
  c}(\CC)$ with support contained in $B(0,A)$, as $N\ra+\infty$ and
uniformly on $\Lambda$ varying in a compact subset of
$\operatorname{Grid}_2$, we have
\begin{equation}\label{eq:errortermlambdaphiinfty}
\frac{1}{\psi'(N)}\;\R^{\L_\Lambda,\psi}_N(f_N)=\bigO\Big(\;
\frac{A^4\;\|df\|_\infty\;N^4}{\covol_{\vec\Lambda}^{\;2}\,\psi(N)^3\psi'(N)}
\;\Big)\;.
\end{equation}

\medskip
\noindent{\bf Case 2. } Let us now assume that $\lambda_\psi= 0$, that
is, ${\displaystyle \lim_{N\ra+\infty}} \frac{\psi(N)}{N}=0$. 

For all $z\in \CC^+-\{0\}$, by Equations \eqref{eq:defithetaN} and
\eqref{eq:gaussRquat} for $k=2$, we have
\begin{equation}\label{eq:controlONt}
\frac{\psi(N)^4}{N^4}\,\theta_N(z)= \big(\frac{\psi(N)}{N\,|z|}\big)^4
\sum_{p\in\vec\Lambda:\,|p|\leq \frac{N\,|z|}{\psi(N)}}\;|p|^2
=\frac{\pi}{2\covol_{\vec\Lambda}}+
\bigO\big(\,\frac{(1+\diam_{\vec\Lambda})\,\psi(N)}
{\covol_{\vec\Lambda}\,N\,|z|}\big)\;.
\end{equation}
In particular, if $|z|\geq \frac{\psi(N)\operatorname{Sys}_{\vec\Lambda}}
{N}$, then $\frac{\psi(N)^4}{N^4}\,\theta_N(z)$ is uniformly
bounded. Since $\theta_N(z)$ vanishes if $|z|< \frac{\psi(N)
  \operatorname{Sys}_{\vec\Lambda}}{N}$, this proves that the function
$\frac{\psi(N)^4}{N^4}\,\theta_N$ is uniformly bounded on
$\CC^+-\{0\}$, and pointwise converges to the constant function
$\frac{\pi}{2\covol_{\vec\Lambda}}$. Hence by Equation
\eqref{eq:equation5dePCPC} and by the Lebesgue dominated convergence
theorem, we have, with a convergence which is uniform on every compact
subset of $\Lambda$ in $\operatorname{Grid}_2$,
\begin{equation}\label{eq:convergmupluscas2}
\frac{\psi(N)^2}{N^4}\;\mu^+_N\;\;\weakstar\;\;
\frac{\pi}{2\covol_{\vec\Lambda}^2}\;\Leb_{\CC^+}\;.
\end{equation}

More precisely, for every $A\geq 1$, for every $f\in C^1_{\rm
  c}(\CC^+-\{0\})$ with support in $B(0,A)$, and for every $\Lambda$
in a compact subset of $\operatorname{Grid}_2$, by Equations
\eqref{eq:equation5dePCPC} and \eqref{eq:controlONt}, using the
equality ${\displaystyle
  \int_{-\pi/2}^{\pi/2}\int_0^A\frac{1}{\rho}\,\rho\,
  d\rho\,d\theta=\pi A}$ in order to integrate the error term in
Equation \eqref{eq:controlONt}, and since $\psi(N)\leq N$ for $N$ large
enough, we have
\begin{align*}
  \frac{\psi(N)^2}{N^4}\;\mu_N^+(f) &=
  \frac{\psi(N)^4}{\covol_{\vec\Lambda}\,N^4}
  \;\int_{z\in\CC^+} f(z)\;\theta_N(z)\;d\Leb_\CC(z)
\\&\;\;\;+\bigO\Big(\;\frac{A^4\,\diam_{\vec\Lambda}\,\|df\|_\infty}
  {\covol_{\vec\Lambda}^2\,\operatorname{Sys}_{\vec\Lambda}\,\psi(N)}
  +\frac{A^2(\diam_{\vec\Lambda} +A)\,
    \|f\|_\infty}{\covol_{\vec\Lambda}^2\,\psi(N)}\;\Big)
\\&=\frac{\pi}{2\covol_{\vec\Lambda}^2}
  \;\int_{z\in\CC^+} f(z)\;d\Leb_\CC(z)+
  \bigO\Big(\;\frac{A^4\,\diam_{\vec\Lambda}\,\|df\|_\infty}
  {\covol_{\vec\Lambda}^2\,\operatorname{Sys}_{\vec\Lambda}\,\psi(N)}
\\&\;\;\;\;\;\;\;\;\;\;\;\;\;\;\; +\frac{A^2(\diam_{\vec\Lambda} +A)\,
    \|f\|_\infty}{\covol_{\vec\Lambda}^2\,\psi(N)}+
  \frac{A\,(1+\diam_{\vec\Lambda})\,\psi(N)\|f\|_\infty}
       {\covol_{\vec\Lambda}^2\,N}\;\Big)\;.
\end{align*}

If $\psi'(N)=\frac{N^4}{\psi(N)^2}$ as assumed in the second case of
Equation \eqref{eq:troiscas}, it follows from Formula
\eqref{eq:convergmupluscas2} and Lemma \ref{lem:relatRetmu} by
symmetry that
\[
\frac{1}{\psi'(N)}\;\R^{\L_\Lambda,\psi}_N\;\;\weakstar\;\;
\frac{\pi}{2\covol_{\vec\Lambda}^2}\;\Leb_{\CC}\;.
\]
This proves Formula \eqref{eq:logpaircorrelscal} when $\lambda_\psi=
0$, with a convergence which is uniform on every compact subset of
$\Lambda$ in $\operatorname{Grid}_2$, as well as the case $0<\alpha<1$
in Theorem \ref{theo:intro1}. Furthermore, for every $f\in C^1_{\rm
  c}(\CC)$ with support contained in $B(0,A)$, as $N\ra+\infty$ and
uniformly on $\Lambda$ varying in a compact subset of
$\operatorname{Grid}_2$, using the error term in Lemma
\ref{lem:relatRetmu} with the fact that that
$\operatorname{Sys}_{\vec\Lambda}\leq \diam_{\vec\Lambda}$, we have
\begin{align}
\frac{1}{\psi'(N)}\;\R^{\L_\Lambda,\psi}_N(f_N)&=
\frac{\pi}{2\covol_{\vec\Lambda}^2}\int_{\CC} f\;d\Leb_{\CC}+
\bigO\Big(\;\frac{A^4\,\diam_{\vec\Lambda}\,\|df\|_\infty}
{\covol_{\vec\Lambda}^2\,\operatorname{Sys}_{\vec\Lambda}\,\psi(N)}
\nonumber\\&\;\;\;\;\;
+\frac{A^2(\diam_{\vec\Lambda} +A)\,\|f\|_\infty}
{\covol_{\vec\Lambda}^2\,\psi(N)}+
\frac{A\,(1+\diam_{\vec\Lambda})\,\psi(N)\,\|f\|_\infty}
     {\covol_{\vec\Lambda}^2\,N}
\;\Big)\;. \label{eq:errortermlambdaphizero}
\end{align}

\medskip
\noindent{\bf Case 3. }
Let us finally assume that ${\displaystyle \lim_{N\ra+\infty}}
\frac{\psi(N)}{N}=\lambda_\psi$ belongs to $]0,+\infty[\,$.
    
We consider the function $\theta_\infty:\CC\ra[0,+\infty[$ defined by 
\[
z\mapsto\frac{1}{|z|^4}
\sum_{p\in\vec\Lambda\,:\,|p|\leq \frac{|z|}{\lambda_\psi}}\;\;|p|^2\;, 
\]
where by convention $\theta_\infty(0)=0$, and replacing $p\in
\vec\Lambda$ by $p\in\vec\Lambda-\{0\}$ makes no difference. Note that
$\theta_\infty$ vanishes on the open disc $\stackrel{\circ}{B}
\!\!(0,\operatorname{Sys}_{\vec\Lambda})$, is uniformly bounded and
tends to $\frac{\pi}{2\covol_{\vec\Lambda}}$ as $t\ra+\infty$ by
Equation \eqref{eq:gaussRquat} with $k=0$.  Furthermore,
$\theta_\infty$ is piecewise continuous, with discontinuities along
each circle $S(0,|p|)$ centered at $0$ passing through a nonzero
lattice point $p\in\vec\Lambda$.  See the picture in the introduction
representing the graph of $\theta_\infty$ when $\Lambda= \vec\Lambda
=\ZZ[i]$ and $\lambda_\psi=1$.

By Equation \eqref{eq:defithetaN}, the sequence of uniformly bounded
maps $(\theta_N)_{N\in\NN}$ converges almost everywhere to
$\theta_\infty$ (more precisely, it converges at least outside
$\bigcup_{p\in\vec\Lambda-\{0\}}S(0,|p|)$). Hence by Equation
\eqref{eq:equation5dePCPC} and by the Lebesgue dominated convergence
theorem, we have
\begin{equation}\label{eq:convergmupluscas3}
\frac{1}{\psi(N)^2}\;\mu^+_N\;\;\weakstar\;\;
\frac{1}{\covol_{\vec\Lambda}}\;\theta_\infty\;\Leb_{\CC^+}\;.
\end{equation}

Let $A\geq 1$. Note that $|z|\leq A$ implies that
$\frac{|z|}{\lambda_\psi}\leq \frac{A}{\lambda_\psi}\leq
\frac{2A}{\lambda_\psi}$.  If $N$ is large enough so that
$\frac{\psi(N)}{N}\geq \frac{\lambda_\psi}{2}$, then $|z|\leq A$
implies that $\frac{N|z|}{\psi(N)}\leq \frac{2A}{\lambda_\psi}$.
Hence for every $z\in\CC^+\cap B(0,A)$, if $N$ is large enough, we have
$$
|\;\theta_\infty(z)-\theta_N(z)\,|\leq \frac{1}{|z|^4}
\sum_{p\in\vec\Lambda\,:\, |p|\leq \frac{2A}{\lambda_\psi}} |p|^2\;
\Big|\,\mathbbm{1}_{|z|\geq \lambda_\psi|p|}(z)-
\mathbbm{1}_{|z|\geq \frac{\psi(N)|p|}{N}}(z)\,\Big|\;.
$$
Note that if $N$ is large enough, the left term vanishes if $|z|<
\frac{\lambda_\psi}{2}\;\operatorname{Sys}_{\vec\Lambda}$.

Let $f\in C^1_{\rm c}(\CC^+)$ with support in $B(0,A)$. By integration
on annuli and Equation \eqref{eq:gaussRquat} with $k=3$, we have
\begin{align*}
\Big|\int_{\CC^+} f\;(\theta_\infty-\theta_N)\;d\Leb_\CC\,\Big|&=
\bigO\Big(\frac{\|f\|_\infty}
         {(\lambda_\psi\,\operatorname{Sys}_{\vec\Lambda})^4}
\sum_{p\in\vec\Lambda\,:\, |p|\leq \frac{2A}{\lambda_\psi}} |p|^2\;2\pi\;
\big|\;\lambda_\psi |p|-\frac{\psi(N)|p|}{N}\;\big|\; \Big)
\\&=\bigO\Big(\frac{A^5\;\|f\|_\infty}{\lambda_\psi^9\,
  \operatorname{Sys}_{\vec\Lambda}^4\,\covol_{\vec\Lambda}}
\;\big|\,\lambda_\psi -\frac{\psi(N)}{N}\,\big|\, \Big)\;.
\end{align*}
Hence by Equation \eqref{eq:equation5dePCPC}, we have
{\small \begin{align}
&\frac{1}{\psi(N)^2}\mu_N^+(f)= \frac{1}{\covol_{\vec\Lambda}}
\;\int_{z\in\CC^+} f(z)\;\theta_\infty(z)\;d\Leb_\CC(z)+
\bigO\Big(\frac{A^5\;\|f\|_\infty\,\big|\,\lambda_\psi -
  \frac{\psi(N)}{N}\,\big|}{\lambda_\psi^9\,
\operatorname{Sys}_{\vec\Lambda}^4\,\covol_{\vec\Lambda}^2} \;
 \Big)\nonumber
\\&\;\;+\bigO\Big(\;\frac{A^4\,\diam_{\vec\Lambda}\,\|df\|_\infty\,N^4}
  {\covol_{\vec\Lambda}^2\,\operatorname{Sys}_{\vec\Lambda}\,\psi(N)^5}
+\frac{A^2(\diam_{\vec\Lambda} +\frac{AN}{\psi(N)})\,
  \|f\|_\infty\,N^3}{\covol_{\vec\Lambda}^2\,\psi(N)^4}\;\Big)
\nonumber\\&\;\;=\frac{1}{\covol_{\vec\Lambda}}
\;\int_{z\in\CC^+} f(z)\;\theta_\infty(z)\;d\Leb_\CC(z)
\nonumber\\&\;\;+
\bigO\Big(\;\frac{A^5\;\|f\|_\infty\,\big|\,\lambda_\psi -
  \frac{\psi(N)}{N}\,\big|}{\lambda_\psi^9\,
  \operatorname{Sys}_{\vec\Lambda}^4\,\covol_{\vec\Lambda}^2}+
\frac{A^4\,\diam_{\vec\Lambda}\,\|df\|_\infty}
     {\lambda_\psi^4\covol_{\vec\Lambda}^2\,
       \operatorname{Sys}_{\vec\Lambda}\,\psi(N)}
+\frac{A^2(\diam_{\vec\Lambda} +\frac{A}{\lambda_\psi})\,
\|f\|_\infty}{\lambda_\psi^3\covol_{\vec\Lambda}^2\,\psi(N)}\;\Big)\;.
\label{eq:equationsemifincas3}
\end{align}}

If $\psi'(N)=\psi(N)^2$ as assumed in the third case of Equation
\eqref{eq:troiscas}, it follows from Formula
\eqref{eq:convergmupluscas3} and Lemma \ref{lem:relatRetmu} by
symmetry that
\[
\frac{1}{\psi'(N)}\;\R^{\L_\Lambda,\psi}_N\;\;\weakstar\;\;
\frac{1}{\covol_{\vec\Lambda}}\;\theta_\infty\;\Leb_{\CC}\;.
\]
This proves Formula \eqref{eq:logpaircorrelscal} when
$\lambda_\psi\neq 0,\infty$, with a convergence which is uniform on
every compact subset of $\Lambda$ in $\operatorname{Grid}_2$, as well
as the case $\alpha=1$ in Theorem \ref{theo:intro1} (since if
$\psi(N)=N$, then $\lambda_\psi=1$ and $\psi'(N)=\psi(N)^2=N^2=
N^{4-2\alpha}$). Furthermore, for every $f\in C^1_{\rm c} (\CC^+)$
with support contained in $B(0,A)$, as $N\ra+\infty$ and uniformly on
$\Lambda$ varying in a compact subset of $\operatorname{Grid}_2$,
using Equation \eqref{eq:equationsemifincas3} and the error term in
Lemma \ref{lem:relatRetmu} with the fact that that
$\operatorname{Sys}_{\vec\Lambda}\leq \diam_{\vec\Lambda}$, we have
\begin{align*}
 & \frac{1}{\psi'(N)}\;(\R^{\L_\Lambda,\psi}_N)_{\mid E^+_N}(f_N)=
  \frac{1}{\psi(N)^2}\;\mu_N^+(f)+
  \bigO\Big(
  \frac{A^4\;\|df\|_\infty\;N^4}{\covol_{\vec\Lambda}^2\,\psi(N)^5}\Big)
  \\ &=
  \int_{z\in\CC^+} f(z)\;\frac{\theta_\infty(z)}{\covol_{\vec\Lambda}}\;
  d\Leb_\CC(z)+\\
 & \hspace{1cm}\bigO\Big(\;\frac{A^5\;\|f\|_\infty\,\big|\,\lambda_\psi -
  \frac{\psi(N)}{N}\,\big|}{\lambda_\psi^9\,
  \operatorname{Sys}_{\vec\Lambda}^4\,\covol_{\vec\Lambda}^2} 
  +
\frac{A^4\,\diam_{\vec\Lambda}\,\|df\|_\infty}
     {\lambda_\psi^4\covol_{\vec\Lambda}^2\,
       \operatorname{Sys}_{\vec\Lambda}\,\psi(N)}
+\frac{A^2(\diam_{\vec\Lambda} +\frac{A}{\lambda_\psi})\,
\|f\|_\infty}{\lambda_\psi^3\covol_{\vec\Lambda}^2\,\psi(N)}\;\Big)\;.
\end{align*}
By symmetry, this concludes the proof of Theorem
\ref{theo:logpaircorrelscal}.
\cqfd

\medskip
Let us give a numerical illustration of Theorem
\ref{theo:logpaircorrelscal} when $\Lambda=\vec\Lambda=\ZZ[i]$ and
$\psi(N)=N$.  The following figure shows the points $60\log m-60\log
n$ contained in the ball of radius $5$ centered at $0$ for $(m,n)\in
I_{60}$.
\begin{center}
\includegraphics[width=11cm]{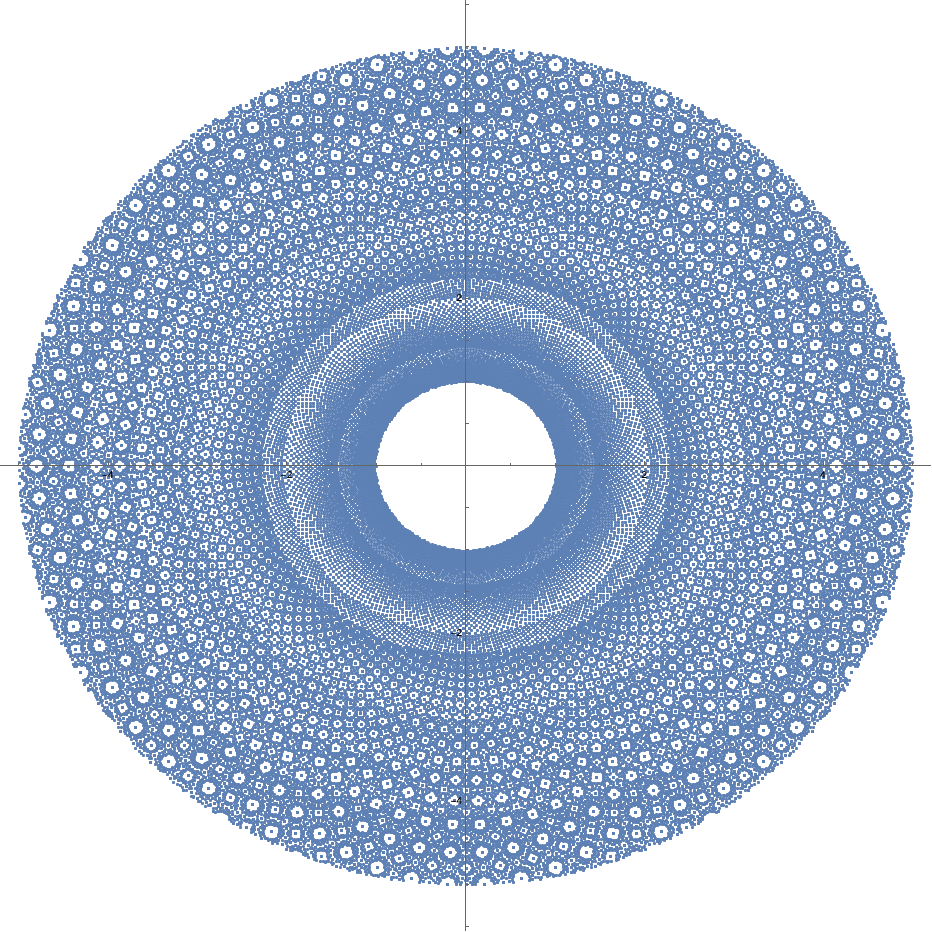}
\end{center}

The second figure shows an approximation (given by Mathematica and its
smoothing process) of the pair correlation function
$g_{\L_\Lambda,\psi}$ computed using the empirical measure
$\frac{1}{60^2}\R^{\L_\Lambda,\psi}_{60}$ in the ball of center $0$
and radius $5$.  We refer to the first picture in the introduction for
the graph of the pair correlation function $g_{\L_\Lambda,\psi}$.

\begin{center}
  \includegraphics[trim = 0mm 40mm 0mm 40mm,clip,width=13cm]
                  {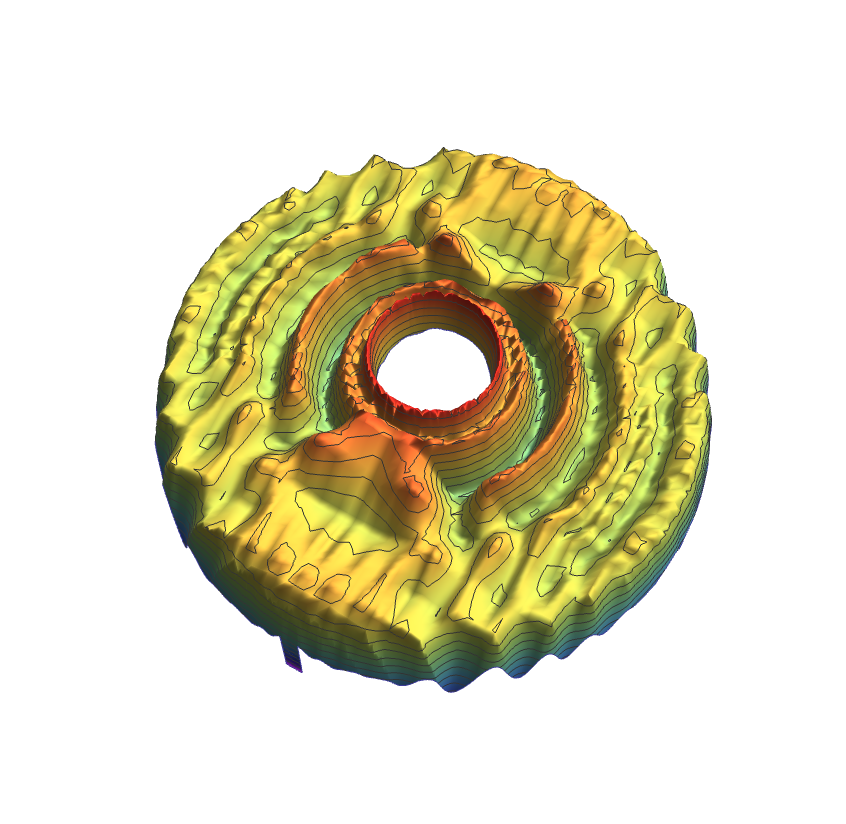}
\end{center}

The figure below gives on the left the graph of the pair correlation
function $g_{\L_\Lambda,\,\psi}$ of the $\ZZ$-lattice $\Lambda=
\vec\Lambda= \ZZ[\frac{1+i\sqrt 3}2]$ of the Eisenstein integers at
the linear scaling $\psi: N\mapsto N$ in the ball of center $0$ and
radius $5$. The blue lines on the bounding box represent the limit
$\frac{2}{\pi\covol^2_{\vec\Lambda}} =\frac{2\pi}{3}$ at $+\infty$ of
$g_{\L_\Lambda,\,\psi}$, given by Equation \eqref{eq:gaussRquat} with
$k=2$.  On the right, we have the approximation of the pair
correlation function computed with the empirical measure
$\frac{1}{60^2}\R^{\L_\Lambda,\psi}_{60}$

\begin{center}
\includegraphics[trim = 10mm 20mm 16mm 30mm,clip,width=7.6cm]{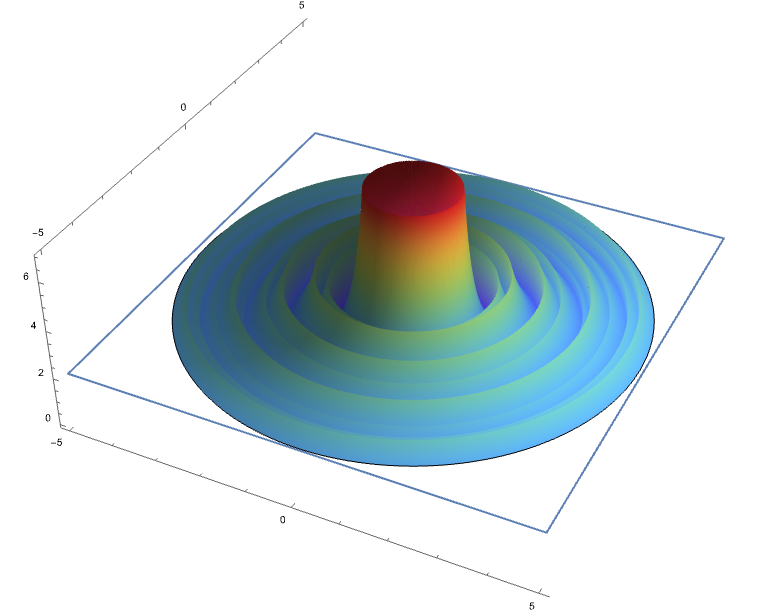}\hspace{-0.2cm}
\includegraphics[trim = 25mm 20mm 30mm 30mm,clip,width=7.4cm]{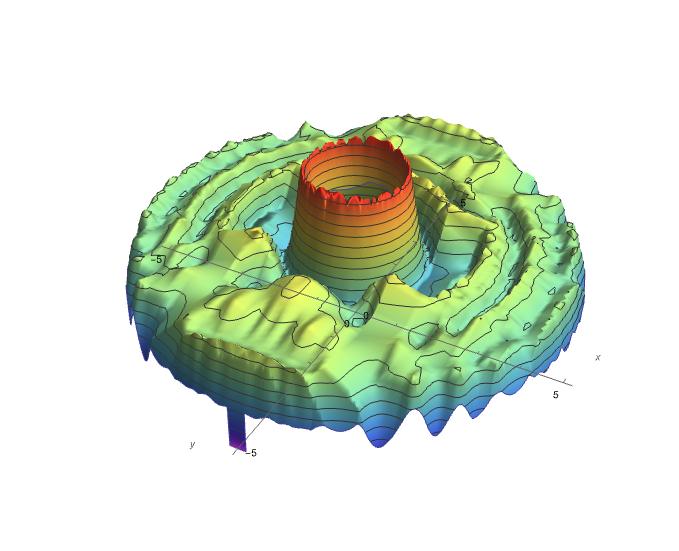}
\end{center}

\section{Mertens and Mirsky formulae for algebraic number fields}
\label{sect:MerMir}

In this short section, we recall the notation and statements of
\cite{ParPau22c} that we will use in Sections
\ref{sect:latlogpaircor3} and \ref{sect:latlogpaircor4}.

Let $K$ be an imaginary quadratic number field (with $D_K$, $\OOO_K$,
$\zeta_K$, $\I^+_K$, $\Nr$ the notation introduced before Corollary
\ref{coro:mmm}). We assume in Sections \ref{sect:MerMir},
\ref{sect:latlogpaircor3} and \ref{sect:latlogpaircor4} that $\OOO_K$
is principal (or equivalently factorial (UFD)). This implies,, see for
instance \cite{Narkiewicz04}, that $D_K\in\{-4,-8,-3,-7,-11, -19, -43,
-67, -163\}$.  For all $I,J\in \I^+_K$, we write $J\mid I$ if
$I\subset J$, we denote by $(I,J)=I+J$ the greatest common ideal
divisor of $I$ and $J$, and by $IJ$ the product ideal of $I$ and $J$.

We denote by $\varphi_K:\I^+_K\ra \NN$ the {\it Euler function} of $K$,
defined (see for instance \cite[page 13]{Narkiewicz04}) equivalently
by
\[
\forall\;\aaa\in\I_K^+,\;\;\;
\varphi_K  (\aaa)=\card\big((\OOO_K/\aaa)^\times\big)=
\Nr (\aaa)\prod_{\ppp\mid \aaa}\big(1-\frac{1}{\Nr (\ppp)}\big)\;,
\]
where, here and thereafter, $\ppp$ ranges over the prime ideals of
$\OOO_K$. For every $a\in \OOO_K-\{0\}$, we define $\varphi_K(a)=
\varphi_K(a\OOO_K)$.

\medskip
We first give a version in angular sectors of the Mertens formula on
the average of the Euler function that will be needed in the proof of
Theorem \ref{theo:latlogpaircorrelphi}.  For all $z\in\CC^\times$,
$\theta\in\;]0,2\pi]$ and $R\geq 0$, we consider the truncated angular
sector
\begin{equation}\label{eq:defiCzthetaR}
C(z,\theta, R)=\big\{\rho\,e^{it}z:
t\in\;\big]-\frac{\theta}{2},\frac{\theta}{2}\big],\;
0<\rho \leq \frac{R}{|z|}\big\}\;.
\end{equation}
Note that for every $z'\in\CC^\times$, we have
\begin{equation}\label{eq:dilatCzthetaR}
  z'C(z,\theta,R)=C(zz',\theta, R\,|z'|\,)\;.
\end{equation}
It is important that the function $\bigO(\cdot)$ in the following
result is uniform in $\mmm$, $z$ and $\theta$.
 For every $\mmm\in \I_K^+$, let
$$
c_\mmm=\Nr (\mmm) \prod_{\ppp\mid\mmm} (1+\frac{1}{\Nr (\ppp)})\,.
$$

\blemm[A Sectorial Mertens formula] \label{lem:sectcount}  Assume that $K$
is imaginary quadratic  with $\OOO_K$ principal. For all
$\mmm\in \I_K^+$, $z\in\CC^\times$ and $\theta\in\;]0,2\pi]$, as
$x\ra+\infty$, we have
\[
\sum_{a\in\mmm\cap C(z,\theta,x)}\varphi_K(a)=
\frac{\theta}{2\,\sqrt{|D_K|}\,\zeta_K(2)\,c_\mmm}\;x^4+\bigO(x^3)\;.
\]
\elemm
\dem See \cite[Thm.~1.2]{ParPau22c}. 
\cqfd

\medskip
We now give an asymptotic formula for the sum in angular sectors of
the products of shifted Euler functions with congruences, which is
used in the proof of Theorems \ref{theo:latlogpaircorrelphi} and
\ref{theo:logpaircorrelphipsi}. When $K=\QQ$ (the sectorial
restriction is then meaningless), this formula is due to Mirsky
\cite[Thm.~9, Eq.~(30)]{Mirsky49} without congruences, and to Fouvry
\cite[Appendix]{ParPau22a} with congruences.

For all $z\in\CC^\times$, $\theta\in\;]0,2\pi]$,
$k\in\OOO_K$, $\mmm\in \I_K^+$ and $x\geq 1$, let
\begin{equation}\label{eq:defiSzthetakmmmx}
S_{z,\theta,k,\mmm}(x)=
\sum_{a\in\mmm\cap C(z,\theta,x)}\varphi_K(a)\,\varphi_K(a+k)\;.
\end{equation}
Let
\begin{equation}\label{eq:simplifcmk}
c_{\mmm,k}= \frac{1}{\Nr(\mmm)}\prod_{\substack{\ppp\\(\ppp,\mmm)\,\mid\, k\OOO_K}}
\Big(1-\frac{\Nr((\ppp,\mmm))}{\Nr(\ppp)^2}\Big)\;
\prod_{\ppp} \Big(1-\frac{\kappa_{\mmm,k}(\ppp)\;
  \kappa'_k(\ppp)\,\Nr((\ppp,\mmm))}{\Nr(\ppp)^2}\Big)\;,
\end{equation}
where
\begin{equation}\label{eq:defikappa}
\kappa_{\mmm,k}(\ppp)=\begin{cases}
(1-\frac{\Nr((\ppp,\mmm))}{\Nr(\ppp)^2})^{-1} &\rm{if }~~
(\ppp,\mmm)\mid k\OOO_K\\ \; 1
&{\rm otherwise} \end{cases} {\rm ~~~and~~~}
\kappa'_{k}(\ppp)=\begin{cases}
1-\frac{1}{\Nr(\ppp)} & {\rm if}~~\ppp\mid k\OOO_K\\ 1 &{\rm otherwise.}
\end{cases}
\end{equation}
For instance, if $\mmm=\OOO_K$ then by \cite[Eq.~(18)]{ParPau22c}, we
have
\begin{equation}\label{eq:computCOKk}
  c_{\OOO_K,k}=\prod_{\ppp}\big(1-\frac{2}{\Nr(\ppp)^2}\big)
\prod_{\ppp\,\mid\, k\OOO_K}\big(1+\frac{1}{\Nr(\ppp)(\Nr(\ppp)^2-2)}\big)\;.
\end{equation}
Since it will be useful in Section \ref{sect:latlogpaircor4}, by
\cite[Lem.~4.2]{ParPau22c}, we have
\begin{equation}\label{eq:mincmk}
  c'_\mmm=\min_{k\in\OOO_K}c_{\mmm,k}>0\;.
\end{equation}

\btheo [A Sectorial Mirsky Formula]\label{theo:mirsky} Assume that $K$
is imaginary quadratic with $\OOO_K$ principal. There exists a
universal constant $C>0$ such that for all $k\in\OOO_K$, $\mmm\in
\I_K^+$, $z\in\CC^\times$, $\theta\in\;]0,2\pi]$ and $x\geq 1$, we
have
$$
\Big|\;S_{z,\theta,k,\mmm}(x)-
\frac{\theta\,c_{\mmm,k}}{3\,\sqrt{|D_K|}}\;x^6\;\Big|\leq
C\big((1+\sqrt{\Nr(k)}\,)\,x^5+\Nr(k)\,x^4\big)\;.
$$
\etheo
\dem See \cite[Thm.~4.1 and Lem.~4.2]{ParPau22c}. 
\cqfd

\section{Pair correlation of integral lattice points with Euler weight
  and no scaling}
\label{sect:latlogpaircor3}

In this section, we fix an imaginary quadratic number field $K$ whose
ring of integers $\OOO_K$ is principal.  We fix a (nonzero integral)
ideal $\Lambda\in \I_K^+$.  Note that $\Lambda=\vec\Lambda$ is a
$\ZZ$-lattice (hence a $\ZZ$-grid) in $\CC$, with $\covol_\Lambda= \Nr
(\Lambda)\,\frac{\sqrt{|D_K|}}{2}$ as seen in Equation
\eqref{eq:calccocoldiamideal}. As in Section
\ref{sect:latlogpaircor1}, we work on the constant cylinder
$E=\CC/(2\pi i\ZZ)$ in this section.

\medskip
Recall that $\L_\Lambda^{\varphi_K}$ is the family defined in Equation
\eqref{eq:defiLLambdavarphi}. For every $N\in\NN-\{0\}$, the (not
normalised, empirical) pair correlation measure of the logarithms of
nonzero elements in $\Lambda$, with trivial scaling function and
multiplicities given by the Euler function, is the measure on $E$ with
finite support defined, with $I_N=I_{N,\Lambda}$ by
\[
\wt\nu_N =\R_N^{\L_\Lambda^{\varphi_K},1}=
\sum_{(m,\,n)\in I_N}\;\varphi_K(m)\;\varphi_K(n)\;\Delta_{\log m-\log n}\;.
\]

\btheo\label{theo:latlogpaircorrelphi} As $N\ra+\infty$, the measures
$\wt\nu_N$ on $E$, renormalized to be probability measures, weak-star
converge to the measure absolutely continuous with respect to the
Lebesgue measure on $E$, with Radon-Nikodym derivative the function
$g_{\L_\Lambda^{\varphi_K},1}:z'\mapsto \frac{1}{\pi}\,
e^{-\,4\,|\Re\;z'|}$, which is independent of $\Lambda$ and $K$:
$$
\frac{\wt\nu_N}{\|\,\wt\nu_N\|}\;\;\;\weakstar\;\;\;
g_{\L_\Lambda^{\varphi_K},1}\;\Leb_E\;.
$$
Furthermore, for all $f\in C^1_{\rm c}(E)$ and $\alpha\in \;]0,\frac{1}{2}[$,
with ${\displaystyle c_\Lambda=\Nr(\Lambda)\prod_{\ppp\,\mid\,\Lambda}
(1+\frac{1}{\Nr(\ppp)})}$, we have
$$
    \frac{\wt\nu_N}{\|\,\wt\nu_N\|}(f)=
    \int_{z'\in E}\;
\frac{1}{\pi}\,e^{-\,4\,|\Re \;z'|}\,f(z')\;d\operatorname{Leb}_E(z') +
\bigO\big(\frac{c_\Lambda\,\|f\|_\infty}{N^{1-2\alpha}}+
\frac{\|e^{z'}df(z')\|_\infty}{N^{\alpha}}\big)\;.
$$
\etheo

This result gives the first assertion of Theorem \ref{theo:intro2} in
the introduction.

\medskip
\dem In this proof, all functions $\bigO(\cdot)$ are absolute, since
there are finitely many $K$. The first assertion of Theorem
\ref{theo:latlogpaircorrelphi} follows from the second one, by the
density of $C^1_c(E)$ in $C^0_c(E)$ for the uniform convergence.

For all $N\in\NN$ and $q\in\Lambda$ with $0<|q|\leq N$, let $J_q$ be
given by the equation on the left in Formula \eqref{eq:defJsubq}. We
now define
$$
\wt\omega_q=\sum_{p\in J_q}\varphi_K(p)\;\Delta_{\frac{p}{q}}\;,
$$
which is a finitely supported measure on the closed unit disc $\DD=B(0,1)$
of $\CC$, and is nonzero since $-q\in J_q$.

\blemm\label{lem:totalmassomegaq}
As $|q|\ra +\infty$, we have ${\displaystyle\|\,\wt\omega_q\|=
\frac{\pi}{\sqrt{|D_K|}\,\zeta_K(2)\,c_\Lambda}\;|q|^4+\bigO(|q|^3)}$.
\elemm

\dem This follows from Lemma \ref{lem:sectcount} applied with
$\mmm=\Lambda$, $z=1$, $\theta=2\pi$ and $x=|q|$, since
$\varphi_K(q)=\bigO(\Nr(q))$ and
\begin{align*}
  \|\,\wt\omega_q\|&=\sum_{p\in\Lambda\;:\;0<|p|\leq |q|,\;p\neq q}\varphi_K(p)
  =\big(\sum_{p\in\Lambda\cap C(1,2\pi,|q|)}\varphi_K(\aaa)\;\big)
  -\varphi_K(q)\;. \;\;\;\Box
\end{align*}

\blemm \label{lem:calcwtomega} For all $f\in C^1(\DD)$ and $\alpha\in
\;]0,\frac{1}{2}\,[\,$, as $|q|\ra+\infty$, we have
$$
\frac{\wt\omega_q}{\|\,\wt\omega_q\|}(f)=
\int_{\DD}\;\frac{2}{\pi}\,|z|^2\,f(z)\;\;d\operatorname{Leb}_\CC(z)
\;+\;\bigO\Big(\frac{c_\Lambda\,\|f\|_\infty}{|q|^{1-2\alpha}}+
\frac{\|df\|_\infty}{|q|^\alpha}\Big)\;.
$$    
\elemm

\dem Note that $c_\Lambda\geq 1$ and let us define
\[
c''_\Lambda= 2\,\sqrt{|D_K|}\,\zeta_K(2)\,c_\Lambda\;.
\]
By Lemma \ref{lem:totalmassomegaq}, as $|q|\ra+\infty$, we have
\begin{equation}\label{eq:invtotalmassomegaq}
  \frac{1}{\|\,\wt\omega_q\|}=\frac{c''_\Lambda}{2\pi\,|q|^4}
  +\bigO\big(\frac{c_\Lambda^2}{|q|^5}\big)\;.
\end{equation}

\medskip
\noindent \begin{minipage}{6.5cm} Let $Q=\lfloor\, |q|^\alpha
  \rfloor\geq 1$, which tends to $+\infty$ as $|q|\ra+\infty$.  For
  all elements $m$ and $n$ in $\{0,\dots, Q-1\}$, let
\begin{align*}
A_{n,m}=\Big\{\rho\,e^{2i\pi\,t}:\;
&\rho\in\;\big]\frac{n}{Q},\frac{n+1}{Q}\big],\;\\&
t\in\;\big]\frac{m}{Q},\frac{m+1}{Q}\big]\Big\}\;,
\end{align*}
so that $\DD-\{0\}$ is the disjoint union of the sets $A_{n,m}$ for
$m,n\in\{0,\dots, Q-1\}$.
\end{minipage}
\begin{minipage}{7.3cm}
\begin{center}
\begin{picture}(0,0)%
\includegraphics{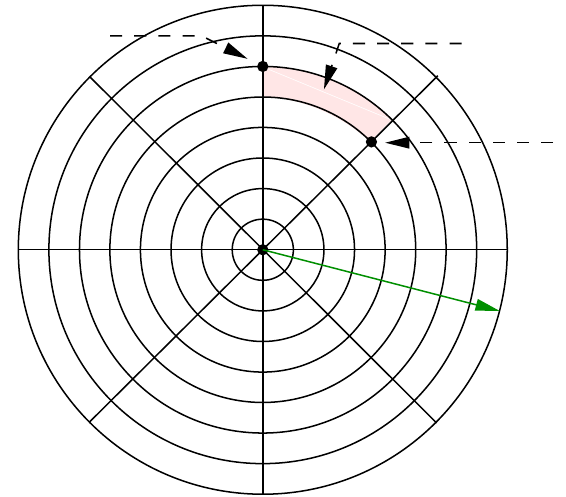}%
\end{picture}%
\setlength{\unitlength}{3812sp}%
\begingroup\makeatletter\ifx\SetFigFont\undefined%
\gdef\SetFigFont#1#2#3#4#5{%
  \reset@font\fontsize{#1}{#2pt}%
  \fontfamily{#3}\fontseries{#4}\fontshape{#5}%
  \selectfont}%
\fi\endgroup%
\begin{picture}(2802,2454)(6196,-2863)
\put(8528,-649){\makebox(0,0)[lb]{\smash{{\SetFigFont{9}{10.8}{\rmdefault}{\mddefault}{\updefault}{\color[rgb]{0,0,0}$A_{n,m}$}%
}}}}
\put(8983,-1142){\makebox(0,0)[lb]{\smash{{\SetFigFont{9}{10.8}{\rmdefault}{\mddefault}{\updefault}{\color[rgb]{0,0,0}$\frac{n}{Q}\,e^{2i\pi\frac{m}{Q}}$}%
}}}}
\put(6211,-573){\makebox(0,0)[lb]{\smash{{\SetFigFont{9}{10.8}{\rmdefault}{\mddefault}{\updefault}{\color[rgb]{0,0,0}$\frac{n+1}{Q}\,e^{2i\pi\frac{m+1}{Q}}$}%
}}}}
\put(7464,-1788){\makebox(0,0)[lb]{\smash{{\SetFigFont{9}{10.8}{\rmdefault}{\mddefault}{\updefault}{\color[rgb]{0,0,0}$0$}%
}}}}
\put(8398,-2016){\makebox(0,0)[lb]{\smash{{\SetFigFont{9}{10.8}{\rmdefault}{\mddefault}{\updefault}{\color[rgb]{0,.56,0}$1$}%
}}}}
\end{picture}%

\end{center}
\end{minipage}

\medskip
\noindent With the notation of Equation \eqref{eq:defiCzthetaR}, we
have
\begin{equation}\label{eq:relatAnmCzthetaR}
  A_{n,m}=C\big(e^{2i\pi\frac{m+1/2}{Q}},\frac{2\pi}{Q},
  \frac{n+1}{Q}\big) - C\big(e^{2i\pi\frac{m+1/2}{Q}},\frac{2\pi}{Q},
  \frac{n}{Q}\big)\;.
\end{equation}
Note that since $n+1\leq Q$, as $Q$ tends to $+\infty$, we have
\[
\diam(\,\overline{A_{n,m}}\,)\leq
\Big|\,\frac{n+1}{Q}e^{2i\pi\,\frac{m+1}{Q}}-
\frac{n+1}{Q}e^{2i\pi\,\frac{m}{Q}}\,\Big| +
\Big|\,\frac{n+1}{Q}e^{2i\pi\,\frac{m}{Q}}-
\frac{n}{Q}e^{2i\pi\,\frac{m}{Q}}\,\Big|=\bigO\big(\frac{1}{Q}\big)\;.
\]

Hence for every $z\in A_{n,m}$, we have by the mean value theorem
\begin{equation}\label{eq:meanvalue1}
f(z)=f\big(\frac{n}{Q}\,e^{2i\pi \frac{m}{Q}}\big)+
\bigO\big(\frac{\|df\|_\infty}{Q}\big)\;.
\end{equation}
Since
\[\int_{A_{n,m}} |z|^2\;d\operatorname{Leb}_\CC(z)
  =\int_{2\pi\frac{m}{Q}}^{2\pi\frac{m+1}{Q}}
  \int_{\frac{n}{Q}}^{\frac{n+1}{Q}} \rho^3\;d\rho \;d\theta=
  \frac{2\pi}{Q}\frac{(n+1)^4-n^4}{4\,Q^4} =\bigO
  \big(\frac{1}{Q^2}\big)\,,
\]
we have therefore
\begin{align}
  \int_{A_{n,m}}\frac{1}{\pi}\,|z|^2\,f(z)\;d\operatorname{Leb}_\CC(z)
  &=\big(f\big(\frac{n}{Q}\,e^{2i\pi \frac{m}{Q}}\big)+
\bigO\big(\frac{\|df\|_\infty}{Q}\big)\big)\int_{A_{n,m}}
\frac{1}{\pi}\,|z|^2\;d\operatorname{Leb}_\CC(z) \nonumber 
\\ & =\frac{1}{Q^2}\Big(\frac{(n+1)^4-n^4}{2\,Q^3}
f\big(\frac{n}{Q}\,e^{2i\pi \frac{m}{Q}}\big) +
\bigO\big(\frac{\|df\|_\infty}{Q}\big)\Big)\;.
\label{eq:controlint}
\end{align}
By Equations \eqref{eq:relatAnmCzthetaR} and \eqref{eq:dilatCzthetaR},
we have $qA_{n,m}=C\big(q\,e^{2i\pi\frac{m+1/2}{Q}},\frac{2\pi}{Q},
\frac{(n+1)|q|}{Q}\big) -
C\big(q\,e^{2i\pi\frac{m+1/2}{Q}},\frac{2\pi}{Q},
\frac{n|q|}{Q}\big)$. By Equations \eqref{eq:meanvalue1},
\eqref{eq:invtotalmassomegaq}, applying twice Lemma
\ref{lem:sectcount} with $\mmm=\Lambda$, $\theta= \frac{2\pi}{Q}$ and
$x=|q|\,\frac{n+1}{Q},\;|q|\,\frac{n}{Q}$, and using the fact that
$\frac{|q|}{Q}$ tends to $+\infty$ as $|q|\ra+\infty$ since
$\alpha<1$, we have, as $|q|\ra+\infty$,
\begin{align}
 &\sum_{p\,\in\,qA_{n,m}\,\cap J_q} f\big(\frac{p}{q}\big)\;
 \frac{1}{\|\,\wt\omega_q\|}\;\varphi_K(p) \nonumber \\  =\;\;&
 \Big(f\big(\frac{n}{Q}\,e^{2i\pi \frac{m}{Q}}\big)+\bigO
 \big(\frac{\|df\|_\infty}{Q}\big)\Big)\;\frac{1}{\|\,\wt\omega_q\|}\;
 \sum_{p\,\in\,\Lambda\cap\, qA_{n,m}\;:\;p\neq q}\varphi_K(p) \nonumber \\  =\;\;&
 \Big(f\big(\frac{n}{Q}\,e^{2i\pi \frac{m}{Q}}\big)+
 \bigO\big(\frac{\|df\|_\infty}{Q}\big)\Big)
\Big(\frac{c''_\Lambda}{2\pi\,|q|^4}+\bigO
\big(\frac{c_\Lambda^2}{|q|^5}\big)\Big)\nonumber \\  &\;\;\;\;\times
\Big(\frac{2\pi \,|q|^4}{Q\,c''_\Lambda}\frac{(n+1)^4-n^4}{Q^4}+
\bigO\big(\frac{|q|^3(n+1)^3}{Q^3}\big)\Big)\nonumber
\\ =\;\;& \frac{1}{Q^2}\Big(\frac{(n+1)^4-n^4}{Q^3}\,
f\big(\frac{n}{Q}\,e^{2i\pi \frac{m}{Q}}\big) +
\bigO\big(\frac{\|df\|_\infty \,n^3}{Q^4}\big)
+\bigO\big(\frac{\|f\|_\infty \,c_\Lambda\,n^3}{|q|\,Q}\big)
\Big)\;.\label{eq:controlmoyphi}
\end{align}
Note that $q\DD=B(0,|q|)$. By cutting the sum defining $\wt\omega_q$
and the integral over $\DD$ into $Q^2$ subparts, by using Equations
\eqref{eq:controlint} and \eqref{eq:controlmoyphi}, and since $n\leq
Q\leq q^\alpha$, as $|q|\ra+\infty$, we have
\begin{align*}
  &\Big|\;\frac{\wt\omega_q}{\|\,\wt\omega_q\|}(f)-
  \int_\DD\;\frac{2}{\pi}\,|z|^2\,f(z)\;d\operatorname{Leb}_\CC(z)
  \;\Big| \\=\; & \Big|\;\sum_{n,m=0}^{Q-1} \Big(
  \sum_{p\,\in\,qA_{n,m}\,\cap J_q}
\frac{\varphi_K(p)}{\|\,\wt\omega_q\|}\;f\big(\frac{p}{q}\big)-
\int_{A_{n,m}}\frac{2}{\pi}\,|z|^2\,f(z)\;
d\operatorname{Leb}_\CC(z)\Big)\;\Big|
\\ =\; &\bigO\big(\frac{\|df\|_\infty}{Q}\big)+
\bigO\big(\frac{\|f\|_\infty \,c_\Lambda}{|q|^{1-2\alpha}}\big)\;.
\end{align*}
This proves Lemma \ref{lem:calcwtomega}.
\cqfd

\medskip
For every $N\in\NN-\{0\}$, let us define
\[
\wt\mu^-_N =\sum_{(m,\,n)\in I^-_N}\;\varphi_K(m)\;\varphi_K(n)\;
\Delta_{\frac{m}{n}}\;\;
=\sum_{q\in\Lambda-\{0\}\;:\;|q|\leq N}\;\varphi_K(q)\;\wt\omega_q\;,
\]
which is a finitely supported measure on $\DD$.  By Lemma
\ref{lem:sectcount} and Theorem \ref{theo:mirsky} both with
$\mmm=\Lambda$, $\theta=2\pi$, $x=N$ and the second one with $k=0$,
since $c_\Lambda\geq 1$ and $c_{\Lambda,0}\leq 1$ by Equation
\eqref{eq:simplifcmk}, and since there are finitely many such fields
$K$), its total mass is
\begin{align*}
  \|\wt\mu^-_N\|&=
  \sum_{q\in\Lambda-\{0\}\;:\;|q|\leq N}\;\varphi_K(q)\;\|\wt\omega_q\|
  =\sum_{(m,n)\in I_N^-} \;\varphi_K(m)\;\varphi_K(n)\\ &=
\frac{1}{2}\Big(
\big(\sum_{\substack{p\in\Lambda-\{0\}\\|p|\leq N}} \varphi_K(p)\big)^2-
\sum_{\substack{p\in\Lambda-\{0\}\\|p|\leq N}}\varphi_K(p)^2\Big)
=\frac{2\,\pi^2}{(c''_{\Lambda})^2}\,N^8+\bigO(N^7)\;.
\end{align*}
For every $f\in C^1(\DD)$, by Lemmas \ref{lem:calcwtomega} and
\ref{lem:totalmassomegaq}, again by Lemma \ref{lem:sectcount} with
$\mmm=\Lambda$, $\theta=2\pi$ and $x=N$, we have
\begin{align*}
  \frac{\wt\mu^-_N(f)}{\|\wt\mu^-_N\|}& = \frac{1}{\|\wt\mu^-_N\|}
  \sum_{q\in\Lambda-\{0\}\;:\;|q|\leq N}\;\varphi_K(q)
\;\|\wt\omega_q\|\;\frac{\wt\omega_q(f)}{\|\wt\omega_q\|}\\ &=
\int_\DD\;\frac{2}{\pi}\,|z|^2\,f(z)\;d\operatorname{Leb}_\CC(z)
+\frac{1}{\|\wt\mu^-_N\|}\sum_{\substack{q\in\Lambda-\{0\}\\|q|\leq N}}
\;\varphi_K(q)\;\|\wt\omega_q\|\,\bigO\big(
\frac{\|df\|_\infty}{|q|^\alpha}+
\frac{c_\Lambda\,\|f\|_\infty }{|q|^{1-2\alpha}}\big)
\\ & =\int_\DD\;\frac{2}{\pi}\,|z|^2\,f(z)\;d\operatorname{Leb}_\CC(z)
\\ & \;\;\;\;\;+\bigO\Big(\frac{(c''_\Lambda)^2}{2\pi^2\,N^{8}}
\sum_{\substack{q\in\Lambda-\{0\}\\|q|\leq N}}\;\varphi_K(q)\frac{2\pi}{c''_\Lambda}
\big(N^{3+2\alpha}c_\Lambda\,\|f\|_\infty+N^{4-\alpha}\,\|df\|_\infty\Big) \\ &
=\int_\DD\;\frac{2}{\pi}\,|z|^2\,f(z)\;d\operatorname{Leb}_\CC(z)
+\bigO\big(\frac{c_\Lambda\,\|f\|_\infty}{N^{1-2\alpha}}+
\frac{\|df\|_\infty}{N^{\alpha}}\big)\;.
\end{align*}

For every $N\in\NN-\{0\}$, let us define
\[
\wt\nu^\pm_N =\sum_{(m,\,n)\in I^\pm_N}\;\varphi_K(m)\;\varphi_K(n)\;
\Delta_{\log m- \log n}\;,
\]
which is a measure with finite support on $E^\pm=
(\pm[0,\infty[\,+i\RR) /(2\pi i \ZZ)$, so that $\wt\nu^-_N=
    \log_*\wt\mu_N^-=\wt\nu_N\!\mid_{E^-}$, and $\|\wt\nu_N^-\|
    =\|\wt\mu_N^-\|$. For every $f\in C^1_{\rm c}(E^-)$, the function
    $f\circ\log$ is a $C^1$ function on $\DD$ which vanishes on a
    neighborhood of $0$. By Equation \eqref{eq:poussmeslog}, we have
\begin{align*}
  \frac{\wt\nu^-_N(f)}{\|\wt\nu^-_N\|}& =
  \frac{\wt\mu^-_N(f\circ\log)}{\|\wt\mu^-_N\|} \\ & = \int_\DD\;
  \frac{2}{\pi}\,|z|^2\,f\circ\log(z)\;d\operatorname{Leb}_\CC(z)
  +\bigO\big(\frac{c_\Lambda\,\|f\circ\log\|_\infty}{N^{1-2\alpha}}+
\frac{\|d(f\circ\log)\|_\infty}{N^{\alpha}}\big)\\ & = \int_{E^-}\;
\frac{2}{\pi}\,e^{4\,\Re(z')}\,f(z')\;d\operatorname{Leb}_E(z') +
\bigO\big(\frac{c_\Lambda\,\|f\|_\infty}{N^{1-2\alpha}}+
\frac{\|e^{z'}df(z')\|_\infty}{N^{\alpha}}\big)\;.
\end{align*}
Since $\wt\nu_N=\wt\nu^-_N + \wt\nu^+_N$, since $\wt\nu_N^+=
\sg_*\wt\nu^-_N$ where $\sg:E\mapsto E$ is the map $x'+iy'\mapsto
-x'+iy'$, we have $\|\wt\nu^\pm_N\|=\frac{1}{2}\;\|\wt\nu_N\|$ and the
last claim of Theorem \ref{theo:latlogpaircorrelphi} follows by
symmetry.  \cqfd

\section{Pair correlation of integral lattice points with scaling
  and Euler weight}
\label{sect:latlogpaircor4}

As in Section \ref{sect:latlogpaircor3}, we fix an imaginary
quadratic number field $K$ whose ring of integers $\OOO_K$ is
principal, and a (nonzero integral) ideal $\Lambda=\vec\Lambda\in
\I_K^+$. We also study the pair correlations of the family
$\L_\Lambda^{\varphi_K}$ defined in the introduction, but now with the
linear scaling function $\psi=\id^1:N\mapsto N$. We leave to the
reader the study of a general scaling $\psi$, assumed to converge to
$+\infty$, proving a Poissonian behaviour for sublinear scalings and
total loss of mass behaviour for superlinear scalings. We also leave
to the reader a statement similar to Theorem
\ref{theo:logpaircorrelphipsi}, replacing the above $\ZZ$-lattice
$\Lambda$ by a $\ZZ$-grid $a+\Lambda$ for any $a\in\OOO_K$.

As in Section \ref{sect:latlogpaircor2}, we work on the family of
varying cylinders $(E_N=\CC/ (2\pi i\,N\,\ZZ))_{N\in \NN-\{0\}}$. As
in Section \ref{sect:latlogpaircor2}, for every $f\in C^1_{\rm c}
(\CC)$, for every $N$ large enough such that the support of $f$ is
contained in $\stackrel{\circ}{B}\!\!(0,\pi N)$, we denote by $f_N\in
C^1_{\rm c} (E_N)$ the map which coincides with $f$ on $B(0,\pi N)$
modulo $2\pi i\,N\,\ZZ$ and vanishes elsewhere. For every
$N\in\NN-\{0\}$, we consider the measure on $E_N$ with finite support
defined with $I_N=I_{N,\Lambda}$ by
\[
\wt\R_N=\R_N^{\L_\Lambda^{\varphi_K},\;\id^1}=
\sum_{(m,\,n)\in I_N}\;\varphi_K(m)\;\varphi_K(n)\;\Delta_{N(\log m-\log n)}\;,
\]
which is the (not normalised) empirical pair correlation measure at
time $N$ of the complex logarithms of the elements of $\Lambda$ with
multiplicities given by the Euler function and with linear scaling
$N$.

\btheo\label{theo:logpaircorrelphipsi} As $N\ra+\infty$, the family
$\big(\frac{1}{N^6}\; \wt\R_N\big)_{N\in\NN}$ of measures on $E_N$
converges (for the pointed Hausdorff-Gromov weak-star convergence) to
the measure absolutely continuous with respect to the Lebesgue measure
on $\CC$, with Radon-Nikodym derivative the function
$$
g_{\L_\Lambda^{\varphi_K},\,\id^1}:z\mapsto \frac{1}{|z|^8}
\sum_{k\in\Lambda\;:\;|k|\leq |z|}
\frac{2\,c_{\Lambda,k}}{\sqrt{|D_K|}}\,|k|^6\;,
$$
that is, as $N\ra+\infty$,
$$
\frac{1}{N^6}\;\wt\R_N\;\;\;\weakstar\;\;\;
g_{\L_\Lambda^{\varphi_K},\,\id^1}\;\Leb_\CC\;.
$$ Furthermore, for all $A\geq 1$ and $f\in C^1(\CC)$ with compact
support contained in $B(0,A)$, as $N\ra+\infty$, we have
$$
\frac{1}{N^6}\;\wt\R_N(f_N)=
\int_{z\in\CC}\;f(z)\;g_{\L_\Lambda^{\varphi_K},\,\id^1}(z)\;d\Leb_\CC(z)
\;+\;\bigO\Big(\frac{A^4}{\covol_\Lambda\,c'_{\Lambda}\,N^{1/2}}
\big(\|df\|_\infty+\|f\|_\infty\big)\Big)\;.
$$
\etheo

The above result with $\Lambda=\OOO_K$ gives the second assertion of
Theorem \ref{theo:intro2} in the introduction, using the value of
$c_{\OOO_K,k}$ given in Equation \eqref{eq:computCOKk}.

Note that, as the proof below shows, the total mass of $\wt\R_N$ is
equivalent to $c\,N^8$ as $N\ra+\infty$, for some constant
$c>0$. Hence renormalising $\wt\R_N$ to be a probability measure would
make it converge to the zero measure on $\CC$.

\medskip
\dem We proceed as in the beginning of the proof of Theorem
\ref{theo:logpaircorrelscal} : We only have to prove the second
assertion above; We define $E^\pm_N=(\pm [0,\infty[\,+i\RR)/(2\pi
    i\,N\, \ZZ)$; We only study the convergence of the measures
    $\frac{1}{N^6} \;\wt\R_N$ on the half-cylinder $E^+_N$ to the
    measure $g_{\L_\Lambda^{\varphi_K},\id^1}\Leb_{\CC^+}$ on the
    half-plane $\CC^+= \{z\in\CC: \Re(z)\geq 0\}$ as $N\ra+\infty$;
    And we deduce the global result by the symmetry of
    $g_{\L_\Lambda^{\varphi_K},\id^1}$ under $z\mapsto -z$.
    
For all $N\in\NN-\{0\}$ and $p\in\Lambda-\{0\}$, let $J_{p,\,N}$ be
given par Equation \eqref{eq:defiJsubpN}.  Note that
\begin{equation}\label{eq:encadJpN}
  (\Lambda-\{0\})\cap B(0,N-|p|)\subset J_{p,\,N}\subset
  (\Lambda-\{0\})\cap B(0,N)\;.
\end{equation}
We now define the key auxiliary measure by
\[
\wt\omega_{p,\,N}=
\sum_{q\in J_{p,\,N}}\varphi_K(q)\;\varphi_K(q+p)\;\Delta_{\frac{q}{N\,p}}\;.
\]
Then $\wt\omega_{p,\,N}$ is a measure with finite support on $B(0,
\frac{1}{|p|})-\{0\}$, which is nonzero if $N\geq 2|p|$ (which is the
case if $p$ is bounded and $N\ra+\infty$), and vanishes if $|p|>2N$.
If $N\geq 2|p|$, by Theorem \ref{theo:mirsky} with $\mmm=\Lambda$,
$k=p$ and $\theta=2\pi$, by Formula \eqref{eq:encadJpN}, since $
|p|\geq 1$, and since $c_{\Lambda,p}\leq 1$ (see Equation
\eqref{eq:simplifcmk}), we have
\begin{align}
  \|\,\wt\omega_{p,\,N}\|&=\sum_{q\in J_{p,\,N}}
  \varphi_K(q)\;\varphi_K(q+p)\nonumber\\ &=
  \frac{2\,\pi\,c_{\Lambda,p}}{3\,\sqrt{|D_K|}}\;(N+\bigO(|p|))^6
  +\bigO\big(\,|p|\,(N+\bigO(|p|))^5+|p|^2\,(N+\bigO(|p|))^4\big)
  \nonumber\\&= \frac{2\,\pi\,c_{\Lambda,p}}{3\,\sqrt{|D_K|}}
  \;N^6 +\bigO(\,|p|\,N^5)\;.\label{eq:totmasswtomsubpN}
\end{align}
In particular, if $N\geq 2|p|$, since $c'_{\Lambda}>0$ by Equation
\eqref{eq:mincmk}, we have
\begin{equation}\label{eq:unsurnormomegpN}
\frac{1}{\|\,\wt\omega_{p,\,N}\|}=
\frac{3\,\sqrt{|D_K|}}{2\,\pi\,c_{\Lambda,p}\,N^6}
\Big(1+\bigO\big(\,\frac{|p|}{c'_{\Lambda}\,N}\big)\Big)\;.
\end{equation}

The next result implies that the measures $\wt\omega_{p,N}$, once
normalized to be probability measures, weak-star converge to the
measure $d\mu(z)=\frac{6}{\pi}\;|p|^6\;|z|^4\;
  d\Leb_{B(0,\frac{1}{|p|})}(z)$ on $B(0,\frac{1}{|p|})$ as
  $N\ra+\infty$, uniformly on $p\in\Lambda-\{0\}$ bounded.

\blemm\label{lem:distribomegasubpN} For all $p\in\Lambda-\{0\}$,
$\alpha\in \;]\,0,1\,[\,$ and $f\in C^1_{\rm c}(\CC)$, as
    $N\ra+\infty$, we have
\[
\frac{\wt\omega_{p,N}}{\|\wt\omega_{p,N}\|}(f)=
\int_{B(0,\frac{1}{|p|})}\frac{3}{\pi}\;|p|^6\;|z|^4\;f(z)\;
d\Leb_{\CC}(z)+\bigO\Big(\,\frac{\|df\|_\infty}{N^\alpha\,|p|}+
  \frac{|p|\,\|f\|_\infty}{c'_{\Lambda}\,N^{1-\alpha}}+
  \frac{\|f\|_\infty}{N^\alpha}\,\Big)\;.
\]
\elemm

\dem As in the proof of Lemma \ref{lem:calcwtomega}, we will estimate
the difference of the main terms in the above centered formula by
cutting the sum defining the renormalized measure $\wt\omega_{p,N}$
and by cutting similarly the integral on $B(0,\frac{1}{|p|})$. We
assume, as we may, that $N\geq 2|p|$.

Let $Q=\lfloor\, N^\alpha\rfloor\geq 1$, which tends to $+\infty$ as
$N\ra+\infty$.  For all $m,n\in\{0,\dots, Q-1\}$, let
\begin{equation}\label{eq:defiAnmprim}
A'_{n,m}=\Big\{\rho\,e^{2i\pi\,t}:\;
\rho\in\;\Big]\frac{n}{Q\,|p|},\frac{n+1}{Q\,|p|}\Big],\;
t\in\;\Big]\frac{m}{Q},\frac{m+1}{Q}\Big]\Big\}\;,
\end{equation}
so that $B(0,\frac{1}{|p|})-\{0\}$ is the disjoint union of the sets
$A'_{n,m}$ for $m,n\in\{0,\dots, Q-1\}$. With the notation of Equation
\eqref{eq:defiCzthetaR}, we have
\begin{equation}\label{eq:relatAprimnmCzthetaR}
  A'_{n,m}=C\big(e^{2i\pi\frac{m+1/2}{Q}},\frac{2\,\pi}{Q},
  \frac{n+1}{Q\,|p|}\big) - C\big(e^{2i\pi\frac{m+1/2}{Q}},\frac{2\,\pi}{Q},
  \frac{n}{Q\,|p|}\big)\;.
\end{equation}
Note that $\diam(\,\overline{A'_{n,m}}\,)=\bigO\big(\frac{1}{Q\,|p|}
\big)$.  Hence for every $z\in A'_{n,m}$, we have by the mean value
theorem
\begin{equation}\label{eq:meanvalue2}
f(z)=f\big(\frac{n}{Q\,|p|}\,e^{2i\pi \frac{m}{Q}}\big)+
\bigO\big(\frac{\|df\|_\infty}{Q\,|p|}\big)\;.
\end{equation}

If $|p|\leq N^{1-\alpha}$ (which is the case if $p$ is bounded and
$N\ra+\infty$) and if $n\leq Q-2$, then
\[
N|p|\frac{n+1}{Q\,|p|}\leq N\frac{Q-1}{Q}\leq N-N^{1-\alpha}
\leq N-|p|\;.
\]
Hence for all $m,n\in\{0,\dots, Q-1\}$, by Formula
\eqref{eq:encadJpN}, if $|p|\leq N^{1-\alpha}$ and if $n\neq Q-1$, we
have
\begin{equation}\label{eq:relatJLambda}
  \big(N\,p\,A'_{n,m}\big)\cap J_{p,N}=
  \Lambda \cap \big(N\,p\,A'_{n,m}\big)\;.
\end{equation}
For all $m,n\in\{0,\dots, Q-1\}$, let
\[
S_{n,m}=
\sum_{q\,\in\,(N\,p\,A'_{n,m})\cap J_{p,N}} f\big(\frac{q}{N\,p}\big)
\;\frac{1}{\|\,\wt\omega_{p,N}\|}\;\varphi_K(q)\;\varphi_K(q+p)\;.
\]
If $n\neq Q-1$, by Equations \eqref{eq:meanvalue2} and
\eqref{eq:relatJLambda} for the first equality, and for the second
one, by Equations \eqref{eq:unsurnormomegpN},
\eqref{eq:relatAprimnmCzthetaR} and \eqref{eq:dilatCzthetaR}, by
Theorem \ref{theo:mirsky} applied twice with $\mmm=\Lambda$, $k=p$,
$\theta=\frac{2\,\pi}{Q}$ and $x=\frac{N(n+1)}{Q}, \frac{N\,n}{Q}$, we
have, as $N\ra+\infty$ (so that in particular $N\geq 2|p|$),
\begin{align}
  S_{n,m}&= \Big(f\big(\frac{n}{Q\,|p|}\,e^{2i\pi \frac{m}{Q}}\big)+
\bigO\big(\frac{\|df\|_\infty}{Q\,|p|}\big)\Big)\;
\frac{1}{\|\,\wt\omega_{p,N}\|}\;\sum_{q\,\in\,\Lambda\cap(N\,p\,A'_{n,m})}
\varphi_K(q)\;\varphi_K(q+p)\nonumber\\ &=
  \Big(f\big(\frac{n}{Q\,|p|}\,e^{2i\pi \frac{m}{Q}}\big)+
\bigO\big(\frac{\|df\|_\infty}{Q\,|p|}\big)\Big)
\frac{3\,\sqrt{|D_K|}}{2\,\pi\,c_{\Lambda,p}\,N^6}
\Big(1+\bigO\big(\,\frac{|p|}{c'_{\Lambda}\,N}\big)\Big)
\nonumber\\ & \;\;\;\;\times
\frac{\frac{2\,\pi}{Q}\,c_{\Lambda,p}}{3\,\sqrt{|D_K|}}\;
\Big(\big(\frac{N(n+1)}{Q}\big)^6-\big(\frac{N\,n}{Q}\big)^6
+\bigO\big(\frac{|p|}{c'_{\Lambda}}\,\big(\frac{N\,n}{Q}\big)^5+
\frac{|p|^2}{c'_{\Lambda}}\,\big(\frac{N\,n}{Q}\big)^4\big)\Big)
\nonumber\\ & =\frac{1}{Q^2}\Big(\frac{(n+1)^6-n^6}{Q^5}
f\big(\frac{n}{Q\,|p|}\,e^{2i\pi \frac{m}{Q}}\big)+
\bigO\big(\frac{\|df\|_\infty}{Q\,|p|}+
\frac{Q\,|p|\,\|f\|_\infty}{c'_{\Lambda}\,N}\big)\Big)\;.
\label{eq:controlmoyphipsi}
\end{align}
Note that by Equations \eqref{eq:relatAprimnmCzthetaR},
\eqref{eq:encadJpN} and \eqref{eq:dilatCzthetaR} for the first
inequality, and for the second one, by Equations
\eqref{eq:unsurnormomegpN} and twice \eqref{eq:totmasswtomsubpN}, as
$N\ra+\infty$, we have
\begin{align}
\Big|\sum_{0\leq m\leq Q-1} S_{Q-1,m}\;\Big| & \leq \|f\|_\infty\;
\frac{1}{\|\,\wt\omega_{p,N}\|}\;
\sum_{q\,\in\,\Lambda\cap(B(0,N)-B(0,N-\frac{N}{Q}))}
\varphi_K(q)\;\varphi_K(q+p)\nonumber\\ &=\|f\|_\infty
\frac{3\,\sqrt{|D_K|}}{2\,\pi\,c_{\Lambda,p}\,N^6}
\Big(1+\bigO\big(\,\frac{|p|}{c'_{\Lambda}\,N}\big)\Big)
\nonumber\\ &\;\;\;\;\times
\frac{2\,\pi\,c_{\Lambda,p}}{3\,\sqrt{|D_K|}}\Big(
\;N^6-(N-\frac{N}{Q})^6 +\bigO(\,\frac{|p|}{c'_{\Lambda}}\,N^5)\Big)
=\bigO\big(\frac{\|f\|_\infty}{Q}\big)\;.\label{eq:controlSend}
\end{align}

For all $m,n\in\{0,\dots, Q-1\}$, let
\[
I_{n,m}=\int_{A'_{n,m}} \frac{3}{\pi}\;|p|^6\;|z|^4\;f(z)\; d\Leb_{\CC}(z)\;.
\]
By Equations \eqref{eq:meanvalue2} and \eqref{eq:defiAnmprim}, we have
\begin{align}
I_{n,m}&=\Big(f\big(\frac{n}{Q\,|p|}\,e^{2i\pi \frac{m}{Q}}\big)+
\bigO\big(\frac{\|df\|_\infty}{Q\,|p|}\big)\Big)\;
\int_{\frac{2\pi m}{Q}}^{\frac{2\pi (m+1)}{Q}}
\int_{\frac{n}{Q\,|p|}}^{\frac{n+1}{Q\,|p|}} \frac{3}{\pi}\;|p|^6\;\rho^5\,
d\rho\;d\theta\nonumber\\ & =
\frac{1}{Q^2}\Big(\frac{(n+1)^6-n^6}{Q^5}
f\big(\frac{n}{Q\,|p|}\,e^{2i\pi \frac{m}{Q}}\big)+
\bigO\big(\frac{\|df\|_\infty}{Q\,|p|}\big)\Big)\;.
\label{eq:controlintpsi}
\end{align}
Furthermore,
\begin{equation}\label{eq:controlIend}
  \Big|\sum_{0\leq m\leq Q-1} I_{Q-1,m}\;\Big| \leq \|f\|_\infty
  \int_{0}^{2\pi}\int_{\frac{1}{|p|}-\frac{1}{Q\,|p|}}^{\frac{1}{|p|}}
  \frac{3}{\pi}\;|p|^6\;\rho^5\,d\rho\;d\theta=
  \bigO\big(\frac{\|f\|_\infty}{Q}\big)\;.
\end{equation}
Since ${\displaystyle B(0,\frac{1}{|p|})-\{0\} =
  \bigsqcup_{n,m=0}^{Q-1} A'_{n,m}}$, putting together Equations
\eqref{eq:controlmoyphipsi}, \eqref{eq:controlintpsi},
\eqref{eq:controlSend} and \eqref{eq:controlIend}, and since
$Q=\lfloor\, N^\alpha\rfloor \in [\frac{ N^\alpha}{2}, N^\alpha]$ for
$N$ large enough, we have
\begin{align*}
 &\Big|\;\frac{\wt\omega_{p,N}}{\|\wt\omega_{p,N}\|}(f)-
  \int_{B(0,\frac{1}{|p|})}
  \frac{3}{\pi}\;|p|^6\;|z|^4\;f(z)\; d\Leb_{\CC}(z)\;\Big|\\=\;&
  \Big|\sum_{n,m=0}^{Q-2} (S_{n,m}-I_{n,m}) +
  \sum_{m=0}^{Q-1} S_{Q-1,m}-\sum_{m=0}^{Q-1} I_{Q-1,m}\;\Big|
  \\\leq \;&
  \sum_{n,m=0}^{Q-2} |S_{n,m}-I_{n,m}| +
  \big|\sum_{m=0}^{Q-1} S_{Q-1,m}\big|+\big|\sum_{m=0}^{Q-1} I_{Q-1,m}\;\Big|
  \\=\;&\bigO\big(\frac{\|df\|_\infty}{Q\,|p|}+
  \frac{Q\,|p|\,\|f\|_\infty}{c'_{\Lambda}\,N}\big)+
  \bigO\big(\frac{\|df\|_\infty}{Q\,|p|}\big)+
  \bigO\big(\frac{\|f\|_\infty}{Q}\big)+
  \bigO\big(\frac{\|f\|_\infty}{Q}\big)
  \\=\;&\bigO\big(\frac{\|df\|_\infty}{N^\alpha\,|p|}+
  \frac{|p|\,\|f\|_\infty}{c'_{\Lambda}\,N^{1-\alpha}}+
  \frac{\|f\|_\infty}{N^\alpha}\big)\;.
\end{align*}
This proves Lemma \ref{lem:distribomegasubpN}.
\cqfd

\medskip
Now, let us introduce the finitely supported measure on $\CC-\{0\}$
defined by
\[
\wt\mu^+_N=\sum_{p\in\Lambda-\{0\}}\;\iota_*\wt\omega_{p,N}\;
=\sum_{p,q\in\Lambda-\{0\}\;:\;|q|\leq |q+p|\leq N}
\;\varphi_K(q)\;\varphi_K(q+p)\;\Delta_{N\frac{p}{q}}\;,
\]
where as previously $\iota:z\mapsto \frac{1}{z}$ (noting that the
measure $\wt\omega_{p,\,N}$ vanishes if $|p|> 2N$ and has finite
support contained in $B(0,\frac{1}{|p|})-\{0\}$).

\blemm\label{lem:comparRsubNmusubN} For all $A\geq 1$ and $f\in
C^1(\CC^+)$ with compact support contained in $B(0,A)$, as
$N\ra+\infty$, we have
$$
\big|\;\wt \R_N\!\mid_{E^+_N}\!(f_N)-\wt\mu_N^+(f)\,\big|=
\bigO\Big(\frac{A^4\|df\|_\infty\,N^5}{\covol_\Lambda}\Big)\;.
$$
\elemm

\dem Let us assume that $N> \frac{A}{\pi}$, so that the ball
$B(0,A)$ injects by the canonical projection $\CC\ra E_N=\CC/(2\pi
i\,N\,\ZZ)$. Note that $f_N$ has support in $E_N^+$. Using the change
of variables $(p,q)\mapsto (m=p+q,n=q)$, we have
\begin{align*}
\wt\R_N(f_N)&=\sum_{(m,n)\in I_N^+}
\;\varphi_K(m)\;\varphi_K(n)\; f_N(N\log m- N\log n)\;
\\ &=\sum_{p,q\in\Lambda-\{0\}\;:\;|q|\leq |q+p|\leq N}
\;\varphi_K(q)\;\varphi_K(q+p)\; f_N(N\log (p+q)- N\log q)\;.
\end{align*}
As in the proof of Lemma \ref{lem:relatRetmu} (see Formulas
\eqref{eq:controlabsp} and \eqref{eq:meanvalscal}), if a pair $(p,q)$
occurs in the index of the sum defining either $\wt\R_N(f_N)$ or
$\wt\mu^+_N(f)$ with nonzero corresponding summand, then
$\frac{|p|}{|q|}=\bigO\big(\frac{A}{N}\big)$, $|p|=\bigO(A)$, and
\[
\big|\,f_N(N\log (p+q)- N\log q) -f(N\frac{p}{q})\,\big|=
\bigO\big(\frac{A^2\|df\|_\infty}{N}\big)\;.
\]
Hence, by Equation \eqref{eq:totmasswtomsubpN}, since
$c_{\Lambda,p}\leq 1$ (see Equation \eqref{eq:simplifcmk}) and by
Lemma \ref{lem:liminftytheta} with $k=0$, we have
\begin{align*}
\big|\;\wt \R_N(f_N)-\wt\mu_N^+(f)\,\big|&\leq
\sum_{p\in\Lambda-\{0\}\;:\;|p|=\bigO(A),\;q\in J_{p,N}}
\;\varphi_K(q)\;\varphi_K(q+p)
\bigO\big(\frac{A^2\|df\|_\infty}{N}\big)\;\\ &=
\sum_{p\in\Lambda-\{0\}\;:\;|p|=\bigO(A)}
\bigO(N^6)\bigO\big(\frac{A^2\|df\|_\infty}{N}\big)=
\bigO\Big(\frac{A^4\|df\|_\infty\,N^5}{\covol_\Lambda}\Big)\;
\end{align*}
This proves Lemma \ref{lem:comparRsubNmusubN}.
\cqfd

\blemm\label{lem:distribmusubN} For all $A\geq 1$ and $f\in
C^1(\CC^+)$ with compact support contained in $B(0,A)$,  as
$N\ra+\infty$, we have
$$
\frac{1}{N^6}\;\wt\mu_N^+(f)=
\int_{\CC^+} f(z)\;g_{\L_\Lambda^{\varphi_K},\,\id^1}(z)\;d\Leb_\CC(z)+
\bigO\Big(\frac{A^4}{\covol_\Lambda\,c'_{\Lambda}\,N^{1/2}}
\big(\|df\|_\infty+\|f\|_\infty\big)\Big)\;.
$$
\elemm

\dem Let $A$ and $f$ be as in the statement, let $N$ be large enough,
and let $\alpha\in\;]0,1[\,$. Since the support of $\wt\omega_{p,\,N}$
is contained in $B(0,\frac{1}{|p|})$, the support of
$\iota_*\wt\omega_{p,\,N}$ is contained in $\{z\in\CC:|z|\geq|p| \}$.
Since a nonzero element of $\OOO_K$ has norm, hence absolute value, at
least $1$, the measures $\wt\mu_N^+$ and $g_{\L_\Lambda^{\varphi_K},\,\id^1}
(z) \;d\Leb_\CC(z)$ both vanish on $\stackrel{\circ}{B}\!\!(0,1)$.
Hence we may assume that the support of $f$ is contained in
$\{z\in\CC:|z|\geq 1\}$, so that the support of $f\circ\iota$ is
compact.  Note that $\|f\circ \iota\|_\infty=\|f\|_\infty$ and as
the support of $f$ is contained in $B(0,A)$, that
\[
\|d(f\circ \iota)\|_\infty\leq A^2\|df\|_\infty\;.
\]

By Equation \eqref{eq:totmasswtomsubpN} and by Lemma
\ref{lem:distribomegasubpN}, by Equation \eqref{eq:jaciota}, since
$1\leq |p|=\bigO(A)$ and $c_{\Lambda,p}\leq 1$, as $N\ra+\infty$, we
hence have
\begin{align*}
  \wt\mu^+_N(f)&=\sum_{p\in\Lambda-\{0\}}\;\iota_*\wt\omega_{p,N}\;
  =\sum_{p\in\Lambda-\{0\}}\;\|\wt\omega_{p,\,N}\|\;
  \frac{\wt\omega_{p,N}}{\|\wt\omega_{p,\,N}\|}(f\circ \iota)
  \\ & =\sum_{p\in\Lambda-\{0\}}
  \Big(\frac{2\,\pi\,c_{\Lambda,p}}{3\,\sqrt{|D_K|}}
  \;N^6 +\bigO(\,|p|\,N^5)\Big)\\ & \;\;\;\;\times
  \Big(\int_{B(0,\frac{1}{|p|})}\frac{3}{\pi}\;|p|^6\;|z|^4\;
  f\circ \iota(z)\;d\Leb_{\CC}(z)+\bigO\Big(\,
  \frac{A^2\|df\|_\infty}{N^\alpha\,|p|}+
  \frac{|p|\,\|f\|_\infty}{c'_{\Lambda}\,N^{1-\alpha}}+
  \frac{\|f\|_\infty}{N^\alpha}\,\Big)\Big)\\ & =N^6\Big(
  \sum_{p\in\Lambda}\frac{2\,c_{\Lambda,p}\,|p|^6}{\sqrt{|D_K|}}
  \int_{|z|\geq |p|} \frac{1}{|z|^8}\;f(z)\;d\Leb_{\CC}(z)
  \\ & \;\;\;\;+\bigO\Big(\sum_{p\in\Lambda\;:\;|p|=\bigO(A)}
 \frac{|p|\,\|f\|_\infty}{N}+ \frac{A^2\|df\|_\infty}{N^\alpha}+
  \frac{A\,\|f\|_\infty}{c'_{\Lambda}\,N^{1-\alpha}}+
  \frac{\|f\|_\infty}{N^\alpha}\Big)\Big)\;.
\end{align*}
By Lemma \ref{lem:liminftytheta} with $k=0$, we hence have
\begin{align*}
  \frac{\wt\mu^+_N(f)}{N^6}&=
  \int_{\CC}\frac{1}{|z|^8}\sum_{p\in\Lambda\;:\;|p|\leq |z|}
  \frac{2\,c_{\Lambda,p}\,|p|^6}{\sqrt{|D_K|}}\;f(z)\;d\Leb_{\CC}(z)
  \\ & \;\;\;\;+\bigO\Big(
  \frac{A^2}{\covol_\Lambda}\big(\frac{A^2\,\|df\|_\infty}{N^\alpha}
+\frac{A\,\|f\|_\infty}{c'_{\Lambda}\,N^{1-\alpha}}+
\frac{\|f\|_\infty}{N^\alpha}\big)\Big)\;.
\end{align*}
Taking $\alpha=\frac{1}{2}$, this proves Lemma \ref{lem:distribmusubN}
since $c'_\lambda\leq 1$ and $A\geq 1$. \cqfd

\medskip
Theorem \ref{theo:logpaircorrelphipsi} now follows from Lemmas
\ref{lem:comparRsubNmusubN} and \ref{lem:distribmusubN}, as explained
in the beginning of the proof.  \cqfd

\medskip
The following figure illustrates Theorem
\ref{theo:logpaircorrelphipsi} when $K=\QQ(\frac{1+i\sqrt 3}2)$ and
$\Lambda=\OOO_K=\ZZ[\frac{1+i\sqrt 3}2]$.  It shows an approximation
of the pair correlation function $g_{\L_\Lambda^{\varphi_K},\,\id^1}$
computed using the empirical measure $\frac{1}{50^6}\wt\R_{50}$ in the
ball of radius $5$ centered at the origin, to be compared with the
orange radial profile of $g_{\L_\Lambda^{\varphi_K},\,\id^1}$ in the
second figure of the introduction.

\begin{center}
  \includegraphics[trim = 0mm 40mm 0mm 60mm, clip,width=15cm]
                  {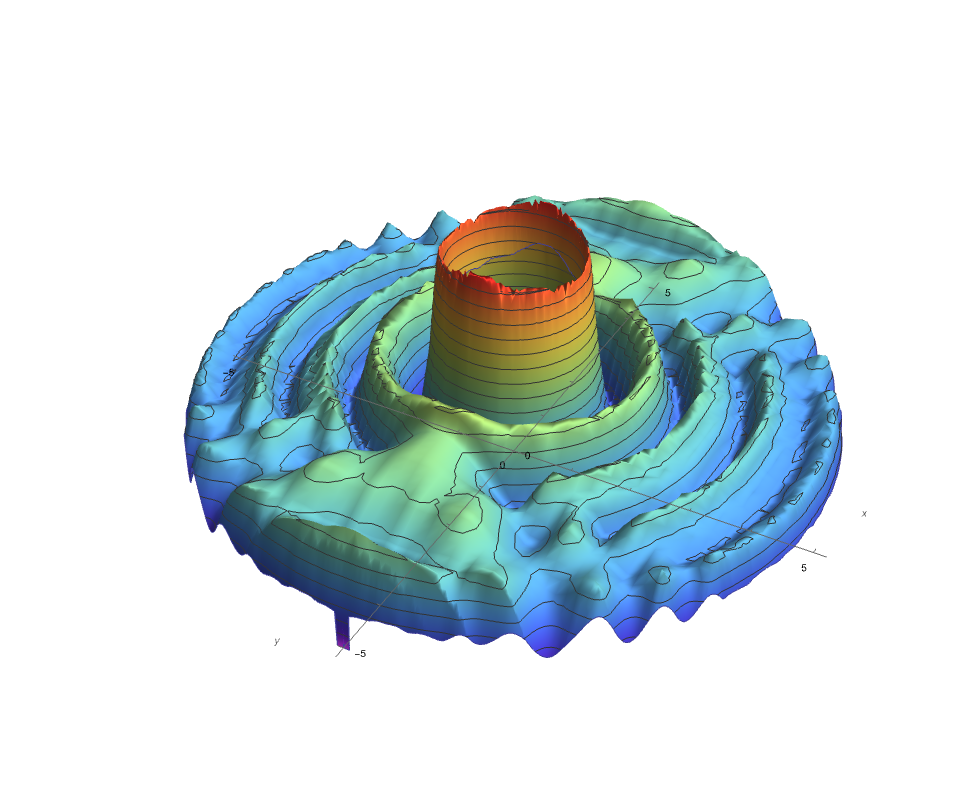}
\end{center}

\medskip

The graph of $g_{\L_\Lambda^{\varphi_K},\,\id^1}$ is bounded by Lemma
\ref{lem:liminftytheta} with $k=6$ since $c_{\Lambda,p}\leq 1$. It is
asymptotic to a horizontal plane at infinity, by the following
result. In its proof, we use the {\it Möbius function}
$\mu_K:\I^+_K\ra \ZZ$ of $K$, defined by
\[
\forall\;\aaa\in\I_K^+,\;\;\;
\mu_K(\aaa)=\begin{cases} 
0 & {\rm if~}\ppp^2\mid \aaa
   {\rm ~for~some~prime~ideal~}\ppp\\
   (-1)^m &{\rm if~} \aaa=\ppp_1\dots\ppp_m
   {\rm ~for~distinct~prime~ideals~}\ppp_1,\dots,\ppp_m
\end{cases}
\]
(in particular $\mu_K(\OOO_K)=1$), For every $a\in \OOO_K-\{0\}$, we
define $\mu_K(a)= \mu_K(a\OOO_K)$.  We have (see for instance
\cite{Shapiro59}) the {\it Möbius inversion formula}: for all
$f,g:\I^+_K\ra \CC$,
\begin{equation}\label{eq:Mobinvform}
  f(\aaa)=\sum_{\bbb\mid \aaa}\;g(\bbb) {\rm~~~if~and~only~if~~~}
  g(\aaa)=\sum_{\bbb\mid \aaa}\;\mu_K(\bbb)f(\aaa\bbb^{-1})\;.
\end{equation}

\bprop\label{prop:horizasymptot} We have
$$\lim_{|z|\ra\infty}\;
g_{\L_\Lambda^{\varphi_K},\,\id^1}(z) =\frac{\pi}{|D_K|}
\;\prod_{\ppp}\big(1-\frac{2}{\Nr(\ppp)^2}\big)
\big(1+\frac{1}{\Nr(\ppp)^2(\Nr(\ppp)^2-2)}\big)\;.
$$
\eprop

\dem Let us consider the multiplicative\footnote{Recall that a
function $f:\I^+_K\ra\CC^\times$ is {\it multiplicative} if
$f(\OOO_K)=1$ and for all coprime integral ideals $\aaa,\bbb$ in
$\I^+_K$, we have $f(\aaa\bbb)=f(\aaa)f(\bbb)$.} function on $\I^+_K$
defined by
\[f:\aaa\mapsto \prod_{\ppp\,\mid\, \aaa} \big(1+\frac{1}{\Nr(\ppp)
  (\Nr(\ppp)^2-2)}\big)
\]
and the constant $C_1=\frac{2\,\pi}{|\OOO_K^\times|\,\sqrt{|D_K|}}
\prod_{\ppp} \big(1+\frac{1}{\Nr(\ppp)^2 (\Nr(\ppp)^2-2)}\big)$.  Let
us prove that uniformly in $x\geq 1$, we have
\begin{equation}\label{eq:asympcube}
  \sum_{\aaa\in \I^+_K\;:\;\Nr(\aaa)\leq x}\Nr(\aaa)^3f(\aaa)=
  \frac{C_1}{4}\;x^4\;+\bigO(x^{7/2})\;.
\end{equation}
Applying this with $x=|z|^2$, by Equation \eqref{eq:gLvarphiid}, since
the map $k\mapsto k\OOO_K$ from $\OOO_K-\{0\}$ onto $\I^+_K$ is
$|\OOO_K^\times|$-to-$1$, this proves Proposition
\ref{prop:horizasymptot}.

Let $$g=f*\mu_K:\aaa\mapsto \sum_{\bbb\mid \aaa}\;\mu_K(\bbb)f(\aaa\bbb^{-1})$$
be the Dirichlet convolution of $f$ with the Möbius function $\mu_K$
of $K$. Then $g$ is multiplicative. For every prime ideal $\ppp$ of
$\OOO_K$, we have
$$
g(\ppp)=f(\ppp)\,\mu_K(\OOO_K)+f(\OOO_K)\,\mu_K(\ppp)=\frac{1}{\Nr(\ppp)
  (\Nr(\ppp)^2-2)}
$$
and $g(\ppp^k)=f(\ppp^k)\,\mu_K(\OOO_K)+f(\ppp^{k-1})\, \mu_K(\ppp)=0$
for every $k\geq 2$.  Therefore, for every $\bbb\in \I^+_K$, we have
$$
g(\bbb)=\mu_K(\bbb)^2\;\prod_{\ppp\,\mid\, \bbb}\;\frac{1}{\Nr(\ppp)
  (\Nr(\ppp)^2-2)}\;.
$$

\blemm\label{lem:controlpointg} For every $\bbb\in \I^+_K$, we have $0\leq
g(\bbb)\leq  \Nr(\bbb)^{-3}\displaystyle{\prod_{\ppp}}
\big(1-\frac{2}{\Nr(\ppp)^2}\big)^{-1}$. In particular, ${\displaystyle
  \sum_{\substack{\bbb\in \I^+_K\\ \Nr(\bbb)\geq x}}
  \frac{g(\bbb)}{\Nr(\bbb)}= \bigO\big(\frac{1}{x}\big)}$.
\elemm

\dem This is immediate if $\mu_K(\bbb)=0$. Otherwise, $\bbb=
\ppp_1\dots \ppp_k$ with $\ppp_1,\dots ,\ppp_k$ pairwise distinct
prime ideals, and
$$
0\leq \Nr(\bbb)^{3}g(\bbb)=\prod_{i=1}^k
\frac{\Nr(\ppp_i)^3}{\Nr(\ppp_i)(\Nr(\ppp_i)^2-2)}=
\prod_{i=1}^k \big(1-\frac{2}{\Nr(\ppp_i)^2}\big)^{-1}\leq
\prod_{\ppp}\big(1-\frac{2}{\Nr(\ppp)^2}\big)^{-1}<+\infty\;.
$$
The last claim follows from the well known error term in the Dedekind
zeta function summation. \cqfd

\medskip By for instance Equation \eqref{eq:gaussRquat} with
$\Lambda=\OOO_K$, $k=0$ and $x=\sqrt{y}$, by Equation
\eqref{eq:calccocoldiamideal} with $\mmm=\OOO_K$, and again since the
map $k\mapsto k\OOO_K$ is $|\OOO_K^\times|$-to-$1$, as $y\ra+\infty$,
we have
\begin{equation}\label{eq:asympcard}
\card\{\aaa\in\I_K^+:\Nr (\aaa)\leq y\} =
\frac{2\pi}{|\OOO_K^\times|\,\sqrt{|D_K|}}\,y+\bigO(y^{\frac{1}{2}})\;.
\end{equation}

Using the Möbius inversion formula \eqref{eq:Mobinvform} for the first
equality, Equation \eqref{eq:asympcard} for the third equality, Lemma
\ref{lem:controlpointg} for the fifth equality and an Eulerian product
(since $g$ is multiplicative and vanishes on ideals divisible by a
nontrivial square) for the sixth equality, with $S(x)= \sum_{\aaa\in
  \I^+_K\;:\;\Nr(\aaa)\leq x} f(\aaa)$, uniformly in $x\geq 1$, we
have
\begin{align*} S(x)&
  = \sum_{\substack{\bbb,\ccc\in \I^+_K\\\Nr(\bbb\ccc)\leq x}} g(\bbb)
  = \sum_{\substack{\bbb\in \I^+_K\\\Nr(\bbb)\leq x}} g(\bbb)
  \sum_{\substack{\ccc\in \I^+_K\\\Nr(\ccc)\leq x/\Nr(\bbb)}} \!1\\ & =
  \sum_{\substack{\bbb\in \I^+_K\\\Nr(\bbb)\leq x}} g(\bbb)
  \Big(\frac{2\,\pi\,x}{|\OOO_K^\times|\,\sqrt{|D_K|}\,\Nr(\bbb)}+
  \bigO\big(\frac{x^{1/2}}{\Nr(\bbb)^{1/2}}\big)\Big)\\ & =
  \frac{2\,\pi\,x}{|\OOO_K^\times|\,\sqrt{|D_K|}}
  \sum_{\bbb\in \I^+_K}\frac{g(\bbb)}{\Nr(\bbb)} +
  \bigO\Big(x\sum_{\substack{\bbb\in \I^+_K\\\Nr(\bbb)\geq x}}
  \frac{g(\bbb)}{\Nr(\bbb)}\Big)+ \bigO\Big(x^{1/2}
  \sum_{\substack{\bbb\in \I^+_K\\\Nr(\bbb)\leq x}}
  \frac{g(\bbb)}{\Nr(\bbb)^{1/2}}\Big)
  \\ & = \frac{2\,\pi\,x}{|\OOO_K^\times|\,\sqrt{|D_K|}}
  \sum_{\bbb\in \I^+_K}\frac{g(\bbb)}{\Nr(\bbb)} +\bigO(x^{1/2})
  \\ & = \frac{2\,\pi\,x}{|\OOO_K^\times|\,\sqrt{|D_K|}}
  \prod_{\ppp}\Big(1+\frac{1}{\Nr(\ppp)^2
  (\Nr(\ppp)^2-2)}\Big)+\bigO(x^{1/2})=
  C_1\;x +\bigO(x^{1/2})\;.
\end{align*}
By summation by parts, we hence have
\begin{align*}
  \sum_{\aaa\in \I^+_K\;:\;\Nr(\aaa)\leq x} \Nr(\aaa)^3f(\aaa)&
  =\int_1^{x}t^3d[S(t)]=
  \big[t^3(C_1\;t+\bigO(t^{1/2}))\big]_1^x-
  3\int_1^xt^2(C_1\;t+\bigO(t^{1/2}))\;dt\\ &=
  \frac{C_1}{4}\;x^4+\bigO(x^{7/2})\;.
\end{align*}
This proves Equation \eqref{eq:asympcube} and concludes the proof of
Proposition \ref{prop:horizasymptot}.
\cqfd

\section{Pair correlations of common perpendiculars in the
  Bianchi manifolds $\PSL(\OOO_K)\bs\htr$}
\label{sect:geometricmotivation}

We again fix an imaginary quadratic number field $K$ whose ring of
integers $\OOO_K$ is principal, and a (nonzero integral) ideal
$\Lambda=\vec\Lambda\in \I_K^+$. In this section, we give a geometric
motivation for the introduction of the Euler function as
multiplicities in the family $\L_\Lambda^{\varphi_K}$ of complex
logarithms of elements of $\Lambda$ defined in Equation
\eqref{eq:defiLLambdavarphi}, and we give a geometric application of
the results in Section \ref{sect:latlogpaircor3}.

We refer to \cite{ParPau17ETDS,BroParPau19} for more information on
the following notions. Let $Y$ be a nonelementary geodesically
complete connected proper locally $\CAT(-1)$ good orbispace, so that
the underlying space of $Y$ is $\Ga\bs\wt Y$ with $\wt Y$ a
geodesically complete proper $\CAT(-1)$ space and $\Ga$ a discrete
group of isometries of $\wt Y$ preserving no point nor pair of points
in $\wt Y\cup\partial_\infty \wt Y$.  Let $D^-$ and $D^+$ be connected
proper nonempty properly immersed locally convex closed subsets of
$Y$, that is, $D^-$ and $D^+$ are locally finite $\Ga$-orbits of
proper nonempty closed convex subsets $\wt D^-$ and $\wt D^+$ of $\wt
Y$.  A {\it common perpendicular} $\alpha$ between $D^-$ and $D^+$ is
the $\Ga$-orbit of the unique shortest arc $\wt \alpha$ between $\wt
D^-$ and $\ga \wt D^+$ for some $\ga\in\Ga$ such that $d(\wt D^-,\ga
\wt D^+)>0$. The {\it multiplicity} $\mult(\alpha)$ of $\alpha$ is the
ratio $A/B$ where

$\bullet$~ $A$ is the number of elements $(\ga_-,\ga_+)\in
(\Ga/\Ga_{D^-}) \times (\Ga/\Ga_{\ga D^+}) $ such that $\wt \alpha$ is
the unique shortest arc between $\ga_-\wt D^-$ and $\ga_+\ga \wt D^+$,
and

$\bullet$~ $B$ is the cardinality of the pointwise stabilizer of
$\wt \alpha$ in $\Ga$.

\noindent
The {\it length} $\lambda(\alpha)$ of the common perpendicular
$\alpha$ is the length of the geodesic segment $\wt \alpha$ in $\wt
Y$.  For every $\ell$ in the set $\operatorname{OL}^\natural(D^-,D^+)$
of lengths of common perpendiculars, the {\it length multiplicity} of
$\ell$ is the sum of the multiplicities of the common perpendiculars
between $D^-$, $D^+$ having the length $\ell$ :
$$
\omega(\ell)=\sum_{\substack{\alpha{\rm~common~perpendicular}\\
     {\rm beween~} D^-{\rm~and~} D^+{\rm~with~}\lambda(\alpha)=\ell}}
\mult(\alpha)\;.
$$
If $\Perp(D^-,D^+)$ is the set of all common perpendiculars from $D^-$
to $D^+$ with multiplicities, then $(\lambda(\alpha))
_{\alpha\in\Perp(D^-,\,D^+)}$ is the {\em marked ortholength spectrum}
from $D^-$ to $D^+$, and the set $\ols(D^-,\,D^+)=
(\operatorname{OL}^\natural(D^-,D^+),\;\omega)$ of the lengths of the
common perpendiculars endowed with the length multiplicity $\omega$ is
the {\em ortholength spectrum} from $D^-$ to $D^+$.

As defined in \cite[\S 6]{ParPau22a}, the {\it pair correlation
  measure of the common perpendiculars from $D^-$ to $D^+$} is the
pair correlation measure of the family
$$\F_{D^-,D^+}=\big((F_N^{D^-,D^+}=
\operatorname{OL}^\natural(D^-,D^+)\cap[0,2\ln N])_{N\in\NN},\;
\omega\big)\,.
$$

Let us specialize these objets as follows. Let
\[
\wt Y=\htr= \big(\,\{(z=x+iy,t)\in\CC\times\RR : t>0\}, ds^2=
\frac{dx^2+dy^2+dt^2}{t^2}\,\big)
\]
be the upper halfspace model of the real hyperbolic $3$-space with
constant curvature $-1$. We identify as usual its space at infinity
$\partial_\infty\htr=(\CC\times\{0\})\cup\{\infty\}$ with
$\PP^1(\CC)=\CC\cup\{\infty\}$.

For every $\bbb\in\I^+_K$, let $\Ga_0[\bbb]$ be {\it Hecke's
  congruence subgroup modulo $\bbb$} of the {\it Bianchi group}
$\PSL_2(\OOO_K)$, which is the preimage of the upper triangular
subgroup of $\PSL_2(\OOO_K/\bbb)$ under the reduction morphism
$\PSL_2(\OOO_K)\ra\PSL_2(\OOO_K/\bbb)$. It acts faithfully on $\htr$
by Poincaré's extension, and is a lattice in the isometry group of
$\htr$. Let $Y^{\bbb}=\Ga_0[\bbb]\bs\htr$, which is a finite (possibly
ramified) cover of the {\it Bianchi orbifold }
$\PSL_2(\OOO_K)\bs\htr$. Note that since $\OOO_K$ is principal, this
Bianchi orbifold has only one cusp (the number of cusps being the
class number of $K$, see for instance \cite{ElsGruMen98}).

Let $\wt D^-=\wt D^+$ be the horoball $\H_\infty=\{(z,t)\in\htr :
t\geq 1\}$ in $\htr$, whose image $D^-=D^+$ in $Y^\bbb$ is a {\it
  Margulis neighbourhood} of a cusp of $Y^\bbb$.  In order to
emphasize the dependence on the ideal $\bbb$, we will use the notation
$\F^\bbb_{D^-,D^+}= \F_{D^-,D^+}$ for the family of lengths of common
perpendiculars between $D^-$ and $D^+$ in $Y^\bbb$.

The following result relates the pair correlation measures of the
common perpendiculars from this Margulis cusp neighbourhood to itself
to the pair correlation measures of the complex logarithms of the
elements of $\Lambda=\bbb$, with multiplicities given by the Euler
function $\varphi_K$. As explained in Remark \ref{rem:mneqn}, in the
following result, we remove from the index set $I_N$ of the
summations defining $\R^{\F^\bbb_{D^-,D^+},\,1}_{N}$ and
$\R^{\L_\bbb^{\varphi_K},\,1}_N$ the assumption that $m\neq n$. Recall
that the map $2\,\Re:E\ra\RR$ is a continuous proper map.

\bprop\label{prop:corrhoromodular}
For every ideal $\bbb\in\I^+_K$,
$$
\R^{\F^\bbb_{D^-,D^+},\,1}_{N}=\frac{4}{|\OOO_K^\times|^2}\;
(2\,\Re)_*\big(\,\R^{\L_\bbb^{\varphi_K},\,1}_N  \,\big)\;.
$$
\eprop

\dem The orbit of $\H_\infty$ under $\Ga_0[\bbb]$ consists, besides
$\H_\infty$ itself, of the Euclidean $3$-balls $\H_{\frac{p}{q}}$ of
Euclidean radius $\frac 1{2|q|^2}$ tangent to the horizontal plane
$\CC$ at the rational elements $\frac pq$ with $p\in\OOO_K$,
$q\in\bbb-\{0\}$ and $(p,q)=1$.

\medskip
\noindent \begin{minipage}{7cm} ~~~ Every common perpendicular between
  $D^-$ and $D^+$ has a vertical representative in $\htr$ which starts
  from a point in $\CC\times \{1\}$ and ends on the boundary of
  $\H_{\frac{p}{q}}$ with $\frac pq$ as above.  Its hyperbolic length
  is $2\ln |q|$.  In particular, the set $\operatorname{OL}^\natural
  (D^-,D^+)$ is equal to
  $$
  \{2\ln |q|=2\,\Re(\log q) : q\in\bbb-\{0\}\}\;.
  $$
\end{minipage}
\begin{minipage}{7.8cm}
\:\:\:\begin{picture}(0,0)%
\includegraphics{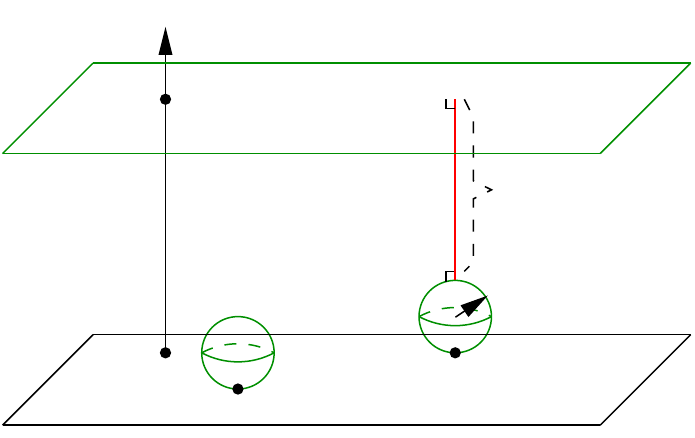}%
\end{picture}%
\setlength{\unitlength}{3812sp}%
\begingroup\makeatletter\ifx\SetFigFont\undefined%
\gdef\SetFigFont#1#2#3#4#5{%
  \reset@font\fontsize{#1}{#2pt}%
  \fontfamily{#3}\fontseries{#4}\fontshape{#5}%
  \selectfont}%
\fi\endgroup%
\begin{picture}(3444,2163)(709,-1381)
\put(1486,-1141){\makebox(0,0)[lb]{\smash{{\SetFigFont{11}{13.2}{\rmdefault}{\mddefault}{\updefault}{\color[rgb]{0,0,0}$0$}%
}}}}
\put(1936,-1231){\makebox(0,0)[lb]{\smash{{\SetFigFont{11}{13.2}{\rmdefault}{\mddefault}{\updefault}{\color[rgb]{0,0,0}$\frac{1}{q}$}%
}}}}
\put(2976,-1141){\makebox(0,0)[lb]{\smash{{\SetFigFont{11}{13.2}{\rmdefault}{\mddefault}{\updefault}{\color[rgb]{0,0,0}$\frac{p}{q}$}%
}}}}
\put(3216,-191){\makebox(0,0)[lb]{\smash{{\SetFigFont{11}{13.2}{\rmdefault}{\mddefault}{\updefault}{\color[rgb]{0,0,0}$2\ln| q|$}%
}}}}
\put(1441,659){\makebox(0,0)[lb]{\smash{{\SetFigFont{11}{13.2}{\rmdefault}{\mddefault}{\updefault}{\color[rgb]{0,0,0}$t$}%
}}}}
\put(1351,299){\makebox(0,0)[lb]{\smash{{\SetFigFont{11}{13.2}{\rmdefault}{\mddefault}{\updefault}{\color[rgb]{0,0,0}$1$}%
}}}}
\put(3836,-1326){\makebox(0,0)[lb]{\smash{{\SetFigFont{11}{13.2}{\rmdefault}{\mddefault}{\updefault}{\color[rgb]{0,0,0}$\CC$}%
}}}}
\put(3151,-691){\makebox(0,0)[lb]{\smash{{\SetFigFont{11}{13.2}{\rmdefault}{\mddefault}{\updefault}{\color[rgb]{0,0,0}$\frac{1}{2|q|^2}$}%
}}}}
\put(2251,579){\makebox(0,0)[lb]{\smash{{\SetFigFont{11}{13.2}{\rmdefault}{\mddefault}{\updefault}{\color[rgb]{0,0,0}$\H_\infty$}%
}}}}
\end{picture}%

\end{minipage}

\medskip
The stabilizer of $\H_\infty$, or equivalently of $\infty$, in
$\Ga_0[\bbb]$ is the upper triangular subgroup of $\Ga_0[\bbb]$, hence
of $\PSL_2(\OOO_K)$. It contains the upper unipotent subgroup
consisting of translations by $\OOO_K$ with finite index, equal to
$\frac{|\OOO_K^\times|}{2}$. Hence given a denominator $q\in\bbb
-\{0\}$, the points at infinity with denominator $q$ of the geodesic
lines containing a lift of a common perpendicular between $D^-$ and
$D^+$ are, modulo translation by $\OOO_K$, exactly the points
$\frac{p}{q}$ where $p$ ranges over a set of representatives of
$(\OOO_K/q\OOO_K)^\times$.

By Equation \eqref{eq:relatlognormRenlog}, the map $z\mapsto \ell=
2\,\Re(z)$ from $L_N=\{\log q: q\in\bbb-\{0\},\; |q|\leq N\}$ to the
set $F_N^{D^-,D^+}$ of lengths of the common perpendiculars between
$D^-$ and $D^+$ with length at most $2\,\ln N$ hence sends the
multiplicity $\frac{2}{|\OOO_K^\times|}\,\varphi_K\circ\exp(\log q)$
of $\log q$ to the multiplicity $\omega(\ell)$ of the corresponding
length of common perpendicular $\ell$.  The claim follows.  \cqfd

\medskip
The following result computes the pair correlation function without
scaling of the lengths of the common perpendiculars from the Margulis
cusp neighbourhood at infinity to itself in the Hecke-Bianchi orbifold
$\Ga_0[\bbb]\bs\hdr$. The maps $\Re:E_N\ra \RR$ for $N\in\NN$ being
not uniformly proper, the case with scalings requires a new analysis,
that we plan to study in another paper.

\bcoro \label{coro:comperphoromodular} For every ideal
$\bbb\in\I^+_K$, as $N\ra+\infty$, the pair correlation measures
$\R^{\F^\bbb_{D^-,D^+},\,1}_N$ on $\RR$, renormalized to be
probability measures, weak-star converge to a measure absolutely
continuous with respect to the Lebesgue measure on $\RR$, with pair
correlation function given by $s\mapsto e^{-\,2|s|}$.  \ecoro

\dem This follows from Theorem \ref{theo:latlogpaircorrelphi} with
$\Lambda=\bbb$ as in the proof of Corollary \ref{coro:applir2d2},
using Proposition \ref{prop:corrhoromodular}. \cqfd

{\small \bibliography{../biblio} }

\bigskip
{\small
\noindent \begin{tabular}{l} 
Department of Mathematics and Statistics, P.O. Box 35\\ 
40014 University of Jyv\"askyl\"a, FINLAND.\\
{\it e-mail: jouni.t.parkkonen@jyu.fi}
\end{tabular}
\medskip

\noindent \begin{tabular}{l}
Laboratoire de mathématique d'Orsay, UMR 8628 CNRS,\\
Universit\'e Paris-Saclay,\\
91405 ORSAY Cedex, FRANCE\\
{\it e-mail: frederic.paulin@universite-paris-saclay.fr}
\end{tabular}}

\end{document}